\documentclass[12pt,twoside]{article}
\usepackage[hmargin=0.8in,vmargin=0.9in]{geometry}
\usepackage{fancyhdr}
\usepackage{graphicx}
\usepackage{subfigure}
\usepackage{amssymb}
\usepackage{amsmath}
\usepackage{amsfonts}
\usepackage{theorem}
\usepackage{mathrsfs}
\usepackage{mathtools}
\usepackage{bm}
\usepackage{color}
\usepackage{setspace}
\usepackage{exscale}
\usepackage{relsize}
\usepackage{epstopdf}
\usepackage{float}
\usepackage{cite}
\usepackage{nicefrac}
\usepackage{setspace}

\usepackage{booktabs,multirow} 
\usepackage{array} 
\usepackage{paralist} 
\usepackage{verbatim} 
\usepackage{subfigure} 

\theoremstyle{plain}			
\newtheorem{thm}{Theorem}[section]

\newtheorem{rmk}[thm]{Remark}
{\theorembodyfont{\rmfamily}}

\usepackage{tikz}
\usetikzlibrary{positioning}
\usepackage{xcolor}

\allowdisplaybreaks[1]

\numberwithin{equation}{section}
\numberwithin{figure}{section}
\numberwithin{table}{section}

\newcommand\eref[1]{(\ref{#1})}

\newcommand*\xbar[1]{%
  \hbox{%
    \vbox{%
      \hrule height 0.5pt 
      \kern0.4ex
      \hbox{%
        \kern-0.05em
        \ensuremath{#1}%
        \kern-0.00em
      }%
    }%
  }%
}

\newcommand{\bmH}{\bm{\mathcal{H}}}

\newcommand{\mK}{\bm{K}}
\newcommand{\mL}{\bm{L}}

\newcommand{\mU}{\bm{U}}

\newcommand{\mV}{\bm{V}}

\newcommand{\dt}{\Delta t}
\newcommand{\dx}{\Delta x}
\newcommand{\dy}{\Delta y}

\newcommand{\hf}{{\frac{1}{2}}}

\newcommand{\jph}{{j+\frac{1}{2}}}
\newcommand{\jmh}{{j-\frac{1}{2}}}
\newcommand{\kph}{{k+\frac{1}{2}}}
\newcommand{\kmh}{{k-\frac{1}{2}}}

\def\softd{{\leavevmode\setbox1=\hbox{d}%
          \hbox to 1.05\wd1{d\kern-0.4ex{\char039}\hss}}}

\title{Low-Dissipation Central-Upwind Schemes for Compressible Multifluids}
\author{Shaoshuai Chu\thanks{Department of Mathematics and Shenzhen International Center for Mathematics, Southern University of Science and
Technology, Shenzhen, 518055, China; {\tt chuss2019@mail.sustech.edu.cn}}, ~Alexander Kurganov\thanks{Department of Mathematics, Shenzhen
International Center for Mathematics, and Guangdong Provincial Key Laboratory of Computational Science and Material Design, Southern
University of Science and Technology, Shenzhen, 518055, China; {\tt alexander@sustech.edu.cn}}, and ~Ruixiao Xin\thanks{Department of
Mathematics, Southern University of Science and Technology, Shenzhen, 518055, China; {\tt xinrx@mail.sustech.edu.cn}}}

\begin{document}

\date{}
\maketitle
\begin{abstract}
We introduce second-order low-dissipation (LD) path-conservative central-upwind (PCCU) schemes for the one- (1-D) and two-dimensional (2-D)
multifluid systems, whose components are assumed to be immiscible and separated by material interfaces. The proposed LD PCCU schemes are
derived within the flux globalization based PCCU framework and they employ the LD central-upwind (LDCU) numerical fluxes. These fluxes have
been recently proposed in [{\sc A. Kurganov and R. Xin}, J. Sci. Comput., 96 (2023), Paper No. 56] for the single-fluid compressible Euler
equations and we rigorously develop their multifluid extensions. In order to achieve higher resolution near the material interfaces, we
track their locations and use an overcompressive SBM limiter in their neighborhoods, while utilizing a dissipative generalized minmod
limiter in the rest of the computational domain. We first develop a second-order finite-volume LD PCCU scheme and then extend it to the
fifth order of accuracy via the finite-difference alternative weighted essentially non-oscillatory (A-WENO) framework. We apply the
developed schemes to a number of 1-D and 2-D numerical examples to demonstrate the performance of the new schemes.
\end{abstract}

\noindent
{\bf Key words:} Low-dissipation central-upwind schemes, path-conservative central-upwind schemes, flux globalization, affine-invariant
WENO-Z interpolation, compressible multifluids.

\medskip
\noindent
{\bf AMS subject classification:} 76M12, 65M08, 76M20, 65M20, 76N30.

\section{Introduction}
In this paper, we focus on the development of highly accurate and conservative finite-volume methods for compressible multifluids, which are
assumed to be immiscible. The studied two-dimensional (2-D) multifluid system reads as
\begin{equation}
\begin{aligned}
&\rho_t+(\rho u)_x+(\rho v)_y=0,\\
&(\rho u)_t+(\rho u^2 +p)_x+(\rho uv)_y=0,\\
&(\rho v)_t+(\rho uv)_x+(\rho v^2+p)_y=0,\\
&E_t+[u(E+p)]_x+[v(E+p)]_y=0.
\end{aligned}
\label{1.1}
\end{equation}
Here, $x$ and $y$ are spatial variables, $t$ is the time, $\rho$ is the density, $u$ and $v$ are the $x$- and $y$-velocities, $p$ is the
pressure, and $E$ is the total energy. The system \eref{1.1} is closed through the following equation of state (EOS) for each of the fluid
components:
\begin{equation}
p=(\gamma-1)\left[E-\frac{\rho}{2}(u^2+v^2)\right]-\gamma\pi_\infty,
\label{1.2}
\end{equation}
where the parameters $\gamma$ and $\pi_\infty$ represent the specific heat ratio and stiffness parameter, respectively. When
$\pi_\infty\equiv0$, the system \eref{1.1}--\eref{1.2} reduces to the ideal gas multicomponent case.

The fluid components can be identified by the variable $\phi$, such as the specific heat ratio $\gamma$ (or a certain function of $\gamma$),
the mass fraction of the fluid component in the fluid mixture, or a level-set function designed to track the interfaces between the fluid
components; see, e.g., \cite{AK_01a,AK_01b,FAMO_99,MOS_92,SA_99,LSDLL17,WS_10,WLS_06,chertock7} and references therein. The state variable
$\phi$ propagates with the fluid velocity and thus satisfies the following advection equation:
\begin{equation}
\phi_t+u\phi_x+v\phi_y=0.
\label{1.3}
\end{equation}

The system \eref{1.1}--\eref{1.3} is a nonlinear hyperbolic system of PDEs and thus its solutions may develop complicated wave structures
including shocks, rarefactions, and contact discontinuities. In the single-fluid regime, that is, when $\gamma\equiv{\rm Const}$ and
$\pi_\infty\equiv{\rm Const}$, the system \eref{1.1}--\eref{1.3} reduces to the Euler equations of gas dynamics, which can be numerically
solved by finite-volume (FV) methods; see, e.g., the monographs \cite{Hesthaven18,LeV02,Tor} and references therein. However, a
straightforward application of single-fluid FV methods to the multifluid system \eref{1.1}--\eref{1.3} may generate spurious pressure and
velocity oscillations, which typically originate near the material interface and then spread all over the computational domain; see, e.g.,
the review paper \cite{AK_01a} and references therein.

In recent years, a variety of FV methods capable of capturing material interfaces in a non-oscillatory manner have been proposed. A fully
conservative approach was first developed in \cite{SA_99a}, where the pressure and velocity remained constant across the material interface.
This approach is robust but may suffer obvious drawbacks when strong shocks pass through the fluid interface. The quasi-conservative
approach was first introduced in \cite{Abgrall96}, where pressure and velocity non-disturbing condition at an isolated material interface
was introduced to analyze and derive the spatial discretization. The resulting schemes reduced the numerical oscillations effectively with
the help of a quasi-conservative discretization. There are also many locally nonconservative approaches designed to prevent
pressure/velocity oscillations by sacrificing the conservation property near material interfaces. The conservation errors in these
approaches are typically small and decay after the mesh is refined. The pressure-based hybrid algorithms \cite{Karni96,CCK_21} are obtained
by switching from the conservative energy equation to the nonconservative pressure one near the interfaces. The ghost-cell methods based on
the single-fluid interpolations leading to two different single-fluid numerical fluxes at the material interfaces (placed at the cell
interfaces at each time step) were introduced in \cite{AK_01b,FAMO_99}. The interface tracking method \cite{chertock7} is based on the
interpolation between the single-fluid data from both sides of the interface and ignoring the ``mixed'' cell data. We note that both the
ghost fluid and interface tracking approaches are very robust in the 1-D case, but their multidimensional extensions are rather cumbersome.
For several high-order WENO schemes for compressible multifluids, we refer the reader to \cite{Don20,HLZT_18,JC_06,Nonomura20,NMTOF_12}.

In this paper, our objective is to develop highly accurate and non-oscillatory numerical schemes for the so-called $\gamma$-based multifluid
systems studied in \cite{Shyue98,CZL20}, which in the 2-D case read as \eref{1.1}--\eref{1.2} together with the equations \eref{1.3} for the
state variables $\Gamma:=1/(\gamma-1)$ and $\Pi:=\gamma\pi_\infty/(\gamma-1)$, which we recast as follows:
\begin{equation}
\Gamma_t+(u\Gamma)_x+(v\Gamma)_y=\Gamma(u_x+v_y),\quad\Pi_t+(u\Pi)_x+(v\Pi)_y=\Pi(u_x+v_y).
\label{1.5}
\end{equation}
The resulting system is nonconservative (in fact, it can be rewritten in the conservative form, but as it was shown in \cite{Shyue98}, a
nonconservative form is preferable for designing an accurate numerical method) and the nonconservative terms on the right-hand side require
a special treatment.

We numerically solve the system \eref{1.1}--\eref{1.2}, \eref{1.5} and its one-dimensional (1-D) version by the Riemann-problem-solver-free
central-upwind (CU) schemes, which were originally introduced in \cite{KLin,KNP,KTcl} for general multidimensional hyperbolic systems of
conservative laws, and then extended to nonconservative hyperbolic systems in \cite{CKM}, where path-conservative CU (PCCU) schemes were
introduced. The PCCU schemes were extended to the flux globalization framework allowing to treat a wider variety of nonconservative systems
in \cite{CKLX_22,CKLZ_23,CKN22a,KLX_21}.

The aforementioned PCCU schemes are based on the CU numerical fluxes from \cite{KNP,KLin}, which have relatively large numerical dissipation
preventing high resolution of contact waves/material interfaces. In the recent work \cite{KX_22}, we have introduced a new way of reducing
the numerical dissipation present in the CU schemes and introduced the low-dissipation CU (LDCU) schemes. In these schemes, the dissipation
is reduced at the projection step, performed after the numerical solution is evolved to the new time level. The novel projection is based on
a subcell resolution technique, which introduces several degrees of freedom that can be utilized to better approximate contact waves and
material interfaces. In \cite{KX_22}, we have designed the LDCU schemes for both the 1-D and 2-D single-fluid compressible Euler equations
and in this paper, we extend the LDCU schemes to the $\gamma$-based multifluid models. The extension is carried out in the flux
globalization PCCU framework and results in new flux globalization based LD PCCU schemes.

We also extend the proposed LD PCCU schemes to the fifth order of accuracy using the framework of the finite-difference alternative WENO
(A-WENO) schemes developed in \cite{JSZ,Liu17,Liu16,WDGK_20,WLGD_18,WDKL,CKX22}. Our new fifth-order schemes are based on the LD PCCU
numerical fluxes, a new, more efficient way to approximate the high-order A-WENO correction terms (see \cite{CKX23}), and the recently
proposed fifth-order affine-invariant WENO-Z (Ai-WENO-Z) interpolation \cite{DLWW22,LLWDW23,WD22} applied to the local characteristic
variables with the help of the local characteristic decomposition (LCD).

This paper is organized as follows. \S\ref{sec2} is devoted to the 1-D LD PCCU scheme. In \S\ref{sec2.1}, we give an overview of the flux
globalization based PCCU schemes and develop such scheme for the $\gamma$-based multifluid model. The LD PCCU scheme is derived in
\S\ref{sec2.2}, where we prove that the new scheme preserves constant velocity and pressure across isolated material interfaces (this
ensures lack of pressure/velocity-based oscillations). In \S\ref{sec2.3}, the second-order LD PCCU scheme is extended to the fifth-order
flux globalization based LD Ai-WENO PCCU scheme. In \S\ref{sec3}, we present the 2-D extensions of the new LD PCCU schemes. In \S\ref{sec4},
we test the proposed LD PCCU schemes together with their fifth-order versions on a number of 1-D and 2-D numerical examples. Finally, in
\S\ref{sec5}, we give some concluding remarks and comments.

\section{One-Dimensional Algorithms}\label{sec2}
In this section, we present the new 1-D flux globalization based LD PCCU schemes for the 1-D compressible multifluids.

\subsection{Flux Globalization Based Path-Conservative Central-Upwind Schemes}\label{sec2.1}
We begin with a brief overview of the flux globalization based PCCU scheme, which was introduced in \cite{KLX_21} for the general
nonconservative system
\begin{equation*}
\bm U_t+\bm F(\bm U)_x=B(\bm U)\bm U_x,
\end{equation*}
which can be rewritten in the following quasi-conservative form:
\begin{equation}
\bm U_t + \bm K(\bm U)_x=0,\quad\bm K(\bm U)=\bm F(\bm U)-\bm R(\bm U),
\label{2.1.2}
\end{equation}
where $\bm U(\bm x,t)\in\mathbb R^d$ is the vector of unknowns, $\bm F:\mathbb R^d\to\mathbb R^d$ is a flux, $B\in\mathbb R^{d\times d}$,
\begin{equation*}
\bm R(\bm U):=\int\limits^x_{\widehat x}B(\bm U)\bm U_\xi(\xi,t)\,{\rm d}\xi,
\end{equation*}
and $\widehat x$ is an arbitrary number.

We first introduce a uniform mesh consisting of the finite-volume cells $C_j:=[x_\jmh,x_\jph]$ of size $x_\jph-x_\jmh\equiv\dx$ centered at
$x_j=(x_\jmh+x_\jph)/2$, $j=1,\ldots,N$ and set $\widehat x=x_\hf$. We assume that at a certain level $t$, an approximate solution, realized
in terms of the cell averages
\begin{equation*}
\xbar{\bm U}_j(t):\approx\frac{1}{\dx}\int\limits_{C_j}\bm U(x,t)\,{\rm d}x,
\end{equation*}
is available. The cell averages $\xbar{\bm U}_j$ are then evolved in time by solving the following system of ODEs (see \cite{KLX_21}):
\begin{equation}
\frac{{\rm d}}{{\rm d}t}\,\xbar{\bm U}_j=-\frac{\bm{{\cal K}}_\jph-\bm{{\cal K}}_\jmh}{\dx},
\label{2.1.4}
\end{equation}
where $\bm{{\cal K}}_\jph$ are the CU numerical fluxes
\begin{equation}
\bm{{\cal K}}_\jph=\frac{a^+_\jph\bm K^-_\jph-a^-_\jph\bm K^+_\jph}{a^+_\jph-a^-_\jph}+\frac{a^+_\jph a^-_\jph}{a^+_\jph-a^-_\jph}
\left(\bm U^+_\jph-\bm U^-_\jph\right).
\label{2.1.5}
\end{equation}
Notice that all of the indexed quantities in \eref{2.1.4} and \eref{2.1.5} as well as other indexed quantities introduced below depend on
$t$, but from now on we will omit this dependence for the sake of brevity. In \eref{2.1.5}, $\bm U^\pm_\jph$ are the left/right-sided point
values of $\bm U$ at the cell interfaces $x_\jph$, which are computed using a piecewise linear reconstruction, which will be discussed in
\S\ref{sec2.1.1}, and $a^\pm_\jph$ are the one-sided local speeds of propagation, which can be estimated using the largest and smallest
eigenvalues of the matrix $\frac{\partial\bm F}{\partial\bm U}(\bm U)-B(\bm U)$. The global fluxes $\bm K^\pm_\jph$ are obtained using the
relation in \eref{2.1.2}, namely, by
\begin{equation}
\bm K^\pm_\jph=\bm F^\pm_\jph-\bm R^\pm_\jph,
\label{2.1.6}
\end{equation}
where $\bm F^\pm_\jph:=\bm F(\bm U^\pm_\jph)$ and the point values of the global variable $\bm R$ are computed as follows. First, we set
$\bm R^-_\hf:=\bm0$ and then evaluate
\begin{equation}
\bm R^+_\hf=\bm B_{\bm\Psi,\hf},
\label{2.1.7}
\end{equation}
and recursively
\begin{equation}
\bm R^-_\jph=\bm R^+_\jmh+\bm B_j,\quad\bm R^+_\jph=\bm R^-_\jph+\bm B_{\bm\Psi,\jph},\quad j=1,\ldots,N.
\label{2.1.9}
\end{equation}
In \eref{2.1.7} and \eref{2.1.9}, $\bm B_j$ and $\bm B_{\bm\Psi,\jph}$ are obtained using a proper quadrature for the integrals in
\begin{equation}
\bm B_j\approx\int\limits_{C_j}B(\bm U)\bm U_x\,{\rm d}x\quad\mbox{and}\quad
\bm B_{\bm\Psi,\jph}=\int\limits^1_0B\big(\bm\Psi_\jph(s)\big)\bm\Psi'_\jph(s)\,{\rm d}s,
\label{2.1.11}
\end{equation}
where, $\bm\Psi_\jph(s):=\bm\Psi\big(s;\bm U^-_\jph,\bm U^+_\jph\big)$ is a sufficiently smooth path connecting the states $\bm U^-_\jph$
and $\bm U^+_\jph$, that is,
\begin{equation*}
\bm\Psi:[0,1]\times\mathbb R^d\times\mathbb R^d\to\mathbb R^d,\quad\bm\Psi(0;\bm U^-_\jph,\bm U^+_\jph)=\bm U^-_\jph,\quad
\bm\Psi(1;\bm U^-_\jph,\bm U^+_\jph)=\bm U^+_\jph.
\end{equation*}

\subsubsection{Application to the Compressible Multifluid System}\label{sec2.1.1}
We apply the flux globalization based PCCU scheme to the 1-D $\gamma$-based multifluid system
\begin{equation}
\begin{aligned}
&\rho_t+(\rho u)_x=0,\\
&(\rho u)_t+(\rho u^2 +p)_x=0,\\
&E_t+[u(E+p)]_x=0,\\
&\Gamma_t+(u\Gamma)_x=\Gamma u_x,\\
&\Pi_t+(u\Pi)_x=\Pi u_x,
\end{aligned}
\label{2.0}
\end{equation}
completed with the following EOS:
\begin{equation}
p=(\gamma-1)\left[E-\frac{\rho}{2}u^2\right]-\gamma\pi_\infty.
\label{2.0.1}
\end{equation}
The system \eref{2.0} can be rewritten in the equivalent quasi-conservative form \eref{2.1.2} with
\begin{equation}
\begin{aligned}
&\bm U:=(\rho,\rho u,E,\Gamma,\Pi)^\top,\quad\bm F(\bm U)=(\rho u,\rho u^2+p,u(E+p),u\Gamma,u\Pi)^\top,\\
&\mbox{and}\quad\bm R(\bm U)=\bigg(0,0,0,\int\limits^x_{\widehat x}\Gamma u_\xi\,{\rm d}\xi,\int\limits^x_{\widehat x}\Pi u_\xi\,{\rm d}\xi
\bigg)^\top.
\end{aligned}
\label{2.1.12}
\end{equation}

We first discuss the reconstruction procedure for recovering the point values $\bm U^\pm_\jph$ out of the cell averages $\xbar{\bm U}_j$.
Since the variables $u$ and $p$ are continuous across material interfaces (contact waves), we reconstruct the primitive variables
$\bm V:=(\rho,u,p,\Gamma,\Pi)^\top$ instead of the conservative ones. To this end, we compute the cell centered values of $u$ and $p$,
\begin{equation}
u_j=\frac{(\xbar{\rho u})_j}{\xbar{\rho}_j},\quad p_j=\frac{1}{\xbar{\Gamma}_j}\left[\xbar E_j-\frac{((\xbar{\rho u})_j)^2}{2\xbar{\rho}_j}-
\xbar{\Pi}_j\right],
\label{2.1.13}
\end{equation}
and then construct the linear pieces
\begin{equation}
\widetilde{\bm V}_j(x)= \bm V_j+(\bm V_x)_j(x-x_j),\quad x\in C_j,
\label{2.1.14}
\end{equation}
where $\bm V_j:=(\xbar\rho_j,u_j,p_j,\xbar\Gamma_j,\xbar\Pi_j)^\top$ and the slopes $(\bm V_x)_j$ are supposed to be computed with the help
of a nonlinear limiter to ensure a non-oscillatory nature of \eref{2.1.14}. In the numerical experiments reported in \S\ref{sec4}, we
implement a simple adaptive limiting strategy and use different limiters near and away from the material interfaces. To this end, we need to
detect the location of the interfaces. In the two-fluid case, this can be done as follows. We first introduce
$\widehat\Gamma:=(\Gamma_{\rm I}+\Gamma_{\rm II})/2$, where $\Gamma_{\rm I}$ and $\Gamma_{\rm II}$ are the values of $\Gamma$ for the first
and second fluid, respectively. We then assume that the interface is located either in cell $C_j$ or $C_{j+1}$ if
\begin{equation}
(\xbar\Gamma_j-\widehat\Gamma)(\xbar\Gamma_{j+1}-\widehat\Gamma)<0.
\label{2.1.15a}
\end{equation}
In these two cells as well as in the neighboring cells $C_{j-1}$ and $C_{j+2}$, we use the overcompressive SBM limiter \cite{LN}:
\begin{equation}
(\bm V_x)_j=\phi^{\rm SBM}_{\theta,\tau}\left(\frac{\,\xbar{\bm V}_{j+1}-\xbar{\bm V}_j}{\,\xbar{\bm V}_j-\xbar{\bm V}_{j-1}}\right)
\frac{\xbar{\bm V}_{j+1}-\xbar{\bm V}_j}{\dx},
\label{2.1.15b}
\end{equation}
where the two-parameter function
\begin{equation}
\phi^{\rm SBM}_{\theta,\tau}(r):=\left\{\begin{aligned}&0&&\mbox{if}~r<0,\\&\min\{r\theta,1+\tau(r-1)\}&&\mbox{if}~0<r\leq1,\\
&r\phi^{\rm SBM}_{\theta,\tau}\Big(\frac{1}{r}\Big)&&\mbox{otherwise},\end{aligned}\right.
\label{2.1.15c}
\end{equation}
is applied in a componentwise manner with $\tau=-0.5$, which belongs to the overcompressive range of values of $\tau$; see \cite{LN}. Away
from the material interfaces, we use a dissipative generalized minmod limiter which is given by the same formulae
\eref{2.1.15b}--\eref{2.1.15c}, but with $\tau=0.5$. In both areas, we use $\theta=1.3$.
\begin{rmk}
Notice that the generalized minmod limiter can be written in a simpler form \cite{LN,NT,Swe} by
\begin{equation*}
(\bm V_x)_j={\rm minmod}\left(\theta\,\frac{\,\xbar{\bm V}_{j+1}-\,\xbar{\bm V}_j}{\dx},\,
\frac{\,\xbar{\bm V}_{j+1}-\,\xbar{\bm V}_{j-1}}{2\dx},\,\theta\,\frac{\,\xbar{\bm V}_j-\,\xbar{\bm U}_{j-1}}{\dx}\right),
\end{equation*}
with the minmod function defined by
\begin{equation*}
{\rm minmod}(c_1,c_2,\ldots)=\left\{\begin{aligned}
&\min(c_1,c_2,\ldots)&&{\rm if}~c_i>0,~\forall i,\\
&\max(c_1,c_2,\ldots)&&{\rm if}~c_i<0,~\forall i,\\
&0&&\mbox{otherwise}.
\end{aligned}\right.
\end{equation*}
\end{rmk}
\begin{rmk}
If the number of fluid components is more than two, detecting interface cells becomes a more complicated task. A reasonable extension of the
strategy used in \eref{2.1.15a} should be developed for the problem at hand.
\end{rmk}

Equipped with $(\bm V_x)_j$, we then use \eref{2.1.14} to obtain
\begin{equation}
\bm V^-_\jph=\lim_{x\to x^-_\jph}\widetilde{\bm V}(x)=\,\xbar{\bm V}_j+\frac{\dx}{2}(\bm V_x)_j,\quad
\bm V^+_\jph=\lim_{x\to x^+_\jph}\widetilde{\bm V}(x)=\,\xbar{\bm V}_{j+1}-\frac{\dx}{2}(\bm V_x)_{j+1},
\label{2.1.17}
\end{equation}
and the corresponding point values of the conservative variables $\bm U$:
\begin{equation}
{\bm U}^\pm_\jph=\left(\rho^\pm_\jph,\rho^\pm_\jph u^\pm_\jph,E^\pm_\jph,\Gamma^\pm_\jph,\Pi^\pm_\jph\right)^\top,\quad
E^\pm_\jph=\Gamma^\pm_\jph p^\pm_\jph+\frac{\rho^\pm_\jph}{2}\big(u^\pm_\jph\big)^2+\Pi^\pm_\jph.
\label{2.1.18}
\end{equation}
We then compute the point values $\bm K^\pm_\jph$ following \eref{2.1.6}--\eref{2.1.9}. First, we use \eref{2.1.12} and \eref{2.1.18} to
obtain
\begin{equation}
\bm F^\pm_\jph=\left(\rho^\pm_\jph u^\pm_\jph,\rho^\pm_\jph\big(u^\pm_\jph\big)^2+p^\pm_\jph,u^\pm_\jph\big(E^\pm_\jph+p^\pm_\jph\big),
\Gamma^\pm_\jph u^\pm_\jph,\Pi^\pm_\jph u^\pm_\jph\right)^\top,
\label{2.1.19a}
\end{equation}
and then evaluate $\bm B_j$ in \eref{2.1.11} by substituting there the piecewise linear reconstructions \eref{2.1.14} of $u$, $\Gamma$, and
$\Pi$, which results in
\begin{equation}
\bm B_j=\bigg(0,0,0,\frac{\Gamma^-_\jph+\Gamma^+_\jmh}{2}\big(u^-_\jph-u^+_\jmh\big),
\frac{\Pi^-_\jph+\Pi^+_\jmh}{2}\big(u^-_\jph-u^+_\jmh\big)\bigg)^\top.
\label{2.1.19}
\end{equation}
Next, in order to obtain $\bm B_{\bm\Psi,\jph}$ in \eref{2.1.11}, a proper path connecting the states $(u^-_\jph,\Gamma^-_\jph,\Pi^-_\jph)$
and $(u^+_\jph,\Gamma^+_\jph,\Pi^+_\jph)$ needs to be used, for instance, a simple linear path:
\begin{equation}
\begin{aligned}
&\bm\Psi^u_\jph(s)=u^-_\jph+s\big(u^+_\jph-u^-_\jph\big),~~\bm\Psi^\Gamma_\jph(s)=\Gamma^-_\jph+s\big(\Gamma^+_\jph-\Gamma^-_\jph\big),\\
&\bm\Psi^\Pi_\jph(s)=\Pi^-_\jph+s\big(\Pi^+_\jph-\Pi^-_\jph\big).
\end{aligned}
\label{2.1.20}
\end{equation}
Substituting \eref{2.1.20} into \eref{2.1.11} then results in
\begin{equation}
\bm B_{\bm\Psi,\jph}=\bigg(0,0,0,\frac{\Gamma^+_\jph+\Gamma^-_\jph}{2}\big(u^+_\jph-u^-_\jph\big),
\frac{\Pi^+_\jph+\Pi^-_\jph}{2}\big(u^+_\jph-u^-_\jph\big)\bigg)^\top.
\label{2.1.21}
\end{equation}

Finally, the one-sided local-speeds of propagation $a^\pm_\jph$ can be estimated by
\begin{equation*}
a^+_\jph=\max\left\{u^-_\jph+c^-_\jph,u^+_\jph+c^+_\jph,0\right\},\quad a^-_\jph=\min\left\{u^-_\jph-c^-_\jph,u^+_\jph-c^+_\jph,0\right\},
\end{equation*}
where $c:=\sqrt{\left[(1+\Gamma)p+\Pi\right]/(\Gamma\rho)}$.
\begin{rmk}
The numerical fluxes \eref{2.1.5} are slightly different from those present in \cite{KLX_21} as one of the goals in \cite{KLX_21} was to
make the resulting scheme well-balanced and thus different values of $\bm U^\pm_\jph$ were used in \eref{2.1.5} and in the computation of
$\bm K^\pm_\jph$.
\end{rmk}

\subsection{Flux Globalization Based Low-Dissipation PCCU Schemes}\label{sec2.2}
In this section, we derived a modified, LD version of the flux globalization based PCCU scheme presented in \S\ref{sec2.1}. We follow the
idea used in \cite{KX_22}, where the LDCU scheme for the single-fluid compressible Euler equations has been introduced. In order to extend
the LDCU scheme to the multifluid case, we now go through all of the derivation steps and begin with the development of the fully discrete
LD PCCU scheme.

\subsubsection{Fully Discrete Scheme}\label{sec2.2.1}
We assume that the computed cell averages $\,\xbar{\bm U}^{\,n}_j:\approx\frac{1}{\dx}\int_{C_j}\bm U(x,t^n)\,{\rm d}x$ are available at a
certain time level $t=t^n$ and use them to reconstruct a second-order piecewise linear interpolant consisting of the linear pieces
$\xbar{\bm U}^{\,n}_j+(\bm U_x)^n_j(x-x_j)$, $x\in C_j$, where the slopes $(\bm U_x)^n_j$ are obtained using a certain nonlinear limiter. We
then estimate the local speeds of propagation $a^\pm_\jph$, introduce the corresponding points $x^n_{\jph,\ell}:=x_\jph+a^-_\jph\dt^n$ and
$x^n_{\jph,r}:=x_\jph+a^+_\jph\dt^n$, and integrate the system \eref{2.1.2} over the space-time control volumes, which consist of the
``smooth'', $[x_{\jmh,r},x_{\jph,\ell}]\times[t^n,t^{n+1}]$, and ``nonsmooth'', $[x_{\jph,\ell},x_{\jph,r}]\times[t^n,t^{n+1}]$, ones, where
$t^{n+1}:=t^n+\dt^n$. This way the solution is evolved in time and upon the completion of the evolution step, we obtain the intermediate
cell averages
\begin{equation}
\hspace*{-0.2cm}\begin{aligned}
\xbar{\bm U}^{\rm\,int}_\jph=\frac{1}{a^+_\jph-a^-_\jph}\Big\{&\bm U^n_{\jph,r}a^+_\jph-\frac{(\bm U_x)^n_{j+1}}{2}\big(a^+_\jph\big)^2\dt^n
-\bm U^n_{\jph,\ell}a^-_\jph+\frac{(\bm U_x)^n_j}{2}\big(a^-_\jph\big)^2\dt^n\\
&-\left[\bm K\big(\bm U^{n+\hf}_{\jph,r}\big)-\bm K\big(\bm U^{n+\hf}_{\jph,\ell}\big)\right]\Big\}
\end{aligned}
\label{2.2.2}
\end{equation}
and
\begin{equation}
\begin{aligned}
\xbar{\bm U}^{\rm\,int}_j&=\,\xbar{\bm U}^{\,n}_j+\frac{(\bm U_x)^n_j}{2}\big(a^+_\jmh+a^-_\jph\big)\dt^n\\
&-\frac{\dt^n}{\dx-\big(a^+_\jmh-a^-_\jph\big)\dt^n}\left[\bm K\big(\bm U^{n+\hf}_{\jph,\ell}\big)-\bm K\big(\bm U^{n+\hf}_{\jmh,r}\big)
\right];
\end{aligned}
\label{2.2.3}
\end{equation}
see \cite{KX_22} for details. In \eref{2.2.2}--\eref{2.2.3}, the point values of $\bm U$ at $\big(x^n_{\jph,\ell},t^n\big)$ and
$\big(x^n_{\jph,r},t^n\big)$ are computed using the piecewise linear reconstruction of $\bm U$, namely,
\begin{equation*}
\begin{aligned}
&\bm U^n_{\jph,\ell}:=\widetilde{\bm U}(x^n_{\jph,\ell},t^n)=\,\xbar{\bm U}^{\,n}_j+(\bm U_x)^n_j\Big(\frac{\dx}{2}+a^-_\jph\dt^n\Big),\\
&\bm U^n_{\jph,r}:=\widetilde{\bm U}(x^n_{\jph,r},t^n)=\,\xbar{\bm U}^{\,n}_{j+1}-(\bm U_x)^n_{j+1}\Big(\frac{\dx}{2}-a^+_\jph\dt^n\Big),
\end{aligned}
\end{equation*}
and the point values of $\bm U$ at $\big(x^n_{\jph,\ell},t^{n+\hf}\big)$ and $\big(x^n_{\jph,r},t^{n+\hf})$, are obtained using the Taylor
expansions about $\big(x_{\jph,r},t^n\big)$ and $\big(x_{\jph,\ell},t^n\big)$, respectively, and the fact that $\bm U_t=-\bm K(\bm U)_x$:
\begin{equation}
\bm U^{n+\hf}_{\jph,\ell}=\bm U^n_{\jph,\ell}-\frac{\dt^n}{2}\bm K\big(\bm U^n_{\jph,\ell}\big)_x,\quad
\bm U^{n+\hf}_{\jph,r}=\bm U^n_{\jph,r}-\frac{\dt^n}{2}\bm K\big(\bm U^n_{\jph,r}\big)_x.
\label{2.2.3.1}
\end{equation}
Here, the slopes $\bm K\big(\bm U^n_{\jph,\ell}\big)_x$ and $\bm K\big(\bm U^n_{\jph,r}\big)_x$ can be computed with the help of a certain
nonlinear limiter; see \cite{KLin} for details.

Next, the intermediate solution, realized in terms of $\{\xbar{\bm U}^{\,\rm int}_j\}$ and $\{\xbar{\bm U}^{\,\rm int}_\jph\}$, is projected
onto the original grid. To this end, we need to construct the interpolant
\begin{equation}
\widetilde{\bm U}^{\,\rm int}(x)=\sum_j\left\{\widetilde{\bm U}^{\,\rm int}_\jph(x){\cal X}_{[x_{\jph,\ell},x_{\jph,r}]}+
\xbar{\bm U}^{\,\rm int}_j{\cal X}_{[x_{\jmh,r},x_{\jph,\ell}]}\right\},
\label{2.2.4}
\end{equation}
where $\cal X$ denotes the characteristic function of the corresponding intervals. We set
\begin{equation}
\widetilde{\bm U}^{\,\rm int}_\jph(x)=\left\{\begin{aligned}
&\xbar{\bm U}^{\,\rm int,L}_\jph,&&x<x_\jph,\\
&\xbar{\bm U}^{\,\rm int,R}_\jph,&&x>x_\jph,
\end{aligned}\right.
\label{2.2.5}
\end{equation}
where the values $\xbar{\bm U}^{\,\rm int,L}_\jph$ and $\xbar{\bm U}^{\,\rm int,R}_\jph$ are determined in several steps. First, according
to the local conservation requirement, the conditions
\begin{equation}
\begin{aligned}
&a^+_\jph\,\xbar{\bm U}^{\,\rm int,R}_\jph-a^-_\jph\,\xbar{\bm U}^{\,\rm int,L}_\jph=(a^+_\jph-a^-_\jph)\,\xbar{\bm U}^{\,\rm int}_\jph
\label{2.2.7}
\end{aligned}
\end{equation}
have to be satisfied, and the rest of the relations on $\xbar{\bm U}^{\,\rm int,L}_\jph$ and $\xbar{\bm U}^{\,\rm int,R}_\jph$ are to be
established for the problem at hand. In fact, the conservation of $\Gamma$ and $\Pi$ components of $\bm U$ is not physically essential, but
the relation \eref{2.2.7} for $\Gamma$ and $\Pi$ are crucial for proving physically relevant properties of the resulting flux globalization
based LD PCCU scheme; see \S\ref{sec2.2.3}.

For the $\gamma$-based multifluid model \eref{2.0}--\eref{2.0.1}, we follow the single-fluid approach from \cite{KX_22} and make the
projection step sharp and accurate for the contact waves, which are linearly degenerate and thus affected by the excessive numerical
dissipation much more than nonlinear shock waves. In order to design such projection step, we consider an isolated contact wave consisting
of the jump discontinuities in $\rho$, $\gamma$, and $\Pi$ propagating in the region with constant $u$ and $p$, and make both $u$ and $p$ to
be constant across the cell interface, namely, we set
\begin{equation}
\begin{aligned}
&\frac{(\xbar{\rho u})^{\rm int,L}_\jph}{\xbar\rho^{\,\rm int,L}_\jph}=
\frac{(\xbar{\rho u})^{\rm int,R}_\jph}{\xbar\rho^{\,\rm int,R}_\jph},\\
&\frac{1}{\xbar\Gamma^{\,\rm int,L}_\jph}\bigg(\xbar E^{\,\rm int,L}_\jph-
\frac{\big((\xbar{\rho u})^{\rm int,L}_\jph\big)^2}{2\,\xbar\rho^{\,\rm int,L}_\jph}-\xbar\Pi^{\,\rm int,L}_\jph\bigg)=
\frac{1}{\xbar\Gamma^{\,\rm int,R}_\jph}\bigg(\xbar E^{\,\rm int,R}_\jph-
\frac{\big((\xbar{\rho u})^{\rm int,R}_\jph\big)^2}{2\,\xbar\rho^{\,\rm int,R}_\jph}-\xbar\Pi^{\,\rm int,R}_\jph\bigg),
\end{aligned}
\label{2.2.8}
\end{equation}
where we have used the EOS \eref{2.0.1}. Next, we solve \eref{2.2.7} and \eref{2.2.8} for $(\xbar{\rho u})^{\rm int,L}_\jph$,
$(\xbar{\rho u})^{\rm int,R}_\jph$, $\xbar E^{\,\rm int,L}_\jph$, and $\xbar E^{\,\rm int,R}_\jph$, and express these quantities in terms of
$\xbar\rho^{\,\rm int,L}_\jph$, $\xbar\rho^{\,\rm int,R}_\jph$, $\xbar\Gamma^{\,\rm int,L}_\jph$, $\xbar\Gamma^{\,\rm int,R}_\jph$,
$\xbar\Pi^{\,\rm int,L}_\jph$, and $\xbar\Pi^{\,\rm int,R}_\jph$:
\begin{equation}
\begin{aligned}
&(\xbar{\rho u})^{\rm int,L}_\jph=\frac{\xbar\rho^{\,\rm int,L}_\jph}{\xbar\rho^{\,\rm int}_\jph}(\xbar{\rho u})^{\rm int}_\jph,\quad
\xbar E^{\,\rm int,L}_\jph=\frac{\xbar\Gamma^{\,\rm int,L}_\jph}{\xbar\Gamma^{\,\rm int}_\jph}\,\xbar E^{\,\rm int}_\jph\\
&\hspace*{1.3cm}+\frac{a^+_\jph\big(\,\xbar\Gamma^{\,\rm int,R}_\jph\,\xbar\rho^{\,\rm int,L}_\jph-
\xbar\Gamma^{\,\rm int,L}_\jph\,\xbar\rho^{\,\rm int,R}_\jph\big)}{2\big(a^+_\jph-a^-_\jph\big)\,\xbar\Gamma^{\,\rm int}_\jph
\big(\,\xbar\rho^{\,\rm int}_\jph\big)^2}\big((\xbar{\rho u})^{\rm int}_\jph\big)^2
+\frac{a^+_\jph\big(\,\xbar\Gamma^{\,\rm int,R}_\jph\,\xbar\Pi^{\,\rm int,L}_\jph-
\xbar\Gamma^{\,\rm int,L}_\jph\,\xbar\Pi^{\,\rm int,R}_\jph\big)}{\big(a^+_\jph-a^-_\jph\big)\,\xbar\Gamma^{\,\rm int}_\jph},\\[0.5ex]
&(\xbar{\rho u})^{\rm int,R}_\jph=\frac{\xbar\rho^{\,\rm int,R}_\jph}{\xbar\rho^{\,\rm int}_\jph}(\xbar{\rho u})^{\rm int}_\jph,\quad
\xbar E^{\,\rm int,R}_\jph=\frac{\xbar\Gamma^{\,\rm int,R}_\jph}{\xbar\Gamma^{\,\rm int}_\jph}\,\xbar E^{\,\rm int}_\jph\\
&\hspace*{1.3cm}+\frac{a^-_\jph\big(\,\xbar\Gamma^{\,\rm int,R}_\jph\,\xbar\rho^{\,\rm int,L}_\jph-
\xbar\Gamma^{\,\rm int,L}_\jph\,\xbar\rho^{\,\rm int,R}_\jph\big)}{2\big(a^+_\jph-a^-_\jph\big)\,\xbar\Gamma^{\,\rm int}_\jph
\big(\,\xbar\rho^{\,\rm int}_\jph\big)^2}\big((\xbar{\rho u})^{\rm int}_\jph\big)^2
+\frac{a^-_\jph\big(\,\xbar\Gamma^{\,\rm int,R}_\jph\,\xbar\Pi^{\,\rm int,L}_\jph-
\xbar\Gamma^{\,\rm int,L}_\jph\,\xbar\Pi^{\,\rm int,R}_\jph\big)}{\big(a^+_\jph-a^-_\jph\big)\,\xbar\Gamma^{\,\rm int}_\jph}.
\end{aligned}
\label{2.2.9}
\end{equation}
Notice that after enforcing \eref{2.2.7} and \eref{2.2.8}, we are left with three degrees of freedom, which we use to make the profiles of
$\rho$, $\Gamma$, and $\Pi$ across the cell interface as sharp as possible and, at the same time, non-oscillatory. This is achieved in the
same way as in \cite{KX_22}, namely, by setting
\begin{equation}
\begin{aligned}
&\xbar\rho^{\,\rm int,L}_\jph=\,\xbar\rho^{\,\rm int}_\jph+\frac{\delta^\rho_\jph}{a^-_\jph},\quad
\xbar\Gamma^{\,\rm int,L}_\jph=\,\xbar\Gamma^{\,\rm int}_\jph+\frac{\delta^\Gamma_\jph}{a^-_\jph},\quad
\xbar\Pi^{\,\rm int,L}_\jph=\,\xbar\Pi^{\,\rm int}_\jph+\frac{\delta^\Pi_\jph}{a^-_\jph},\\
&\xbar\rho^{\,\rm int,R}_\jph=\,\xbar\rho^{\,\rm int}_\jph+\frac{\delta^\rho_\jph}{a^+_\jph},\quad
\xbar\Gamma^{\,\rm int,R}_\jph=\,\xbar\Gamma^{\,\rm int}_\jph+\frac{\delta^\Gamma_\jph}{a^+_\jph},\quad
\xbar\Pi^{\,\rm int,R}_\jph=\,\xbar\Pi^{\,\rm int}_\jph+\frac{\delta^\Pi_\jph}{a^+_\jph},
\end{aligned}
\label{2.2.10}
\end{equation}
where
\begin{equation}
\begin{aligned}
&\delta^\rho_\jph={\rm minmod}\left(-a^-_\jph\big(\,\xbar\rho^{\,\rm int}_\jph-\rho^{\rm int}_{\jph,\ell}\big),
a^+_\jph\big(\rho^{\rm int}_{\jph,r}-\,\xbar\rho^{\,\rm int}_\jph\big)\right),\\
&\delta^\Gamma_\jph={\rm minmod}\left(-a^-_\jph\big(\,\xbar\Gamma^{\,\rm int}_\jph-\Gamma^{\rm int}_{\jph,\ell}\big),
a^+_\jph\big(\Gamma^{\rm int}_{\jph,r}-\,\xbar\Gamma^{\,\rm int}_\jph\big)\right),\\
&\delta^\Pi_\jph={\rm minmod}\left(-a^-_\jph\big(\,\xbar\Pi^{\,\rm int}_\jph-\Pi^{\rm int}_{\jph,\ell}\big),
a^+_\jph\big(\Pi^{\rm int}_{\jph,r}-\,\xbar\Pi^{\,\rm int}_\jph\big)\right),
\end{aligned}
\label{2.2.11}
\end{equation}
and $\bm U^{\,\rm int}_{\jph,\ell}\approx\bm U(x_{\jph,\ell},t^{n+1})$ and $\bm U^{\,\rm int}_{\jph,r}\approx\bm U(x_{\jph,r},t^{n+1})$ are
evaluated similarly to \eref{2.2.3.1}:
\begin{equation}
{\bm U}^{\,\rm int}_{\jph,\ell}=\bm U^n_{\jph,\ell}-\dt^n\bm K\big(\bm U^n_{\jph,\ell}\big)_x,\quad
{\bm U}^{\,\rm int}_{\jph,r}=\bm U^n_{\jph,r}-\dt^n\bm K\big(\bm U^n_{\jph,r}\big)_x.
\label{2.2.12}
\end{equation}
We then complete the construction of $\xbar{\bm U}^{\,\rm int,L}_\jph$ and $\xbar{\bm U}^{\,\rm int,R}_\jph$ by substituting \eref{2.2.10}
into \eref{2.2.9}, which results in
\begin{equation}
\begin{aligned}
&(\xbar{\rho u})^{\rm int,L}_\jph=(\xbar{\rho u})^{\rm int}_\jph+\frac{\delta^\rho_\jph}{a^-_\jph}\,u^{\rm int}_\jph,\quad
(\xbar{\rho u})^{\rm int,R}_\jph=(\xbar{\rho u})^{\rm int}_\jph+\frac{\delta^\rho_\jph}{a^+_\jph}\,u^{\rm int}_\jph,\\
&\xbar E^{\,\rm int,L}_\jph=\bigg(1+\frac{\delta^\Gamma_\jph}{a^-_\jph\xbar\Gamma^{\,\rm int}_\jph}\bigg)\xbar E^{\,\rm int}_\jph+
\frac{\delta^\rho_\jph\xbar\Gamma^{\,\rm int}_\jph-\delta^\Gamma_\jph\,\xbar\rho^{\,\rm int}_\jph}{2a^-_\jph\xbar\Gamma^{\,\rm int}_\jph}
\big(u^{\rm int}_\jph\big)^2+\frac{\delta^\Pi_\jph\xbar\Gamma^{\,\rm int}_\jph-\delta^\Gamma_\jph\xbar\Pi^{\,\rm int}_\jph}
{a^-_\jph\xbar\Gamma^{\,\rm int}_\jph},\\
&\xbar E^{\,\rm int,R}_\jph=\bigg(1+\frac{\delta^\Gamma_\jph}{a^+_\jph\xbar\Gamma^{\,\rm int}_\jph}\bigg)\xbar E^{\,\rm int}_\jph+
\frac{\delta^\rho_\jph\xbar\Gamma^{\,\rm int}_\jph-\delta^\Gamma_\jph\,\xbar\rho^{\,\rm int}_\jph}{2a^+_\jph\xbar\Gamma^{\,\rm int}_\jph}
\big(u^{\rm int}_\jph\big)^2+\frac{\delta^\Pi_\jph\xbar\Gamma^{\,\rm int}_\jph-\delta^\Gamma_\jph\xbar\Pi^{\,\rm int}_\jph}
{a^+_\jph\xbar\Gamma^{\,\rm int}_\jph},
\end{aligned}
\label{2.2.13}
\end{equation}
where $u^{\rm int}_\jph:=(\xbar{\rho u})^{\rm int}_\jph/\,\xbar\rho^{\,\rm int}_\jph$.

Equipped with \eref{2.2.10}--\eref{2.2.13}, we finalize the projection step by integrating the piecewise constant interpolant
\eref{2.2.4}--\eref{2.2.5} over the cell $C_j$. This leads to the following cell averages at the time level $t=t^{n+1}$:
\begin{equation}
\begin{aligned}
\xbar\rho^{\,n+1}_j&=\frac{1}{\dx}\int\limits_{C_j}\widetilde\rho^{\,\rm int}(x)\,{\rm d}x=\,\xbar\rho^{\rm\,int}_j+\frac{\dt^n}{\dx}
\Big[a^+_\jmh\big(\,\xbar\rho^{\,\rm int,R}_\jmh-\xbar\rho^{\,\rm int}_j\big)-a^-_\jph\big(\,\xbar\rho^{\,\rm int,L}_\jph-
\xbar\rho^{\,\rm int}_j\big)\Big]\\
&\stackrel{\eref{2.2.10}}{=}\,\xbar\rho^{\rm\,int}_j+\frac{\dt^n}{\dx}\Big[a^+_\jmh\big(\,\xbar\rho^{\,\rm int}_\jmh-
\xbar\rho^{\,\rm int}_j\big)-a^-_\jph\big(\,\xbar\rho^{\,\rm int}_\jph-\xbar\rho^{\,\rm int}_j\big)+\delta^\rho_\jmh-\delta^\rho_\jph\Big],
\end{aligned}
\label{2.2.14}
\end{equation}
\begin{equation}
\begin{aligned}
(\xbar{\rho u})^{n+1}_j&=\frac{1}{\dx}\int\limits_{C_j}(\widetilde{\rho u})^{\rm int}(x)\,{\rm d}x\\
&=(\xbar{\rho u})^{\rm int}_j+\frac{\dt^n}{\dx}\Big[a^+_\jmh\big((\xbar{\rho u})^{\rm int,R}_\jmh-(\xbar{\rho u})^{\rm int}_j\big)-
a^-_\jph\big((\xbar{\rho u})^{\rm int,L}_\jph-(\xbar{\rho u})^{\rm int}_j\big)\Big]\\[0.5ex]
&\stackrel{\eref{2.2.13}}{=}(\xbar{\rho u})^{\rm int}_j+\frac{\dt^n}{\dx}\Big[a^+_\jmh\big((\xbar{\rho u})^{\rm int}_\jmh-
(\xbar{\rho u})^{\rm int}_j\big)-a^-_\jph\big((\xbar{\rho u})^{\rm int}_\jph-(\xbar{\rho u})^{\rm int}_j\big)\\
&\hspace*{3.5cm}+\delta^\rho_\jmh u^{\rm int}_\jmh-\delta^\rho_\jph u^{\rm int}_\jph\Big],
\end{aligned}
\label{2.2.15}
\end{equation}
\begin{equation}
\begin{aligned}
\xbar E^{\,n+1}_j&=\frac{1}{\dx}\int\limits_{C_j}\widetilde E^{\,\rm int}(x)\,{\rm d}x=\xbar E^{\rm int}_j+\frac{\dt^n}{\dx}
\Big[a^+_\jmh\big(\xbar E^{\,\rm int,R}_\jmh-\xbar E^{\,\rm int}_j\big)-a^-_\jph\big(\xbar E^{\,\rm int,L}_\jph-\xbar E^{\,\rm int}_j
\big)\Big]\\
&\stackrel{\eref{2.2.13}}{=}\,\xbar E^{\rm\,int}_j+\frac{\dt^n}{\dx}\bigg[a^+_\jmh\big(\xbar E^{\,\rm int}_\jmh-\xbar E^{\,\rm int}_j\big)-
a^-_\jph\big(\xbar E^{\,\rm int}_\jph-\xbar E^{\,\rm int}_j\big)\\
&\hspace*{1.85cm}+\frac{\delta^\rho_\jmh\xbar\Gamma^{\,\rm int}_\jph-\delta^\Gamma_\jmh\,\xbar\rho^{\,\rm int}_\jmh}
{2\,\xbar\Gamma^{\,\rm int}_\jmh}\big(u^{\rm int}_\jmh\big)^2-\frac{\delta^\rho_\jph\xbar\Gamma^{\,\rm int}_\jph-
\delta^\Gamma_\jph\,\xbar\rho^{\,\rm int}_\jph}{2\,\xbar\Gamma^{\,\rm int}_\jph}\big(u^{\rm int}_\jph\big)^2\\
&+\frac{\delta^\Pi_\jmh\xbar\Gamma^{\,\rm int}_\jmh-\delta^\Gamma_\jmh\xbar\Pi^{\,\rm int}_\jmh}{\xbar\Gamma^{\,\rm int}_\jmh}-
\frac{\delta^\Pi_\jph\xbar\Gamma^{\,\rm int}_\jph-\delta^\Gamma_\jph\xbar\Pi^{\,\rm int}_\jph}{\xbar\Gamma^{\,\rm int}_\jph}+
\frac{\delta^\Gamma_\jmh}{\xbar\Gamma^{\,\rm int}_\jmh}\,\xbar E^{\,\rm int}_\jmh-\frac{\delta^\Gamma_\jph}{\xbar\Gamma^{\,\rm int}_\jph}
\,\xbar E^{\,\rm int}_\jph\bigg],
\end{aligned}
\label{2.2.16}
\end{equation}
\begin{equation}
\begin{aligned}
\xbar\Gamma^{\,n+1}_j&=\frac{1}{\dx}\int\limits_{C_j}\widetilde\Gamma^{\,\rm int}(x)\,{\rm d}x=\xbar\Gamma^{\rm\,int}_j+
\frac{\dt^n}{\dx}\Big[a^+_\jmh\big(\,\xbar\Gamma^{\,\rm int,R}_\jmh-\xbar\Gamma^{\,\rm int}_j\big)-
a^-_\jph\big(\,\xbar\Gamma^{\,\rm int,L}_\jph-\xbar\Gamma^{\,\rm int}_j\big)\Big]\\
&\stackrel{\eref{2.2.10}}{=}\xbar{\Gamma}^{\rm\,int}_j+\frac{\dt^n}{\dx}\Big[a^+_\jmh\big(\,\xbar\Gamma^{\,\rm int}_\jmh-
\xbar\Gamma^{\,\rm int}_j\big)-a^-_\jph\big(\,\xbar\Gamma^{\,\rm int}_\jph-\xbar\Gamma^{\,\rm int}_j\big)+\delta^\Gamma_\jmh-
\delta^\Gamma_\jph\Big],
\end{aligned}
\label{2.2.17}
\end{equation}
\begin{equation}
\begin{aligned}
\xbar\Pi^{\,n+1}_j&=\frac{1}{\dx}\int\limits_{C_j}\widetilde\Pi^{\,\rm int}(x)\,{\rm d}x=\xbar\Pi^{\rm\,int}_j+\frac{\dt^n}{\dx}
\Big[a^+_\jmh\big(\,\xbar\Pi^{\,\rm int,R}_\jmh-\xbar\Pi^{\,\rm int}_j\big)-
a^-_\jph\big(\,\xbar\Pi^{\,\rm int,L}_\jph-\xbar\Pi^{\,\rm int}_j\big)\Big]\\
&\stackrel{\eref{2.2.10}}{=}\xbar\Pi^{\rm\,int}_j+\frac{\dt^n}{\dx}\Big[a^+_\jmh\big(\,\xbar\Pi^{\,\rm int}_\jmh-\xbar\Pi^{\,\rm int}_j\big)
-a^-_\jph\big(\,\xbar\Pi^{\,\rm int}_\jph-\xbar\Pi^{\,\rm int}_j\big)+\delta^\Pi_\jmh-\delta^\Pi_\jph\Big].
\end{aligned}
\label{2.2.18}
\end{equation}

\subsubsection{Semi-Discrete Scheme}
Finally, we pass to the semi-discrete limit $\dt^n\to0$ in \eref{2.2.14}--\eref{2.2.18} and proceed as in \cite[\S2.2]{KX_22} to end up with
the semi-discretization \eref{2.1.4} with the modified (compared with \eref{2.1.5}) numerical fluxes
\begin{equation}
\bm{{\cal K}}_\jph=\frac{a^+_\jph\bm K^-_\jph-a^-_\jph\bm K^+_\jph}{a^+_\jph-a^-_\jph}+\frac{a^+_\jph a^-_\jph}{a^+_\jph-a^-_\jph}
\left(\bm U^+_\jph-\bm U^-_\jph\right)+\bm q_\jph,
\label{2.2.18a}
\end{equation}
where
\begin{equation}
\begin{aligned}
\bm q_\jph&=q^\rho_\jph\Big(1,u^*_\jph,\frac{\big(u^*_\jph\big)^2}{2},0,0\Big)^\top\\
&+q^\Gamma_\jph\bigg(0,0,\frac{1}{\Gamma^*_\jph}\Big[E^*_\jph-\frac{\big((\rho u)^*_\jph\big)^2}{2\rho^*_\jph}-\Pi^*_\jph\Big],1,0
\bigg)^\top+q^\Pi_\jph(0,0,1,0,1)^\top
\label{2.2.18b}
\end{aligned}
\end{equation}
is a built-in ``anti-diffusion'' term. In \eref{2.2.18b},
\begin{equation}
\begin{aligned}
&\bm U^*_\jph=\frac{a^+_\jph\bm U^+_\jph-a^-_\jph\bm U^-_\jph-\big(\bm K^+_\jph-\bm K^-_\jph\big)}{a^+_\jph-a^-_\jph},\quad
u^*_\jph=\frac{(\rho u)^*_\jph}{\rho^*_\jph},\\
&q^\rho_\jph={\rm minmod}\left(-a^-_\jph\big(\rho^*_\jph-\rho^-_\jph\big),a^+_\jph\big(\rho^+_\jph-\rho^*_\jph\big)\right),\\
&q^\Gamma_\jph={\rm minmod}\left(-a^-_\jph\big(\Gamma^*_\jph-\Gamma^-_\jph\big),a^+_\jph\big(\Gamma^+_\jph-\Gamma^*_\jph\big)\right),\\
&q^\Pi_\jph={\rm minmod}\left(-a^-_\jph\big(\Pi^*_\jph-\Pi^-_\jph\big),a^+_\jph\big(\Pi^+_\jph-\Pi^*_\jph\big)\right).
\end{aligned}
\label{2.2.18c}
\end{equation}
\begin{rmk}
If we replace \eref{2.2.5} with the following limited linear piece:
\begin{equation*}
\widetilde{\bm U}^{\,\rm int}_\jph(x)=\xbar{\bm U}^{\,\rm int}_\jph+
(\bm U_x)^{\rm int}_\jph\left(x-\frac{x_{\jph,r}+x_{\jph,\ell}}{2}\right),
\end{equation*}
where
\begin{equation*}
(\bm U_x)^{\rm int}_\jph={\rm minmod}\Bigg(\,\frac{\bm U^{\rm int}_{\jph,r}-\,\xbar{\bm U}^{\,\rm int}_\jph}{\delta},
\frac{\xbar{\bm U}^{\,\rm int}_\jph-\bm U^{\rm int}_{\jph,\ell}}{\delta}\Bigg),\quad\delta:=\frac{\dt^n}{2}\big(a^+_\jph-a^-_\jph\big),
\end{equation*}
we will end up with an alternative built-in ``anti-diffusion'' term
\begin{equation*}
\bm q_\jph=-\frac{a^+_\jph a^-_\jph}{a^+_\jph-a^-_\jph}\,{\rm minmod}\left(\bm U^*_\jph-\bm U^-_\jph,\bm U^+_\jph-\bm U^*_\jph\right),
\end{equation*}
which was introduced in \cite{KLin}. This will result in another flux globalization based PCCU scheme for the $\gamma$-based multifluid
system. In the numerical results reported in \S\ref{sec4}, we will compare the behavior of this PCCU scheme with the proposed LD PCCU one.
\end{rmk}

\subsubsection{Properties of the Semi-Discrete Scheme}\label{sec2.2.3}
In this section, we establish two important properties the designed semi-discrete scheme satisfies:

1. When initially $\gamma\equiv{\rm Const}$ and $\pi_\infty\equiv{\rm Const}$, that is, if initially the system contains a single fluid,
then $\gamma$ and $\pi_\infty$ will stay constant for all $t$;

2. At the isolated material interface at which initially $u\equiv{\rm Const}$ and $p\equiv{\rm Const}$, both $u$ and $p$ will stay constant
for all $t$.

\noindent
To this end, we prove the following theorem.
\begin{thm}\label{thm1}
1. If $\,\xbar\Gamma_j\equiv\widehat\Gamma={\rm Const}$ and $\,\xbar\Pi_j\equiv\widehat\Pi={\rm Const}$ for all $j$ at a certain time level
$t$, then
\begin{equation*}
\frac{{\rm d}}{{\rm d}t}\,\xbar\Gamma_j=\frac{{\rm d}}{{\rm d}t}\,\xbar\Pi_j\equiv0,\quad\forall j.
\end{equation*}
2. If at a certain time level $t=t^n$, $u^n_j\equiv\widehat u={\rm Const}$ and $p^n_j\equiv\widehat p={\rm Const}$ for all $j$, then at the
next time level $t=t^{n+1}$,
\begin{equation}
u^{n+1}_j\equiv\widehat u\quad{\rm and}\quad p^{n+1}_j\equiv\widehat p,\quad\forall j,
\label{pr2.0}
\end{equation}
provided the ODE system \eref{2.1.4}, \eref{2.2.18a}--\eref{2.2.18c} is discretized using the forward Euler method.
\end{thm}
{\bf Proof:} 1. We only show that $\frac{{\rm d}}{{\rm d}t}\,\xbar\Gamma_j\equiv0$ as $\frac{{\rm d}}{{\rm d}t}\,\xbar\Pi_j\equiv0$ can be
proved in a similar way. Since the point values of $\Gamma$ at the cell interfaces are obtained using the piecewise linear reconstruction
\eref{2.1.17}, we have $\Gamma^+_\jph=\Gamma^-_\jph=\Gamma^+_\jmh\equiv\widehat\Gamma$, which we substitute into \eref{2.1.19a},
\eref{2.1.19}, and \eref{2.1.21} to obtain
\begin{equation}
\big(F^{(4)}\big)^\pm_\jph=\widehat\Gamma u^\pm_\jph,\quad B^{(4)}_j=\widehat\Gamma\big(u^-_\jph-u^+_\jmh\big),\quad
B^{(4)}_{\bm\Psi,\jph}=\widehat\Gamma\big(u^+_\jph-u^-_\jph\big).
\label{pr1.1}
\end{equation}
We then use \eref{2.1.6}--\eref{2.1.9} to compute the flux differences
\begin{equation}
\begin{aligned}
&\big(K^{(4)}\big)^+_\jph-\big(K^{(4)}\big)^-_\jph=\big(F^{(4)}\big)^+_\jph-\big(F^{(4)}\big)^-_\jph-B^{(4)}_{\bm\Psi,\jph}
\stackrel{\eref{pr1.1}}{=}0,\\
&\big(K^{(4)}\big)^-_\jph-\big(K^{(4)}\big)^+_\jmh=\big(F^{(4)}\big)^-_\jph-\big(F^{(4)}\big)^+_\jmh-B^{(4)}_j\stackrel{\eref{pr1.1}}{=}0,
\end{aligned}
\label{pr1.3}
\end{equation}
which we substitute into \eref{2.2.18c} to verify that
\begin{equation}
q^\Gamma_\jph=0.
\label{pr1.5}
\end{equation}
Finally, we substitute \eref{pr1.3}--\eref{pr1.5} into \eref{2.1.4}, \eref{2.2.18a}--\eref{2.2.18b} to end up with
$\frac{{\rm d}}{{\rm d}t}\,\xbar\Gamma_j\equiv0$. This completes the proof of the first part of the theorem.

\smallskip
2. Since the point values of $u$ and $p$ at the cell interfaces are computed using the piecewise linear reconstruction \eref{2.1.17}, we
have $u^+_\jph=u^-_\jph=u^+_\jmh\equiv\widehat u$ and $p^+_\jph=p^-_\jph=p^+_\jmh\equiv\widehat p$ at the time level $t=t^n$. This results
in
\begin{equation}
\begin{aligned}
&\bm U^\pm_\jph=\big(\rho^\pm_\jph,\rho^\pm_\jph\widehat u,E^\pm_\jph,\Gamma^\pm_\jph,\Pi^\pm_\jph\big)^\top,\quad
E^\pm_\jph=\widehat p\,\Gamma^\pm_\jph+\frac{\widehat u^{\,2}}{2}\rho^\pm_\jph+\Pi^\pm_\jph,\\
&\bm F^\pm_\jph=\big(\rho^\pm_\jph\widehat u,\rho^\pm_\jph\widehat u^{\,2}+\widehat p,\widehat u\,(E^\pm_\jph+\widehat p\,),
\Gamma^\pm_\jph\widehat u,\Pi^\pm_\jph\widehat u\big)^\top,
\label{pr2.1}
\end{aligned}
\end{equation}
and hence, after substituting $u^\pm_\jph$ into \eref{2.1.19} and \eref{2.1.21}, we obtain
\begin{equation*}
\bm B_j=\bm B_{\bm\Psi,\jph}\equiv\bm0.
\end{equation*}
The latter equality implies $\bm K^\pm_\jph=\bm F^\pm_\jph$, which together with \eref{pr2.1} results in
\begin{equation*}
\bm F^+_\jph-\bm F^-_\jph=\widehat u\,(\bm U^+_\jph-\bm U^+_\jph),
\end{equation*}
so that the first line in \eref{2.2.18c} can be rewritten as
\begin{equation}
\bm U^*_\jph = \frac{a^+_\jph\bm U^+_\jph-a^-_\jph\bm U^-_\jph-\widehat u\,\big(\bm U^+_\jph-\bm U^-_\jph\big)}{a^+_\jph-a^-_\jph},\quad
u^*_\jph=\widehat u.
\label{pr2.3}
\end{equation}
We then use \eref{2.2.18b}, \eref{2.2.18c}, \eref{pr2.1}, and \eref{pr2.3} to compute the ``anti-diffusion'' term, which reduces to
\begin{equation}
\bm q_\jph=\Big(q^\rho_\jph,\widehat u\,q^\rho_\jph,\widehat p\,q^\Gamma_\jph+\frac{\widehat u^{\,2}}{2}\,q^\rho_\jph+q^\Pi_\jph,
q^\Gamma_\jph,q^\Pi_\jph\Big)^\top,
\label{pr2.4}
\end{equation}
and then substitute \eref{pr2.1} and \eref{pr2.4} into \eref{2.2.18a} to obtain the numerical fluxes
\begin{equation*}
\bm{{\cal K}}_\jph=\Big({\cal K}^{(1)}_\jph,\widehat u\,{\cal K}^{(1)}_\jph+\widehat p,\widehat p\,{\cal K}^{(4)}_\jph+
\frac{\widehat u^{\,2}}{2}{\cal K}^{(1)}_\jph+{\cal K}^{(5)}_\jph+\widehat u\widehat p,{\cal K}^{(4)}_\jph,{\cal K}^{(5)}_\jph\Big)^\top,
\end{equation*}
which, in turn, are substituted into \eref{2.1.4} to obtain the following semi-discrete relations:
\begin{equation*}
\frac{{\rm d}}{{\rm d}t}(\xbar{\rho u})^n_j=\widehat u\,\frac{{\rm d}}{{\rm d}t}\,\xbar\rho^{\,n}_j,\quad
\frac{{\rm d}}{{\rm d}t}\xbar E^{\,n}_j=\widehat p\,\frac{{\rm d}}{{\rm d}t}\xbar\Gamma^{\,n}_j+\frac{\widehat u^{\,2}}{2}
\frac{{\rm d}}{{\rm d}t}\,\xbar\rho^{\,n}_j+\frac{{\rm d}}{{\rm d}t}\,\xbar\Pi^{\,n}_j.
\end{equation*}
These relations are discretized using the forward Euler method, which results in
\begin{equation}
\begin{aligned}
&\frac{(\xbar{\rho u})^{n+1}_j-(\xbar{\rho u})^n_j}{\dt}=\widehat u\,\frac{\xbar\rho^{\,n+1}_j-\xbar\rho^{\,n}_j}{\dt},\\
&\frac{\xbar E^{\,n+1}_j-\xbar E^{\,n}_j}{\dt}=\widehat p\,\frac{\xbar\Gamma^{\,n+1}_j-\xbar\Gamma^{\,n}_j}{\dt}+
\frac{\widehat u^{\,2}}{2}\cdot\frac{\xbar\rho^{\,n+1}_j-\xbar\rho^{\,n}_j}{\dt}+\frac{\xbar\Pi^{\,n+1}_j-\xbar\Pi^{\,n}_j}{\dt}.
\end{aligned}
\label{pr2.6}
\end{equation}
Finally, we substitute \eref{2.1.13} expressed at both time levels $t=t^n$ and $t=t^{n+1}$ into \eref{pr2.6} to end up with \eref{pr2.0}.
This completes the proof of the second part of the theorem.$\hfill\blacksquare$
\begin{rmk}
The second part of Theorem \ref{thm1} is still true if the forward Euler time discretization is replaced with another strong stability
preserving (SSP) ODE solver; see, e.g.,\cite{Gottlieb11,Gottlieb12}.
\end{rmk}
\begin{rmk}
As in \cite{KX_22}, the computation of numerical fluxes in \eref{2.2.18a} should be desingularized to avoid division by zero or very small
numbers. If $a^+_\jph<\varepsilon_0$ and $a^-_\jph>-\varepsilon_0$ for a small positive $\varepsilon_0$, we replace the fluxes
$\bm{{\cal K}}_\jph$ with
\begin{equation*}
\bm{{\cal K}}_\jph=\frac{\bm K\big(\bm U^-_\jph\big)+\bm K\big(\bm U^+_\jph\big)}{2}.
\end{equation*}
In all of the numerical examples reported in \S\ref{sec4}, we have taken $\varepsilon_0=10^{-12}$.
\end{rmk}

\subsection{Flux Globalization Based LD Ai-WENO PCCU Scheme}\label{sec2.3}
In this section, we extend the second-order flux globalization based LD PCCU schemes from \S\ref{sec2.2} to the fifth order of accuracy
within the A-WENO framework.

The semi-discrete fifth-order LD Ai-WENO PCCU scheme for the 1-D quasi-conservative system \eref{2.1.2} reads as
\begin{equation}
\frac{{\rm d}}{{\rm d}t}\mU_j=-\frac{\bmH_\jph-\bmH_\jmh}{\dx},
\label{2.53}
\end{equation}
where $\mU_j:\approx\mU(x_j,t)$ and the fifth-order numerical fluxes $\bmH_\jph$ are defined by
\begin{equation}
\bmH_\jph=\bm{{\cal K}}_\jph-\frac{\dx}{24}(\bm K_{xx})_\jph+\frac{7(\dx)^3}{5760}(\bm K_{xxxx})_\jph.
\label{2.53a}
\end{equation}
Here, $\bm{{\cal K}}_\jph$ are the finite-volume fluxes \eref{2.2.18a}--\eref{2.2.18c}, and $(\bm K_{xx})_\jph$ and $({\bm K_{xxxx}})_\jph$
are approximations of the second- and fourth-order spatial derivatives of $\bm K$ at $x=x_\jph$, which we compute using the
finite-difference approximations recently proposed in \cite{CKX23}:
\begin{equation}
\begin{aligned}
&(\mK_{xx})_\jph=\frac{1}{12(\dx)^2}\left[-\bm{{\cal K}}_{j-\frac{3}{2}}+16\bm{{\cal K}}_\jmh-30\bm{{\cal K}}_\jph+
16\bm{{\cal K}}_{j+\frac{3}{2}}-\bm{{\cal K}}_{j+\frac{5}{2}}\right],\\
&(\mK_{xxxx})_\jph=\frac{1}{(\dx)^4}\left[\bm{{\cal K}}_{j-\frac{3}{2}}-4\bm{{\cal K}}_\jmh+6\bm{{\cal K}}_\jph-
4\bm{{\cal K}}_{j+\frac{3}{2}}+\bm{{\cal K}}_{j+\frac{5}{2}}\right].
\end{aligned}
\label{2.57}
\end{equation}
The resulting semi-discrete scheme \eref{2.53}--\eref{2.57} will be fifth-order accurate provided the point values $\bm U^\pm_\jph$ are
calculated using a fifth-order interpolation. To this end, we apply the recently proposed fifth-order Ai-WENO-Z interpolation
\cite{DLWW22,LLWDW23,WD22}, which we briefly describe in Appendix \ref{appa}.

\section{Two-Dimensional Algorithms}\label{sec3}
In this section, we extend the proposed 1-D flux globalization based LD A-WENO PCCU schemes to the 2-D $\gamma$-based multifluid system
\eref{1.1}, \eref{1.2}, \eref{1.5}. This system can be written in the vector form
\begin{equation*}
\bm U_t+\bm F(\bm U)_x+\bm G(\bm U)_y=B(\bm U)\bm U_x+C(\bm U)\bm U_y,
\end{equation*}
or, equivalently, in the quasi-conservative form
\begin{equation}
\bm U_t+\bm K(\bm U)_x+\bm L(\bm U)_y=\bm0
\label{3.1}
\end{equation}
with
\begin{equation}
\begin{aligned}
&\bm K(\bm U)=\bm F(\bm U)-\bm R(\bm U),&&\bm L(\bm U)=\bm G(\bm U)-\bm S(\bm U),\\
&\bm R(\bm U)=\int\limits^x_{\widehat x}B(\bm U)\bm U_\xi\,{\rm d}\xi,&&\bm S(\bm U)=\int\limits^y_{\widehat y}C(\bm U)\bm U_\eta\,
{\rm d}\eta.
\end{aligned}
\label{3.2}
\end{equation}
Here,
\begin{equation}
\begin{aligned}
&\bm U:=(\rho,\rho u,\rho v,E,\Gamma,\Pi)^\top,&&~\\
&\bm F(\bm U)=(\rho u,\rho u^2+p,\rho uv,u(E+p),u\Gamma,u\Pi)^\top,&&B(\bm U)\bm U_x=\big(0,0,0,0,\Gamma u_x,\Pi u_x\big)^\top,\\
&\bm G(\bm U)=(\rho v,\rho uv,\rho v^2+p,v(E+p),v\Gamma,v\Pi)^\top,&&C(\bm U)\bm U_y=\big(0,0,0,0,\Gamma v_y,\Pi v_y\big)^\top.
\end{aligned}
\label{3.3}
\end{equation}

We first introduce a uniform mesh consisting of the finite-volume cells $C_{j,k}:=[x_\jmh,x_\jph]\times[y_\kmh,y_\kph]$ of the uniform size
$\dx\dy$ with $x_\jph-x_\jmh\equiv\dx$ and $y_\kph-y_\kmh\equiv\dy$ centered at $(x_j,y_k)$ with $x_j=(x_\jmh+x_\jph)/2$ and
$(y_\kmh+y_\kph)/2$, $j=1,\ldots,N_x$, $k=1,\ldots,N_y$.

We assume that at certain time level $t$, an approximate solution, realized in terms of the cell averages
$\xbar{\bm U}_{j,k}:\approx\frac{1}{\dx\dy}\iint_{C_{j,k}}\bm U(x,y,t)\,{\rm d}x\,{\rm d}y$, is available. These cell averages are then
evolved in time by solving the following system of ODEs:
\begin{equation}
\frac{{\rm d}}{{\rm d}t}\,\xbar{\bm U}_{j,k}=-\frac{\bm{{\cal K}}_{\jph,k}-\bm{{\cal K}}_{\jmh,k}}{\dx}-
\frac{\bm{{\cal L}}_{j,\kph}-\bm{{\cal L}}_{j,\kmh}}{\dy},
\label{3.4}
\end{equation}
where the $x$- and $y$-numerical fluxes are
\begin{align}
&\bm{{\cal K}}_{\jph,k}=\frac{a^+_{\jph,k}\bm K^-_{\jph,k}-a^-_{\jph,k}\bm K^+_{\jph,k}}{a^+_{\jph,k}-a^-_{\jph,k}}+
\frac{a^+_{\jph,k}a^-_{\jph,k}}{a^+_{\jph,k}-a^-_{\jph,k}}\left(\bm U^+_{\jph,k}-\bm U^-_{\jph,k}\right)+\bm q_{\jph,k},\label{3.5}\\
&\bm{{\cal L}}_{j,\kph}=\frac{b^+_{j,\kph}\bm L^-_{j,\kph}-b^-_{j,\kph}\bm L^+_{j,\kph}}{b^+_{j,\kph}-b^-_{j,\kph}}+
\frac{b^+_{j,\kph}b^-_{j,\kph}}{b^+_{j,\kph}-b^-_{j,\kph}}\left(\bm U^+_{j,\kph}-\bm U^-_{j,\kph}\right)+\bm q_{j,\kph}.\label{3.6}
\end{align}

The one-sided point values $\bm U^\pm_{\jph,k}$ and $\bm U^\pm_{j,\kph}$ at the cell interfaces $(x_\jph,y_k)$ and $(x_j,y_\kph)$,
respectively, are obtained as follows. We first use the cell averages $\xbar{\bm{U}}_{j,k}$ to compute the point values of $u$, $v$, and $p$
at the cell centers:
\begin{equation*}
u_{j,k}=\frac{(\xbar{\rho u})_{j,k}}{\xbar\rho_{j,k}},\quad v_{j,k}=\frac{(\xbar{\rho v})_{j,k}}{\xbar\rho_{j,k}},\quad
p_{j,k}=\frac{1}{\xbar\Gamma_{j,k}}\left[\xbar E_{j,k}-\frac{\big((\xbar{\rho u})_{j,k}\big)^2+\big((\xbar{\rho v})_{j,k}\big)^2}
{2\,\xbar\rho_{j,k}}-\xbar\Pi_{j,k}\right],
\end{equation*}
and then construct the linear pieces to approximate the primitive variables ${\bm V}=(\rho,u,v,p,\Gamma,\Pi)^\top$:
\begin{equation}
\widetilde{\bm V}_{j,k}(x,y)=\bm V_{j,k}+(\bm V_x)_{j,k}(x-x_j)+(\bm V_y)_{j,k}(y-y_k),\quad(x,y)\in C_{j,k},
\label{3.7}
\end{equation}
where $\bm V_{j,k}:=(\,\xbar\rho_{j,k},u_{j,k},v_{j,k},p_{j,k},\xbar\Gamma_{j,k},\xbar\Pi_{j,k})^\top$ and $(\bm V_x)_{j,k}$ and
$(\bm V_y)_{j,k}$ are the slopes, which are supposed to be computed using a nonlinear limiter.

As in the 1-D case, we use different limiters near and away from the material interfaces, which need to be detected. In the two-fluid case,
we check whether
\begin{equation}
(\xbar\Gamma_{j.k}-\widehat\Gamma)(\xbar\Gamma_{j+1,k}-\widehat\Gamma)<0,
\label{3.7a}
\end{equation}
where, as before, $\widehat\Gamma=(\Gamma_{\rm I}+\Gamma_{\rm II})/2$. If \eref{3.7a} is satisfied, we then use the overcompressive SBM
limiter,
\begin{equation}
(\bm V_x)_{\ell,k}=\phi^{\rm SBM}_{\theta,\tau}\left(\frac{\xbar{\bm V}_{\ell+1,k}-\xbar{\bm V}_{\ell,k}}
{\xbar{\bm V}_{\ell,k}-\xbar{\bm V}_{\ell-1,k}}\right)\frac{\xbar{\bm V}_{\ell+1,k}-\xbar{\bm V}_{\ell,k}}{\dx},
\label{3.7b}
\end{equation}
for $\ell=j-1$, $j$, $j+1$, and $j+2$. In \eref{3.7b}, the function $\phi^{\rm SBM}_{\theta,\tau}(r)$, given by \eref{2.1.15c}, is applied
in a componentwise manner with $\tau=-0.5$ and $\theta=1.3$. Otherwise, that is, away from the material interface, we use a dissipative
generalized minmod limiter which is also given by the same formulae \eref{3.7b}, \eref{2.1.15c}, but with $\tau=0.5$ and $\theta=1.3$. We
proceed similarly in the $y$-direction: we use
\begin{equation}
(\bm V_y)_{j,m}=\phi^{\rm SBM}_{\theta,\tau}\left(\frac{\xbar{\bm V}_{j,m+1}-\xbar{\bm V}_{j,m}}
{\xbar{\bm V}_{j,m}-\xbar{\bm V}_{j,m-1}}\right)\frac{\xbar{\bm V}_{j,m+1}-\xbar{\bm V}_{j,m}}{\dy},
\label{3.7bf}
\end{equation}
with $\tau=-0.5$ and $\theta=1.3$ for $m=k-1$, $k$, $k+1$, and $k+2$ if
\begin{equation}
(\xbar\Gamma_{j,k}-\widehat\Gamma)(\xbar\Gamma_{j,k+1}-\widehat\Gamma)<0,
\label{3.7af}
\end{equation}
is satisfied, and with $\tau=0.5$ and $\theta=1.3$ otherwise.

Equipped with $(\bm V_x)_{j,k}$ and $(\bm V_y)_{j,k}$, we use \eref{3.7} to obtain
\begin{equation*}
\begin{aligned}
&\bm V^-_{\jph,k}=\,\xbar{\bm V}_{j,k}+\frac{\dx}{2}(\bm V_x)_{j,k},\quad
\bm V^+_{\jph,k}=\,\xbar{\bm V}_{j+1,k}-\frac{\dx}{2}(\bm V_x)_{j+1,k},\\
&\bm V^-_{j,\kph}=\,\xbar{\bm V}_{j,k}+\frac{\dy}{2}(\bm V_y)_{j,k},\quad
\bm V^+_{j,\kph}=\,\xbar{\bm V}_{j,k+1}-\frac{\dy}{2}(\bm V_y)_{j,k+1},
\end{aligned}
\end{equation*}
and then the corresponding point values of the conservative variables $\bm U$ are
\begin{equation*}
{\bm U}^\pm_{\ell,m}=\left(\rho^\pm_{\ell,m},\rho^\pm_{\ell,m}u^\pm_{\ell,m},\rho^\pm_{\ell,m}v^\pm_{\ell,m},E^\pm_{\ell,m},
\Gamma^\pm_{\ell,m},\Pi^\pm_{\ell,m}\right)^\top,
\end{equation*}
for $(\ell,m)=(\jph,k)$ and $(j,\kph)$. Here,
$E^\pm_{\ell,m}=\Gamma^\pm_{\ell,m}p^\pm_{\ell,m}+\rho^\pm_{\ell,m}((u^\pm_{\ell,m})^2+(v^\pm_{\ell,m})^2)/2+\Pi^\pm_{\ell,m}$.

The global fluxes $\bm K^\pm_{\jph,k}$ and $\bm{{\cal L}}_{j,\kph}$ in \eref{3.5}--\eref{3.6} are obtained using the relations in
\eref{3.2}, namely, by
\begin{equation*}
\bm K^\pm_{\jph,k}=\bm F^\pm_{\jph,k}-\bm R^\pm_{\jph,k},\quad\bm L^\pm_{j,\kph}=\bm G^\pm_{j,\kph}-\bm S^\pm_{j,\kph},
\end{equation*}
where $\bm F^\pm_{\jph,k}:=\bm F(\bm U^\pm_{\jph,k})$, $\bm G^\pm_{j,\kph}:=\bm G(\bm U^\pm_{j,\kph})$, and the point values of the global
variables $\bm R$ and $\bm S$ are computed as follows. First, we set $\widehat x=x_\hf$ and $\widehat y=y_\hf$ so that
$\bm R^-_{\hf,k}:=\bm0$ and $\bm S^-_{j,\hf}:=\bm0$. We then evaluate $\bm R^+_{\hf,k}=\bm B_{\bm\Psi,\hf,k}$ and
$\bm S^+_{j,\hf}=\bm B_{\bm\Psi,j,\hf}$ and recursively compute the rest of the required point values:
\begin{equation*}
\begin{aligned}
&\bm R^-_{\jph,k}=\bm R^+_{\jmh,k}+\bm B^x_{j,k},\quad\bm R^+_{\jph,k}=\bm R^-_{\jph,k}+\bm B_{\bm\Psi,\jph,k},\quad j=1,\ldots,N_x,\\
&\bm S^-_{j,\kph}=\bm S^+_{j,\kmh}+\bm B^y_{j,k},\quad\bm S^+_{j,\kph}=\bm S^-_{j,\kph}+\bm B_{\bm\Psi,j,\kph},\quad k=1,\ldots,N_y.
\end{aligned}
\end{equation*}
Here, $\bm B^x_{j,k}$, $\bm B_{\bm\Psi,\jph,k}$, $\bm B^y_{j,k}$, and $\bm B_{\bm\Psi,j,\kph,k}$ are evaluated precisely the same way as in
\S\ref{sec2.1.1}:
\begin{equation*}
\begin{aligned}
&\bm B^x_{j,k}=\Big(0,0,0,0,\frac{\Gamma^-_{\jph,k}+\Gamma^+_{\jmh,k}}{2}\big(u^-_{\jph,k}-u^+_{\jmh,k}\big),
\frac{\Pi^-_{\jph,k}+\Pi^+_{\jmh,k}}{2}\big(u^-_{\jph,k}-u^+_{\jmh,k}\big)\Big)^\top,\\
&\bm B_{\bm\Psi,\jph,k}=\Big(0,0,0,0,\frac{\Gamma^+_{\jph,k}+\Gamma^-_{\jph,k}}{2}\big(u^+_{\jph,k}-u^-_{\jph,k}\big),
\frac{\Pi^+_{\jph,k}+\Pi^-_{\jph,k}}{2}\big(u^+_{\jph,k}-u^-_{\jph,k}\big)\Big)^\top,\\
&\bm B^y_{j,k}=\Big(0,0,0,0,\frac{\Gamma^-_{j,\kph}+\Gamma^+_{j,\kmh}}{2}\big(v^-_{j,\kph}-v^+_{j,\kmh}\big),
\frac{\Pi^-_{j,\kph}+\Pi^+_{j,\kmh}}{2}\big(v^-_{j,\kph}-v^+_{j,\kmh}\big)\Big)^\top,\\
&\bm B_{\bm\Psi,j,\kph}=\Big(0,0,0,0,\frac{\Gamma^+_{j,\kph}+\Gamma^-_{j,\kph}}{2}\big(v^+_{j,\kph}-v^-_{j,\kph}\big),
\frac{\Pi^+_{j,\kph}+\Pi^-_{j,\kph}}{2}\big(v^+_{j,\kph}-v^-_{j,\kph}\big)\Big)^\top.
\end{aligned}
\end{equation*}

Finally, the one-sided local-speeds of propagation $a^\pm_{\jph,k}$ and $b^\pm_{j,\kph}$ can be estimated by
\begin{equation*}
\begin{aligned}
&a^+_{\jph,k}=\max\left\{u^-_{\jph,k}+c^-_{\jph,k},u^+_{\jph,k}+c^+_{\jph,k},0\right\},\\
&a^-_{\jph,k}=\min\left\{u^-_{\jph,k}-c^-_{\jph,k},u^+_{\jph,k}-c^+_{\jph,k},0\right\},\\
&b^+_{j,\kph}=\max\left\{v^-_{j,\kph}+c^-_{j,\kph},v^+_{j,\kph}+c^+_{j,\kph},0\right\},\\
&b^-_{j,\kph}=\min\left\{v^-_{j,\kph}-c^-_{j,\kph},v^+_{j,\kph}-c^+_{j,\kph},0\right\}.
\end{aligned}
\end{equation*}

\subsection{``Built-in'' Anti-Diffusion}\label{sec3.2}
In this section, we discuss the derivation of the ``built-in'' anti-diffusion terms $\bm q_{\jph,k}$ and $\bm q_{j,\kph}$ in
\eref{3.5}--\eref{3.6} in a ``dimension-by-dimension'' manner following the idea introduced in \cite{KX_22}.

In order to derive the formula for $\bm q_{\jph,k}$, we consider the 1-D restriction of the system \eref{3.1}--\eref{3.3} along the lines
$y=y_k$:
\begin{equation}
\bm U_t(x,y_k,t)+\bm K\big(\bm U(x,y_k,t)\big)_x=\bm0,\quad k=1,\ldots,N_y.
\label{3.17}
\end{equation}
We then can go through all of the steps in the derivation of the 1-D fully discrete scheme for the systems in \eref{3.17} following
\S\ref{sec2.2.1} up to \eref{2.2.5}, which now reads as
\begin{equation*}
\widetilde{\bm U}^{\,\rm int}_{\jph,k}(x,y_k)=\left\{\begin{aligned}
\xbar{\bm U}^{\,\rm int,L}_{\jph,k},\quad x<x_\jph,\\
\xbar{\bm U}^{\,\rm int,R}_{\jph,k},\quad x>x_\jph,
\end{aligned}\right.
\end{equation*}
and the corresponding local conservation requirements \eref{2.2.7} become
\begin{equation}
a^+_{\jph,k}\,\xbar{\bm U}^{\,\rm int,R}_{\jph,k}-a^-_{\jph,k}\,\xbar{\bm U}^{\,\rm int,L}_{\jph,k}=
\big(a^+_{\jph,k}-a^-_{\jph,k}\big)\,\xbar{\bm U}^{\,\rm int}_{\jph,k}.
\label{3.18}
\end{equation}
In addition to the six conservation constraints given by \eref{3.18}, we have six degrees of freedom, which we use as in the 1-D case to
enforce a sharp approximation of quasi 1-D isolated contact waves propagating in the $x$-direction. To this end, we enforce the continuity
of $u$ and $p$ across the cell interfaces $x=x_\jph$ by setting
\begin{equation*}
\begin{aligned}
\frac{(\xbar{\rho u})^{\,\rm int,L}_{\jph,k}}{\xbar\rho^{\,\rm int,L}_{\jph,k}}=
\frac{(\xbar{\rho u})^{\,\rm int,R}_{\jph,k}}{\xbar\rho^{\,\rm int,R}_{\jph,k}},\qquad&
\frac{1}{\xbar\Gamma^{\,\rm int,L}_{\jph,k}}\bigg(\xbar E^{\,\rm int,L}_{\jph,k}-\frac{\big((\xbar{\rho u})^{\,\rm int,L}_{\jph,k}\big)^2+
\big((\xbar{\rho v})^{\,\rm int,L}_{\jph,k}\big)^2}{2\,\xbar\rho^{\,\rm int,L}_{\jph,k}}-\xbar \Pi^{\,\rm int,L}_{\jph,k}\bigg)\\
=&\,\frac{1}{\xbar\Gamma^{\,\rm int,R}_{\jph,k}}\bigg(\xbar E^{\,\rm int,R}_{\jph,k}-
\frac{\big((\xbar{\rho u})^{\,\rm int,R}_{\jph,k}\big)^2+
\big((\xbar{\rho v})^{\,\rm int,R}_{\jph,k}\big)^2}{2\,\xbar\rho^{\,\rm int,R}_{\jph,k}}-\xbar \Pi^{\,\rm int,R}_{\jph,k}\bigg),
\end{aligned}
\end{equation*}
and then proceed as in \S\ref{sec2.2.1} to enforce sharp (yet, non-oscillatory) jumps of the $\rho$-, $\rho v$-, $\Gamma$-, and
$\Pi$-components. This leads to the following formulae analogous to \eref{2.2.10}--\eref{2.2.11}:
\begin{equation*}
\begin{aligned}
\xbar\rho^{\,\rm int,L}_{\jph,k}&=\,\xbar\rho^{\,\rm int}_{\jph,k}+\frac{\delta^\rho_{\jph,k}}{a^-_{\jph,k}},\quad&
\xbar\rho^{\,\rm int,R}_{\jph,k}&=\,\xbar\rho^{\,\rm int}_{\jph,k}+\frac{\delta^\rho_{\jph,k}}{a^+_{\jph,k}},\\[0.2ex]
(\xbar{\rho v})^{\,\rm int,L}_{\jph,k}&=(\xbar{\rho v})^{\,\rm int}_{\jph,k}+\frac{\delta^{\rho v}_{\jph,k}}{a^-_{\jph,k}},\quad&
(\xbar{\rho v})^{\,\rm int,R}_{\jph,k}&=(\xbar{\rho v})^{\,\rm int}_{\jph,k}+\frac{\delta^{\rho v}_{\jph,k}}{a^+_{\jph,k}},\\[0.2ex]
\xbar\Gamma^{\,\rm int,L}_{\jph,k}&=\,\xbar\Gamma^{\,\rm int}_{\jph,k}+\frac{\delta^\Gamma_{\jph,k}}{a^-_{\jph,k}},\quad&
\xbar\Gamma^{\,\rm int,R}_{\jph,k}&=\,\xbar\Gamma^{\,\rm int}_{\jph,k}+\frac{\delta^\Gamma_{\jph,k}}{a^+_{\jph,k}},\\[0.2ex]
\xbar\Pi^{\,\rm int,L}_{\jph,k}&=\,\xbar\Pi^{\,\rm int}_{\jph,k}+\frac{\delta^\Pi_{\jph,k}}{a^-_{\jph,k}},\quad&
\xbar\Pi^{\,\rm int,R}_{\jph,k}&=\,\xbar\Pi^{\,\rm int}_{\jph,k}+\frac{\delta^\Pi_{\jph,k}}{a^+_{\jph,k}},
\end{aligned}
\end{equation*}
where
\begin{equation*}
\begin{aligned}
&\delta^\rho_{\jph,k}={\rm minmod}\left(-a^-_{\jph,k}\big[\,\xbar\rho^{\,\rm int}_{\jph,k}-\big(\rho^{\rm int}_{\jph,k}\big)_\ell\big],\,
a^+_{\jph,k}\big[\big(\rho^{\rm int}_{\jph,k}\big)_r-\,\xbar\rho^{\,\rm int}_{\jph,k}\big]\right),\\
&\delta^{\rho v}_{\jph,k}={\rm minmod}\left(-a^-_{\jph,k}\big[(\xbar{\rho v})^{\,\rm int}_{\jph,k}-
\big((\rho v)^{\rm int}_{\jph,k}\big)_\ell\big],\,a^+_{\jph,k}\big[\big((\rho v)^{\rm int}_{\jph,k}\big)_r-
(\xbar{\rho v})^{\,\rm int}_{\jph,k}\big]\right),\\
&\delta^\Gamma_{\jph,k}={\rm minmod}\left(-a^-_{\jph,k}\big[\,\xbar\Gamma^{\,\rm int}_{\jph,k}-
\big(\Gamma^{\rm int}_{\jph,k}\big)_\ell\big],\,a^+_{\jph,k}\big[\big(\Gamma^{\rm int}_{\jph,k}\big)_r-
\,\xbar\Gamma^{\,\rm int}_{\jph,k}\big]\right),\\
&\delta^\Pi_{\jph,k}={\rm minmod}\left(-a^-_{\jph,k}\big[\,\xbar\Pi^{\,\rm int}_{\jph,k}-\big(\Pi^{\rm int}_{\jph,k}\big)_\ell\big],\,
a^+_{\jph,k}\big[\big(\Pi^{\rm int}_{\jph,k}\big)_r-\,\xbar\Pi^{\,\rm int}_{\jph,k}\big]\right).
\end{aligned}
\end{equation*}
Here, the values $\big(\rho^{\rm int}_{\jph,k}\big)_\ell$, $\big(\rho^{\rm int}_{\jph,k}\big)_r$,
$\big((\rho v)^{\rm int}_{\jph,k}\big)_\ell$, $\big((\rho v)^{\rm int}_{\jph,k}\big)_r$, $\big(\Gamma^{\rm int}_{\jph,k}\big)_\ell$,
$\big(\Gamma^{\rm int}_{\jph,k}\big)_r$, $\big(\Pi^{\rm int}_{\jph,k}\big)_\ell$, and $\big(\Pi^{\rm int}_{\jph,k}\big)_r$, are obtained
using the Taylor expansions as it was done in \eref{2.2.12}.

We then proceed as in the remaining part of \S\ref{sec2.2.1} and complete the derivation of the fully discrete scheme (not shown here for
the sake of brevity), and after this, we pass to the semi-discrete limit and end up with the LD PCCU flux \eref{3.5} with the following
``built-in'' anti-diffusion term:
\begin{equation}
\bm q_{\jph,k}=
\left(q^\rho_{\jph,k},u^*_{\jph,k}q^\rho_{\jph,k},q^{\rho v}_{\jph,k},q^E_{\jph,k},q^\Gamma_{\jph,k},q^\Pi_{\jph,k}\right)^\top.
\label{3.20}
\end{equation}
Here,
$$
\begin{aligned}
&\bm U^*_{\jph,k}=\frac{a^+_{\jph,k}\bm U^+_{\jph,k}-a^-_{\jph,k}\bm U^-_{\jph,k}-\big(\bm K^+_{\jph,k}-\bm K^-_{\jph,k}\big)}
{a^+_{\jph,k}-a^-_{\jph,k}},\quad u^*_{\jph,k}=\frac{(\rho u)^*_{\jph,k}}{\rho^*_{\jph,k}},\\
&q^\rho_{\jph,k}={\rm minmod}\left(-a^-_{\jph,k}\big(\rho^*_{\jph,k}-\rho^-_{\jph,k}\big),
a^+_{\jph,k}\big(\rho^+_{\jph,k}-\rho^*_{\jph,k}\big)\right),\\
&q^{\rho v}_{\jph,k}={\rm minmod}\left(-a^-_{\jph,k}\big((\rho v)^*_{\jph,k}-(\rho v)^-_{\jph,k}\big),
a^+_{\jph,k}\big((\rho v)^+_{\jph,k}-(\rho v)^*_{\jph,k}\big)\right),\\
&q^\Gamma_{\jph,k}={\rm minmod}\left(-a^-_{\jph,k}\big(\Gamma^*_{\jph,k}-\Gamma^-_{\jph,k}\big),
a^+_{\jph,k}\big(\Gamma^+_{\jph,k}-\Gamma^*_{\jph,k}\big)\right),\\
&q^\Pi_{\jph,k}={\rm minmod}\left(-a^-_{\jph,k}\big(\Pi^*_{\jph,k}-\Pi^-_{\jph,k}\big),
a^+_{\jph,k}\big(\Pi^+_{\jph,k}-\Pi^*_{\jph,k}\big)\right),\\
&q^E_{\jph,k}=\frac{a^+_{\jph,k}a^-_{\jph,k}}{a^+_{\jph,k}-a^-_{\jph,k}}
\left\{\frac{(d^\Gamma)^-_{\jph,k}\Bigg((\rho v)^*_{\jph,k}+\dfrac{q^{\rho v}_{\jph,k}}{a^+_{\jph,k}}\Bigg)^2}
{2\Bigg(\rho^*_{\jph,k}+\dfrac{q^\rho_{\jph,k}}{a^+_{\jph,k}}\Bigg)}-
\frac{(d^\Gamma)^+_{\jph,k}\Bigg((\rho v)^*_{\jph,k}+\dfrac{q^{\rho v}_{\jph,k}}{a^-_{\jph,k}}\Bigg)^2}
{2\Bigg(\rho^*_{\jph,k}+\dfrac{q^\rho_{\jph,k}}{a^-_{\jph,k}}\Bigg)}\right\}\\
&\hspace*{1.05cm}+\frac{\big(u^*_{\jph,k}\big)^2}{2}q^\rho_{\jph,k}+
\frac{q^\Gamma_{\jph,k}}{\Gamma^*_{\jph,k}}\bigg[E^*_{\jph,k}-\frac{\big((\rho u)^*_{\jph,k}\big)^2}{2\rho^*_{\jph,k}}-\Pi^*_{\jph,k}\bigg]+
q^\Pi_{\jph,k},
\end{aligned}
$$
where
$$
(d^\Gamma)^\pm_{\jph,k}=1-\frac{q^\Gamma_{\jph,k}}{a^\pm_{\jph,k}\Gamma^*_{\jph,k}}.
$$
Similarly, the ``built-in'' anti-diffusion term in the LD PCCU flux \eref{3.6} is
\begin{equation}
\bm q_{j,\kph}=
\left(q^\rho_{j,\kph},q^{\rho u}_{j,\kph},v^*_{j,\kph}q^\rho_{j,\kph},q^E_{j,\kph},q^\Gamma_{j,\kph},q^\Pi_{j,\kph}\right)^\top.
\label{3.20f}
\end{equation}
with
\allowdisplaybreaks
\begin{align*}
&\bm U^*_{j,\kph}=\frac{b^+_{j,\kph}\bm U^+_{j,\kph}-b^-_{j,\kph}\bm U^-_{j,\kph}-\big(\bm L^+_{j,\kph}-\bm L^-_{j,\kph}\big)}
{b^+_{j,\kph}-b^-_{j,\kph}},\quad v^*_{j,\kph}=\frac{(\rho v)^*_{j,\kph}}{\rho^*_{j,\kph}},\\
&q^\rho_{j,\kph}={\rm minmod}\left(-b^-_{j,\kph}\big(\rho^*_{j,\kph}-\rho^-_{j,\kph}\big),
b^+_{j,\kph}\big(\rho^+_{j,\kph}-\rho^*_{j,\kph}\big)\right),\\
&q^{\rho u}_{j,\kph}={\rm minmod}\left(-b^-_{j,\kph}\big((\rho u)^*_{j,\kph}-(\rho u)^-_{j,\kph}\big),
b^+_{j,\kph}\big((\rho u)^+_{j,\kph}-(\rho u)^*_{j,\kph}\big)\right),\\
&q^\Gamma_{j,\kph}={\rm minmod}\left(-b^-_{j,\kph}\big(\Gamma^*_{j,\kph}-\Gamma^-_{j,\kph}\big),
b^+_{j,\kph}\big(\Gamma^+_{j,\kph}-\Gamma^*_{j,\kph}\big)\right),\\
&q^\Pi_{j,\kph}={\rm minmod}\left(-b^-_{j,\kph}\big(\Pi^*_{j,\kph}-\Pi^-_{j,\kph}\big),
b^+_{j,\kph}\big(\Pi^+_{j,\kph}-\Pi^*_{j,\kph}\big)\right),\\
&q^E_{j,\kph}=\frac{b^+_{j,\kph}b^-_{j,\kph}}{b^+_{j,\kph}-b^-_{j,\kph}}
\left\{\frac{(d^\Gamma)^-_{j,\kph}\Bigg((\rho u)^*_{j,\kph}+\dfrac{q^{\rho u}_{j,\kph}}{b^+_{j,\kph}}\Bigg)^2}
{2\Bigg(\rho^*_{j,\kph}+\dfrac{q^\rho_{j,\kph}}{b^+_{j,\kph}}\Bigg)}-
\frac{(d^\Gamma)^+_{j,\kph}\Bigg((\rho u)^*_{j,\kph}+\dfrac{q^{\rho u}_{j,\kph}}{b^-_{j,\kph}}\Bigg)^2}
{2\Bigg(\rho^*_{j,\kph}+\dfrac{q^\rho_{j,\kph}}{b^-_{j,\kph}}\Bigg)}\right\}\\
&\hspace*{1.05cm}+\frac{\big(v^*_{j,\kph}\big)^2}{2}q^\rho_{j,\kph}+
\frac{q^\Gamma_{j,\kph}}{\Gamma^*_{j,\kph}}\bigg[E^*_{j,\kph}-\frac{\big((\rho v)^*_{j,\kph}\big)^2}{2\rho^*_{j,\kph}}-\Pi^*_{j,\kph}\bigg]+
q^\Pi_{j,\kph},
\end{align*}
where
$$
(d^\Gamma)^\pm_{j,\kph}=1-\frac{q^\Gamma_{j,\kph}}{b^\pm_{j,\kph}\Gamma^*_{j,\kph}}.
$$
\begin{rmk}
As in \cite{KX_22}, the computation of numerical fluxes in \eref{3.5} should be desingularized to avoid division by zero or very small
numbers. First, if
$a^+_{\jph,k}<\varepsilon_0$ and $a^-_{\jph,k}>-\varepsilon_0$ for a small positive $\varepsilon_0$, we replace the flux
$\bm{{\cal K}}_{\jph,k}$ with
\begin{equation*}
\bm{{\cal K}}_{\jph,k}=\frac{\bm K\big(\bm U^-_{\jph,k}\big)+\bm K\big(\bm U^+_{\jph,k}\big)}{2}.
\end{equation*}
Similarly, if $b^+_{j,\kph}<\varepsilon_0$ and $b^-_{j,\kph}>-\varepsilon_0$, we replace the flux $\bm{{\cal L}}_{j,\kph}$ with
\begin{equation*}
\bm{{\cal L}}_{j,\kph}=\frac{\bm L\big(\bm U^-_{j,\kph}\big)+\bm L\big(\bm U^+_{j,\kph}\big)}{2}.
\end{equation*}
In addition, the computation of the energy numerical fluxes have to be desingularized even in the case when only one of the local speeds is
very small. In particular,
$$
\begin{aligned}
&\mbox{if}~a^+_{\jph,k}<\varepsilon_0~{\rm but}~a^-_{\jph,k}<-\varepsilon_0,&&\mbox{we take}~~{\cal K}^{(4)}_{\jph,k}=
u^-_{\jph,k}\big(E^-_{\jph,k}+p^-_{\jph,k}\big);\\
&\mbox{if}~a^-_{\jph,k}>-\varepsilon_0~{\rm but}~a^+_{\jph,k}>\varepsilon_0,&&\mbox{we take}~~{\cal K}^{(4)}_{\jph,k}=
u^+_{\jph,k}\big(E^+_{\jph,k}+p^+_{\jph,k}\big);\\
&\mbox{if}~b^+_{j,\kph}<\varepsilon_0~{\rm but}~b^-_{j,\kph}<-\varepsilon_0,&&\mbox{we take}~~{\cal L}^{(4)}_{j,\kph}=
v^-_{j,\kph}\big(E^-_{j,\kph}+p^-_{j,\kph}\big);\\
&\mbox{if}~b^-_{j,\kph}>-\varepsilon_0~{\rm but}~b^+_{j,\kph}>\varepsilon_0,&&\mbox{we take}~~{\cal L}^{(4)}_{j,\kph}=
v^+_{j,\kph}\big(E^+_{j,\kph}+p^+_{j,\kph}\big).
\end{aligned}
$$
As in the 1-D case, we take $\varepsilon_0=10^{-12}$ in all of the numerical examples.
\end{rmk}

\subsection{Flux Globalization Based LD Ai-WENO PCCU Scheme}
In this section, we extend the 2-D second-order flux globalization based LD PCCU schemes from \S\ref{sec3.2} to the fifth order of accuracy
within the A-WENO framework.

The semi-discrete fifth-order LD Ai-WENO PCCU scheme for the 2-D quasi-conservative system \eref{3.1}--\eref{3.2} reads as
\begin{equation}
\frac{{\rm d}}{{\rm d}t}\mU_{j,k}=-\frac{\bmH^x_{\jph,k}-\bmH^x_{\jmh,k}}{\dx}-\frac{\bmH^y_{j,\kph}-\bmH^y_{j,\kmh}}{\dy},
\label{3.13}
\end{equation}
where the fifth-order numerical fluxes $\bmH^x_{\jph,k}$ and $\bmH^y_{j,\kph}$ are defined by
\begin{equation}
\begin{aligned}
&\bmH^x_{\jph,k}=\bm{{\cal K}}_{\jph,k}-\frac{\dx}{24}(\bm K_{xx})_{\jph,k}+\frac{7(\dx)^3}{5760}(\bm K_{xxxx})_{\jph,k},\\
&\bmH^y_{j,\kph}=\bm{{\cal L}}_{j,\kph}-\frac{\dy}{24}(\bm L_{yy})_{j,\kph}+\frac{7(\dy)^3}{5760}(\bm L_{yyyy})_{j,\kph}.
\end{aligned}
\label{3.14}
\end{equation}
Here, $\bm{{\cal K}}_{\jph,k}$ and $\bm{{\cal L}}_{j,\kph}$ are the finite-volume fluxes \eref{3.5}, and $(\bm K_{xx})_{\jph,k}$,
$({\bm K_{xxxx}})_{\jph,k}$, $(\bm L_{yy})_{j,\kph}$, and $({\bm L_{yyyy}})_{j,\kph}$ are approximations of the second- and fourth-order
spatial derivatives of $\bm K$ at $(x,y)=(x_\jph,y_k)$ and $\bm L$ at $(x,y)=(x_j,y_\kph)$, respectively. We compute these quantities using
the finite-difference approximations analogous to \eref{2.57}:
\begin{equation*}
\begin{aligned}
&(\mK_{xx})_{\jph,k}=\frac{1}{12(\dx)^2}\left[-\bm{{\cal K}}_{j-\frac{3}{2},k}+16\bm{{\cal K}}_{\jmh,k}-30\bm{{\cal K}}_{\jph,k}+
16\bm{{\cal K}}_{j+\frac{3}{2},k}-\bm{{\cal K}}_{j+\frac{5}{2},k}\right],\\
&(\mK_{xxxx})_{\jph,k}=\frac{1}{(\dx)^4}\left[\bm{{\cal K}}_{j-\frac{3}{2},k}-4\bm{{\cal K}}_{\jmh,k}+6\bm{{\cal K}}_{\jph,k}-
4\bm{{\cal K}}_{j+\frac{3}{2},k}+\bm{{\cal K}}_{j+\frac{5}{2},k}\right],\\
&(\mL_{yy})_{j,\kph}=\frac{1}{12(\dy)^2}\left[-\bm{{\cal L}}_{j,k-\frac{3}{2}}+16\bm{{\cal L}}_{j,\kmh}-30\bm{{\cal L}}_{j,\kph}+
16\bm{{\cal L}}_{j,k+\frac{3}{2}}-\bm{{\cal L}}_{j,k+\frac{5}{2}}\right],\\
&(\mL_{yyyy})_{j,\kph}=\frac{1}{(\dy)^4}\left[\bm{{\cal L}}_{j,k-\frac{3}{2}}-4\bm{{\cal L}}_{j,\kmh}+6\bm{{\cal L}}_{j,\kph}-
4\bm{{\cal L}}_{j,k+\frac{3}{2}}+\bm{{\cal L}}_{j,k+\frac{5}{2}}\right].
\end{aligned}
\end{equation*}

As in the 1-D case, the resulting semi-discrete scheme \eref{3.13}--\eref{3.14} is fifth-order accurate provided the point values
$\bm U^\pm_{\jph,k}$ and $\bm U^\pm_{j,\kph}$ are calculated using a fifth-order interpolation. To this end, we apply the recently proposed
fifth-order Ai-WENO-Z interpolation \cite{DLWW22,LLWDW23,WD22}, which can be performed in the $x$- and $y$-directions similarly to the 1-D
case discussed in Appendix \ref{appa}; we omit the details for the sake of brevity.
\begin{rmk}
As in the 1-D case, one needs to apply the Ai-WENO-Z interpolation procedure in the local characteristic variables to reduce the magnitude
of the numerical oscillations. In Appendix \ref{appc}, we provide a detailed explanation on how to apply the LCD to the 2-D system
\eref{3.1}--\eref{3.3}.
\end{rmk}

\section{Numerical Examples}\label{sec4}
In this section, we apply the developed schemes to several 1-D and 2-D numerical examples and compare the performance of the second-order
flux globalization based PCCU, the second-order flux globalization based LD PCCU, and the fifth-order flux globalization based LD Ai-WENO
PCCU schemes, which will be referred to as the {\em PCCU}, {\em LD PCCU}, and {\em Ai-WENO} schemes, respectively.

In all of the numerical examples, we have solved the ODE systems \eref{2.1.4}, \eref{2.53}, \eref{3.4}, and \eref{3.13} using the
three-stage third-order strong stability preserving (SSP) Runge-Kutta method (see, e.g.,\cite{Gottlieb11,Gottlieb12}) and used the CFL
number 0.45.

\subsection{One-Dimensional Examples}
\subsubsection*{Example 1---``Shock-Bubble'' Interaction}
In the first example, we consider the ``shock-bubble'' interaction problem, which is a two-fluid modification of a single-fluid example from
\cite{KX_22}. The initial data are given by
\begin{equation*}
(\rho,u,p;\gamma,\pi_\infty)(x,0)=\begin{cases}(13.1538,0,1;5/3,0),&|x|<0.25,\\(1.3333,-0.3535,1.5;1.4,0),&x>0.75,\\
(1,0,1;1.4,0),&\mbox{otherwise,}\end{cases}
\end{equation*}
which correspond to a left-moving shock, initially located at $x=0.75$, and a resting ``bubble'' with a radius of 0.25, initially located at
the origin. These inital data are considered in the computational domain $[-1,2]$ subject to the solid wall boundary conditions imposed at
both $x=-1$ and $x=2$.

\begin{figure}[ht!]
\centerline{\includegraphics[trim=1.2cm 0.4cm 0.9cm 0.4cm, clip, width=5.5cm]{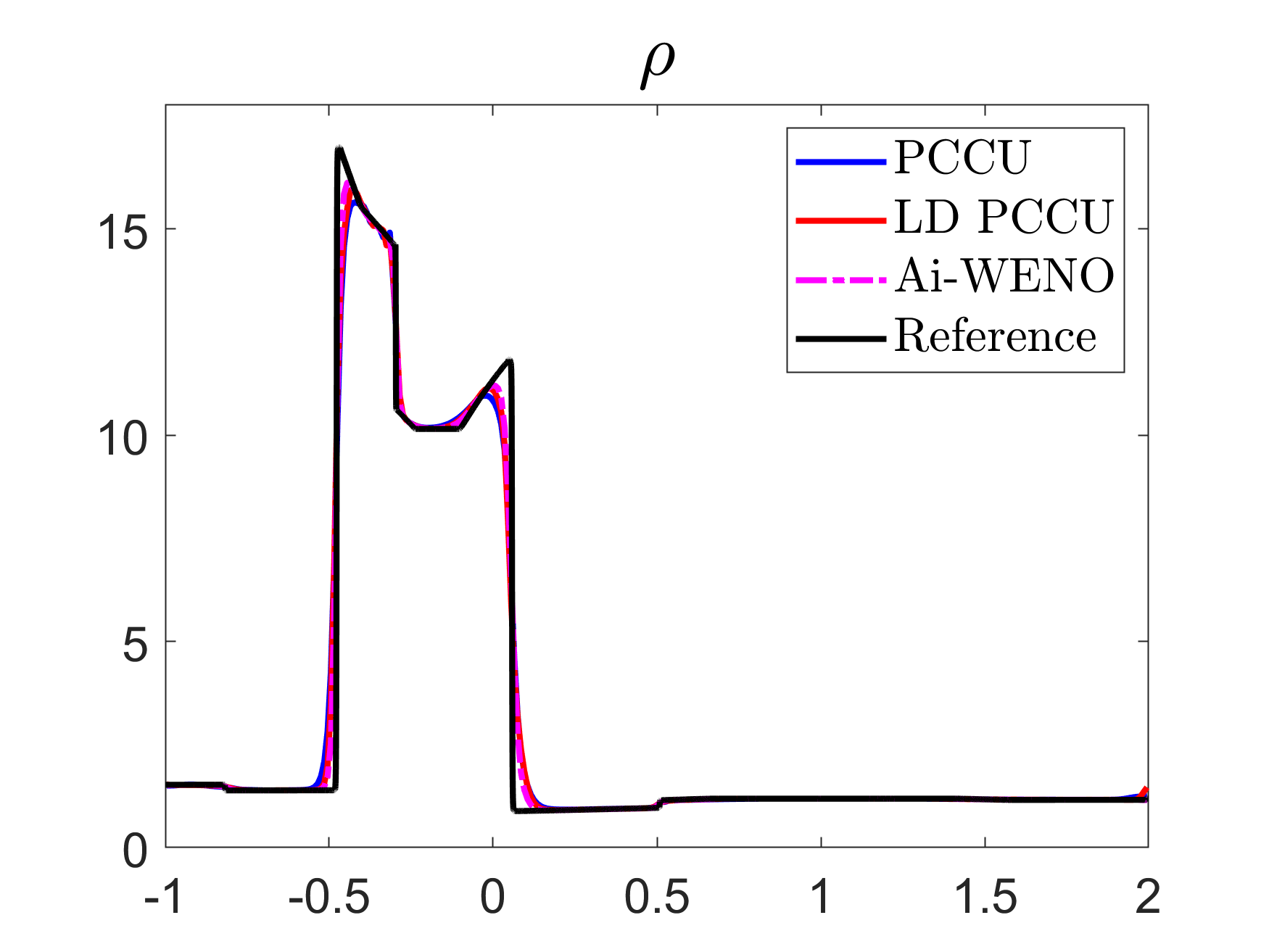}\hspace{1.0cm}
            \includegraphics[trim=1.2cm 0.4cm 0.9cm 0.4cm, clip, width=5.5cm]{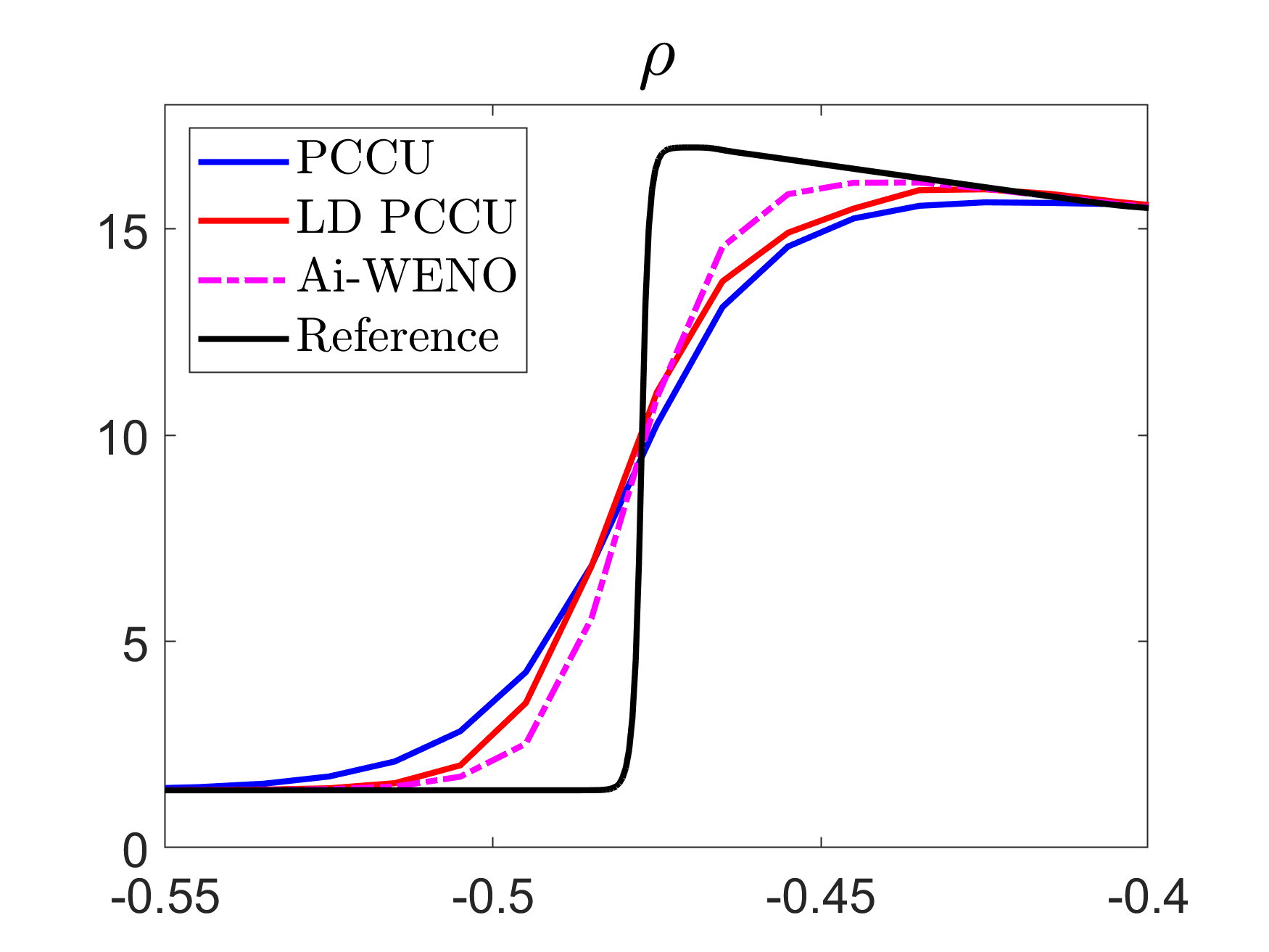}}
\caption{\sf Example 1: Density $\rho$ (left) and zoom at $x\in[-0.55,-0.4]$ (right) .\label{fig1}}
\end{figure}

\begin{figure}[ht!]
\centerline{\includegraphics[trim=0.8cm 0.4cm 1.3cm 0.4cm, clip, width=5.5cm]{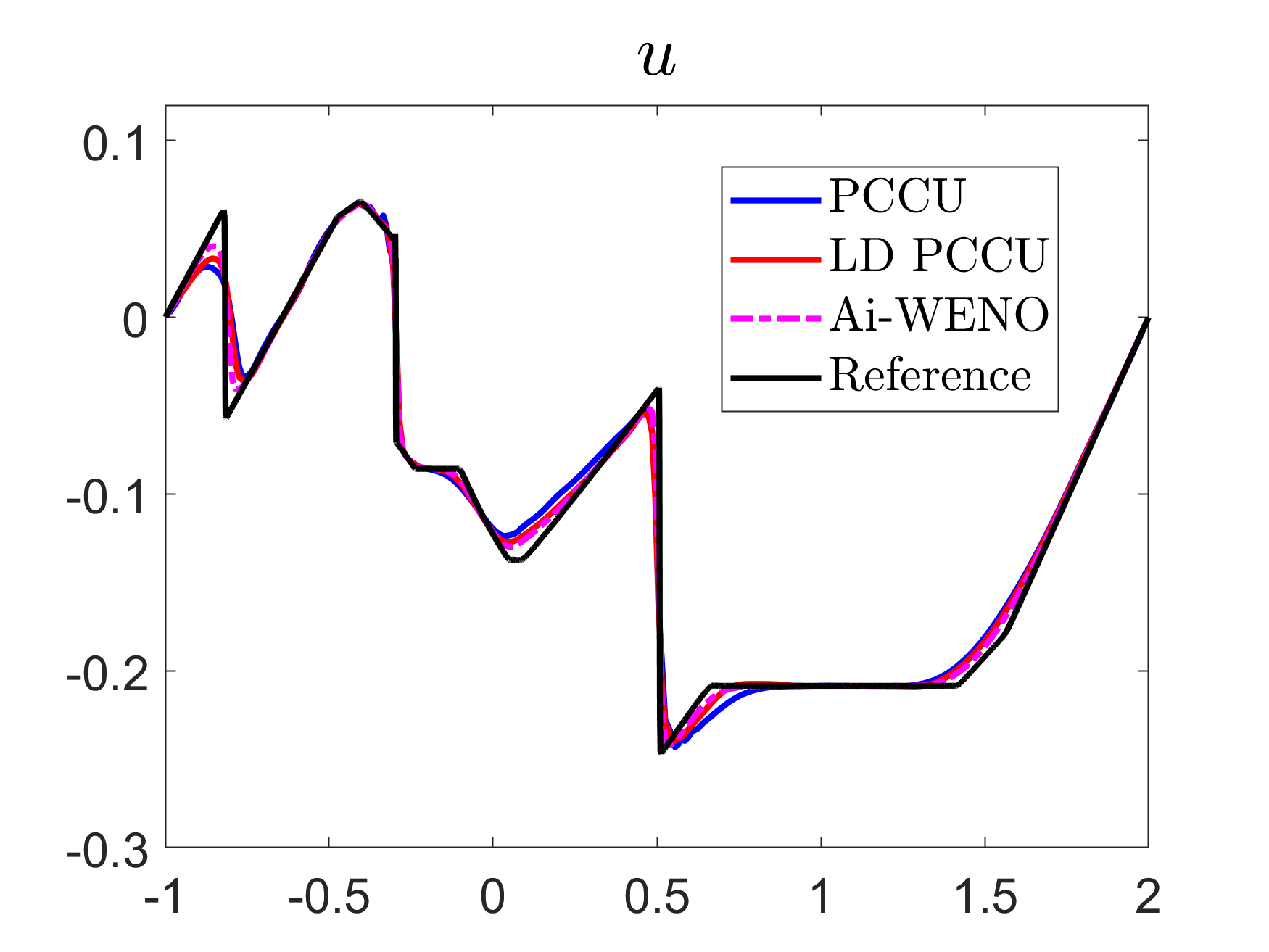}\hspace{1.0cm}
            \includegraphics[trim=0.8cm 0.4cm 1.3cm 0.4cm, clip, width=5.5cm]{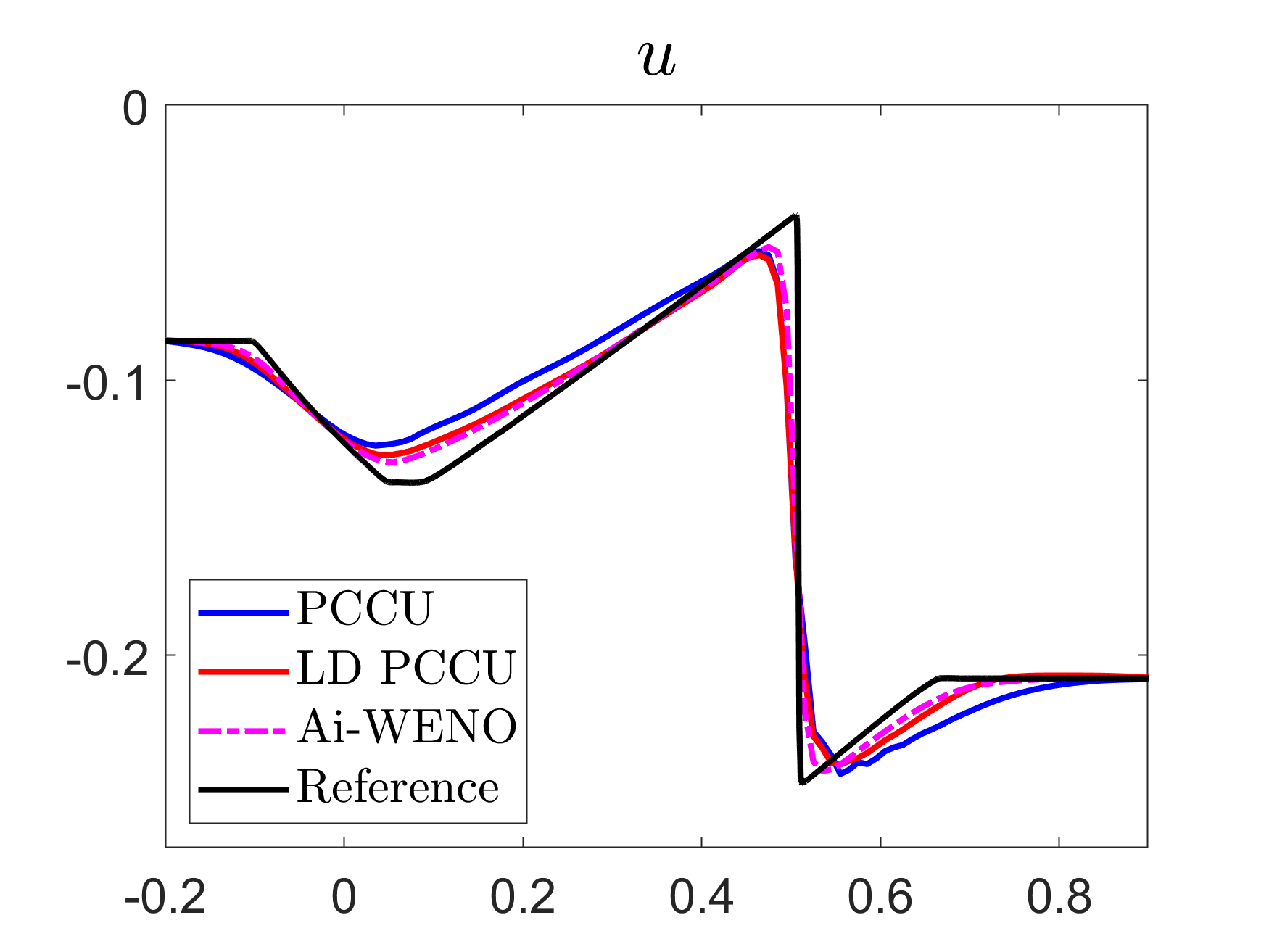}}
\caption{\sf Example 1: Velocity $u$ (left) and zoom at $x\in[-0.2,0.9]$ (right).\label{fig1a}}
\end{figure}

We apply the studied PCCU, LD PCCU, and Ai-WENO schemes to this initial-boundary value problem and compute its numerical solutions until the
final time $t=3$ on a uniform mesh with $\dx=1/100$. The obtained densities and velocities are presented in Figures \ref{fig1} and
\ref{fig1a} together with the reference solution computed by the PCCU scheme on a much finer mesh with $\dx=1/2000$. As one can see, the LD
PCCU scheme achieves sharper resolution than the PCCU one, and the use of a fifth-order Ai-WENO scheme further enhances the resolution.

\subsubsection*{Example 2---Water-Air ``Shock-Bubble'' Interaction}
In the second 1-D example, which is a 1-D modification of an example from \cite{CZL20}, we consider a gas-liquid multifluid system, where
the liquid component is modeled by the EOS \eref{2.0.1} with $\pi_\infty\gg1$. The initial conditions,
\begin{equation*}
(\rho,u,p;\gamma,\pi_\infty)=\begin{cases}(0.05,0,1;1.4,0),&|x-6|<3,\\(1.325,-68.525,19153;4.4,6000),&x>11.4,\\
(1,0,1;4.4,6000),&\mbox{otherwise,}\end{cases}
\end{equation*}
correspond to the left-moving shock, initially located at $x=11.4$, and a resting air ``bubble'' with a radius 3, initially located at
$x=6$. The initial conditions are prescribed in the computational domain $[0,18]$ subject to the free boundary conditions.

\begin{figure}[ht!]
\centerline{\includegraphics[trim=0.9cm 0.4cm 0.7cm 0.4cm, clip, width=5.5cm]{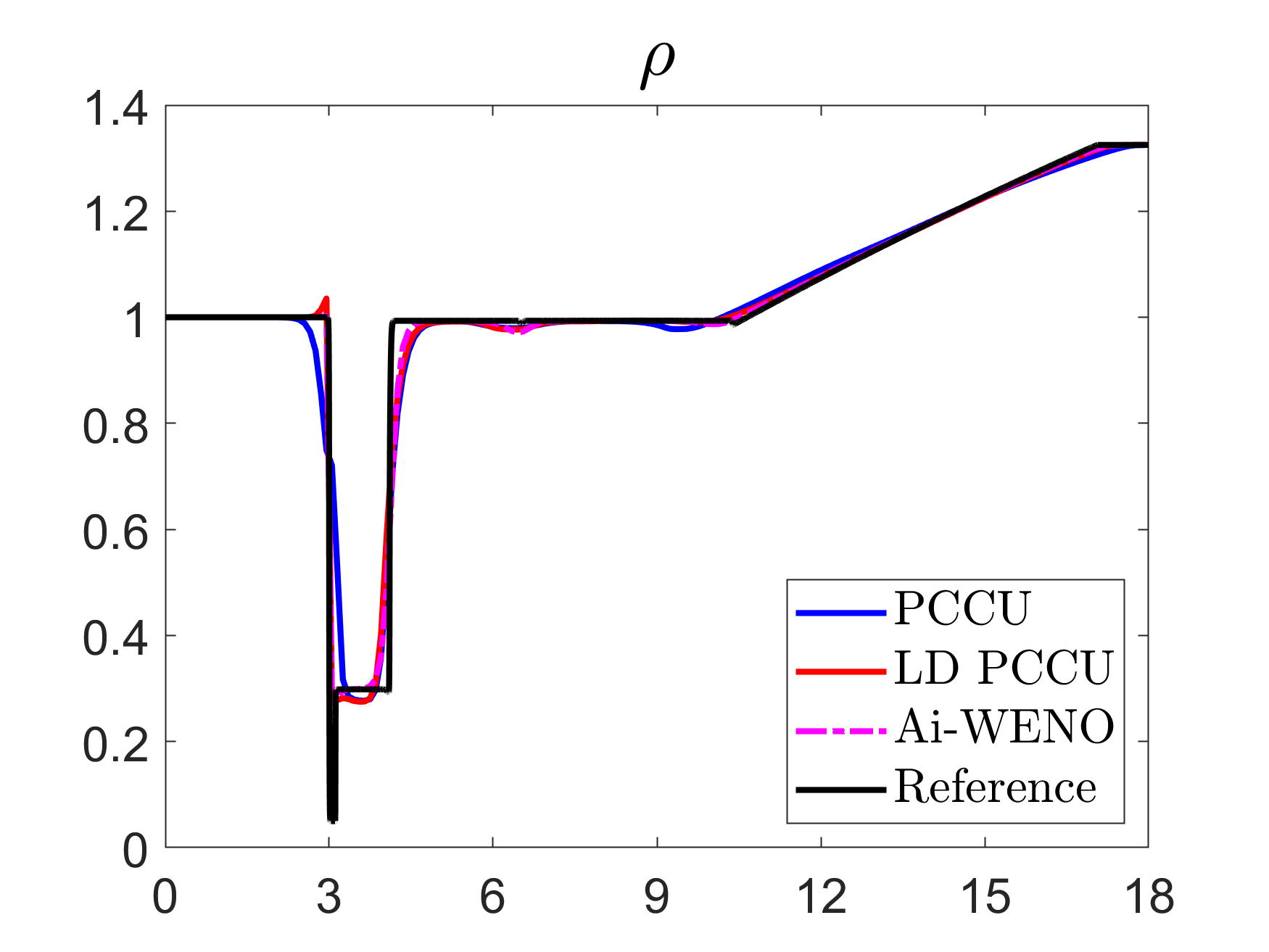}\hspace*{1.0cm}
            \includegraphics[trim=0.9cm 0.4cm 0.7cm 0.4cm, clip, width=5.5cm]{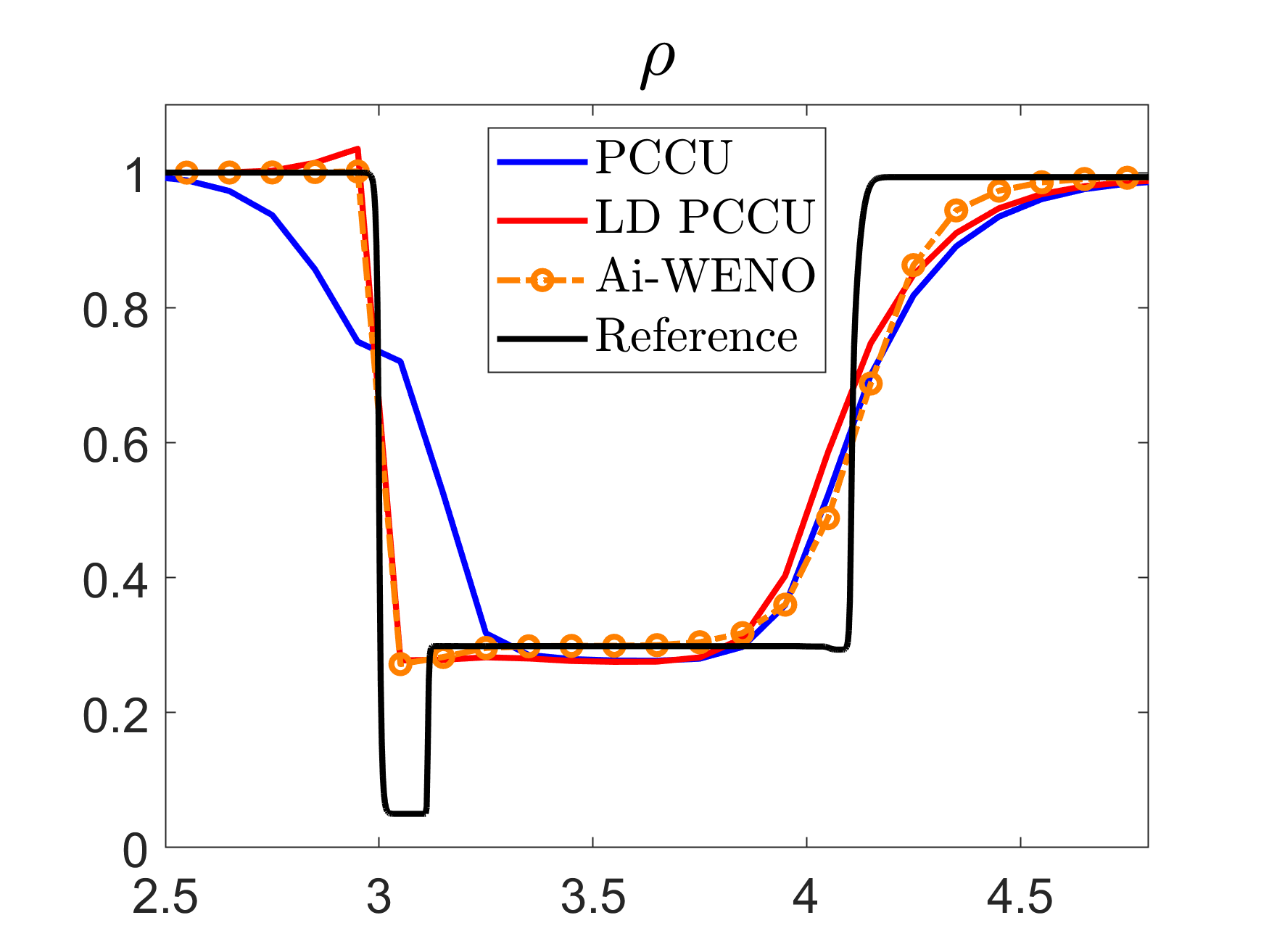}}
\caption{\sf Example 2: Density $\rho$ (left) and zoom at $x\in[2.5,4.8]$ (right).\label{fig2}}
\end{figure}

We compute the numerical solution until the final time $t=0.045$ on a uniform mesh with $\dx=1/10$ by the PCCU, LD PCCU, and Ai-WENO schemes
and plot the obtained density in Figure \ref{fig2} together with the reference solution computed by the PCCU scheme on a much finer mesh
with $\dx=1/400$. As one can observe, the LD PCCU solution achieves sharper resolution (especially of the contact wave located at about
$x=3$) compared with its PCCU counterpart, but it produces certain oscillations. The Ai-WENO scheme, on the other hand, achieves even higher
resolution and the obtained fifth-order results are oscillation-free. We believe that this is thanks to the fact that the Ai-WENO-Z
interpolation is performed using the LCD.

\subsubsection*{Example 3---Water-Air Model With a Very Stiff Equation of State}
In the last 1-D example taken from \cite{AK_01a,chertock7}, we consider another gas-liquid multifluid system with the water component
modeled using even stiffer EOS than the one used in Example 2. The initial conditions that correspond to a severe water-air shock tube
problem,
\begin{equation*}	
(\rho,u,p;\gamma,\pi_\infty)=\begin{cases}(1000,0,10^9;4.4,6\cdot10^8),&x<0.7,\\(50,0,10^5;1.4,0),&x>0.7,\end{cases}
\end{equation*}	
are prescribed in the computational domain $[0,1]$ subject to the free boundary conditions.

\begin{figure}[ht!]
\centerline{\includegraphics[trim=0.5cm 0.4cm 0.9cm 0.4cm, clip, width=5.5cm]{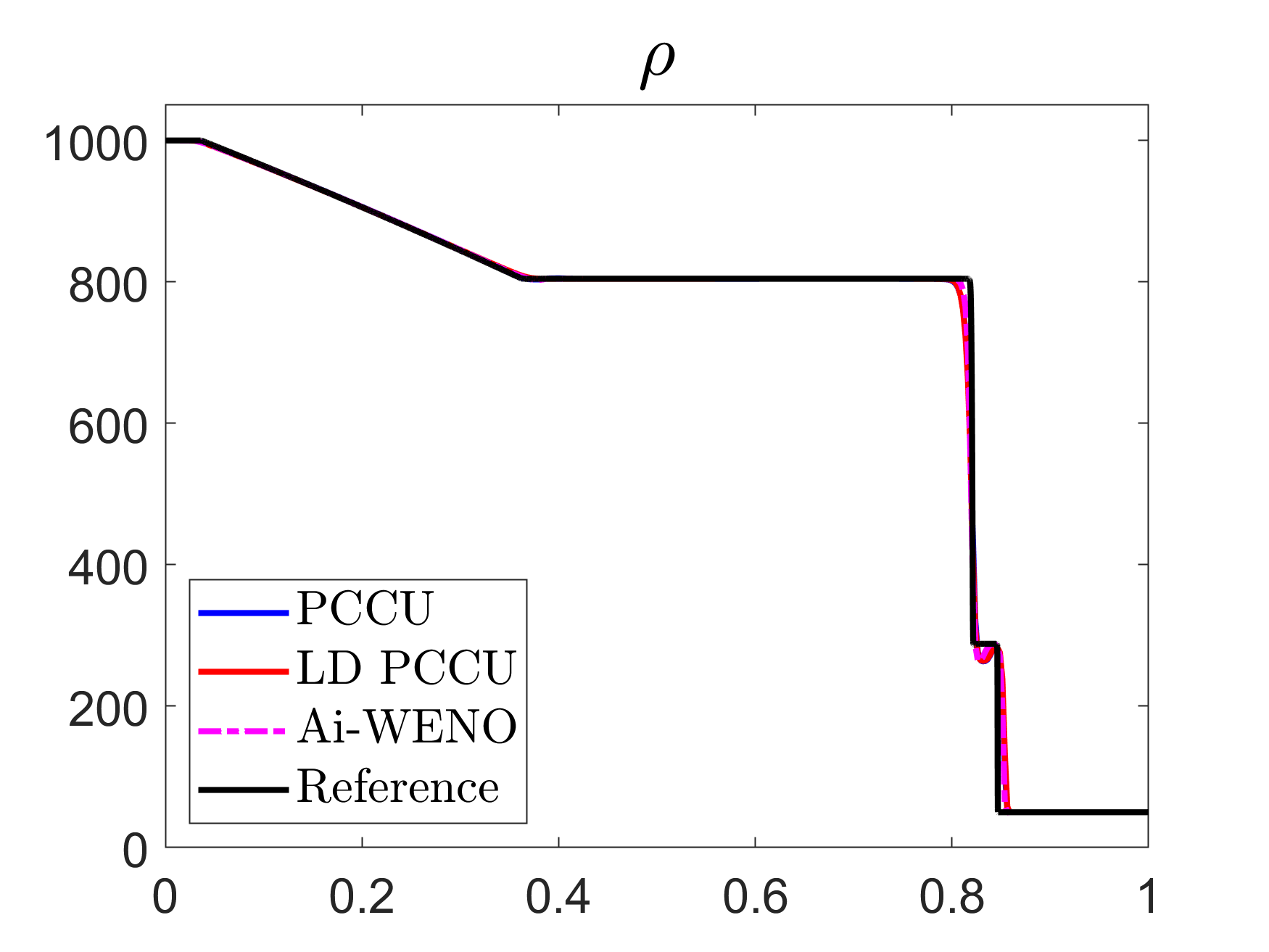}\hspace*{1.0cm}
            \includegraphics[trim=0.5cm 0.4cm 0.9cm 0.4cm, clip, width=5.5cm]{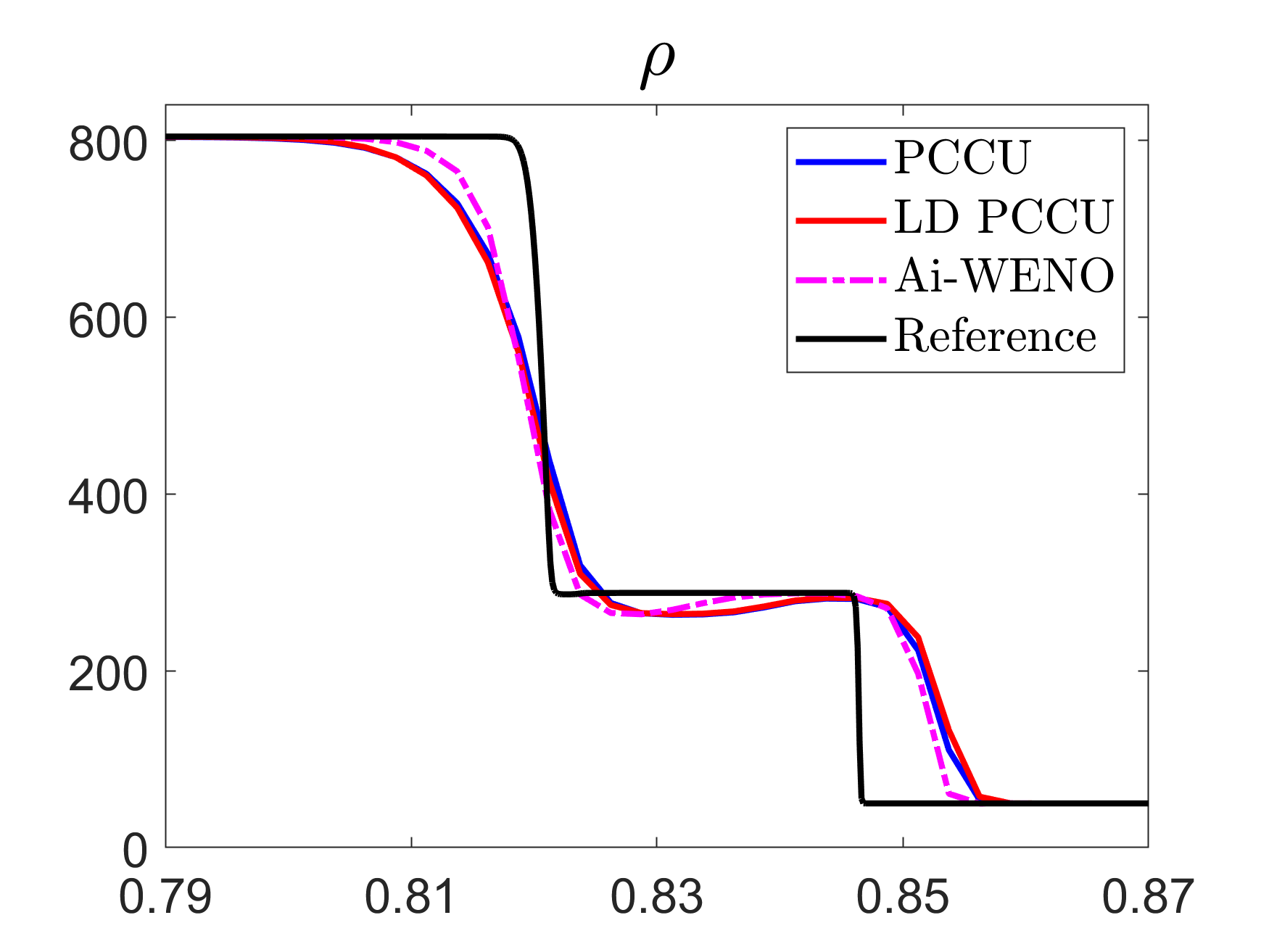}}
\caption{\sf Example 3: Density $\rho$ (left) and zoom at $[0.79, 0.87]$ (right).\label{fig3}}
\end{figure}

We compute the numerical solutions by the studied PCCU, LD PCCU, and Ai-WENO schemes until the final time $t=0.00025$ on a uniform mesh with
$\dx=1/400$. The obtained densities are shown in Figure \ref{fig3} along with the reference solution computed by the PCCU scheme on a much
finer mesh with $\dx=1/6400$. One can observe that all of the studied schemes produce non-oscillatory numerical solutions, and the Ai-WENO
solution is slightly sharper than the solutions computed by the PCCU and LD PCCU schemes.

\subsection{Two-Dimensional Examples}
In this section, we present four 2-D numerical examples. In all of them, we plot the Schlieren images of the magnitude of the density
gradient field, $|\nabla\rho|$. To this end, we have used the following shading function:
\begin{equation*}
\exp\bigg(-\frac{80|\nabla\rho|}{\max(|\nabla\rho|)}\bigg),
\end{equation*}
where the numerical derivatives of the density are computed using standard central differencing.

\subsubsection*{Example 4---Shock-Helium Bubble Interaction}
\begin{figure}[ht!]
\centerline{\includegraphics[trim=1.0cm 0.2cm 1.0cm 0.4cm, clip, width=7.0cm]{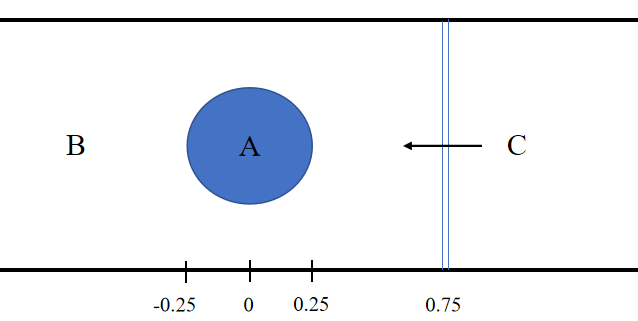}}
\caption{\sf Initial setting for the 2-D numerical examples.\label{fig44a}}
\end{figure}

In the first 2-D example taken from \cite{Quirk1996,chertock7}, a shock wave in the air hits the light resting bubble which contains helium.
We take the following initial conditions:
\begin{equation*}
(\rho,u,v,p;\gamma,p_\infty)=\begin{cases}(4/29,0,0,1;5/3,0),&\mbox{in region A},\\(1,0,0,1;1.4,0),&\mbox{in region B},\\
(4/3,-0.3535,0,1.5;1.4,0),&\mbox{in region C},\end{cases}
\end{equation*}
where regions A, B, and C are outlined in Figure \ref{fig44a}, and the computational domain is $[-3,1]\times[-0.5,0.5]$. We impose the solid
wall boundary conditions on the top and bottom and the free boundary conditions on the left and right edges of the computational domain.

\begin{figure}[ht!]
\centerline{\includegraphics[trim=2.1cm 0.4cm 2.1cm 0.2cm, clip, width=5.1cm]{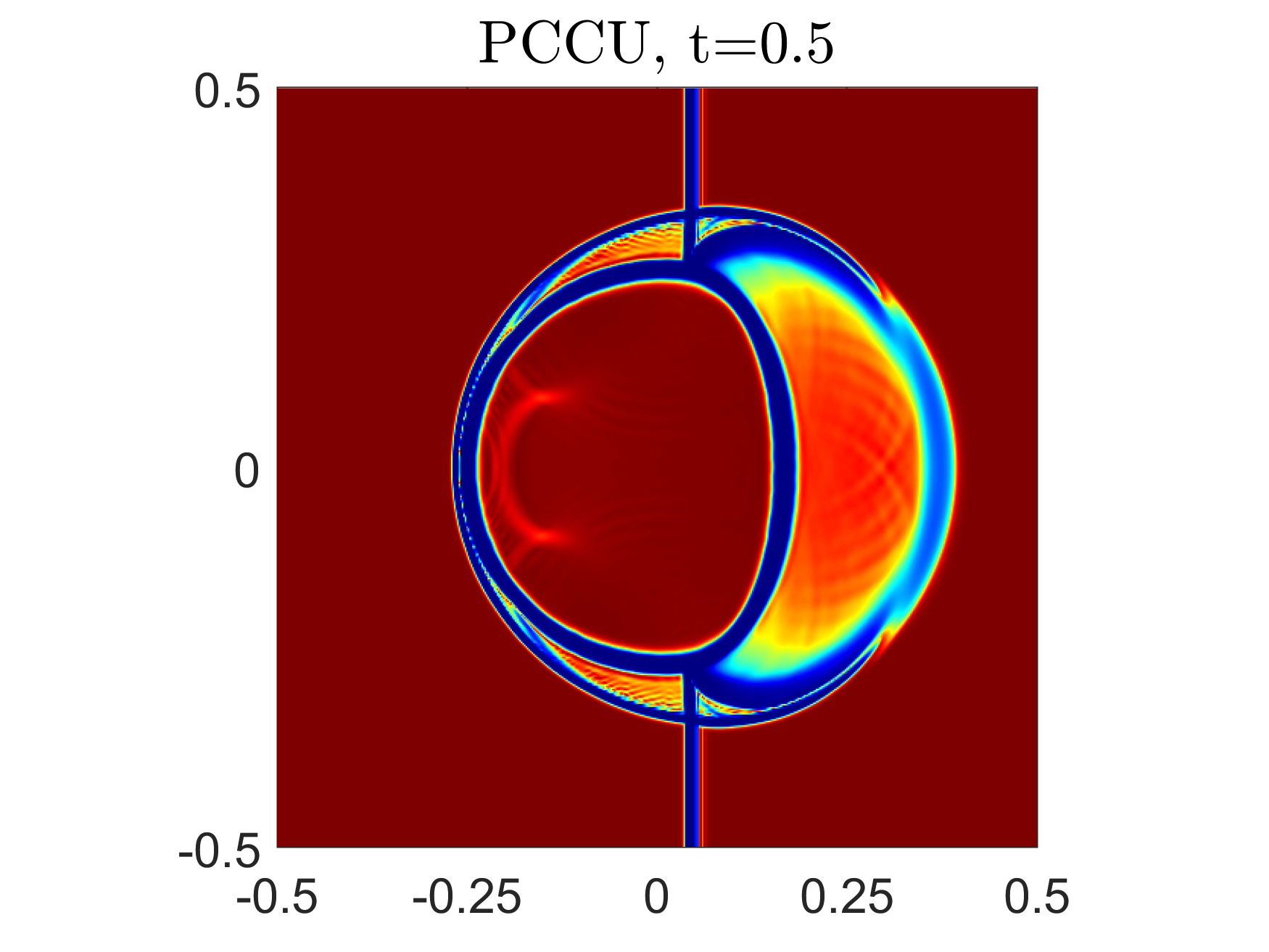}\hspace*{0.5cm}
            \includegraphics[trim=2.1cm 0.4cm 2.1cm 0.2cm, clip, width=5.1cm]{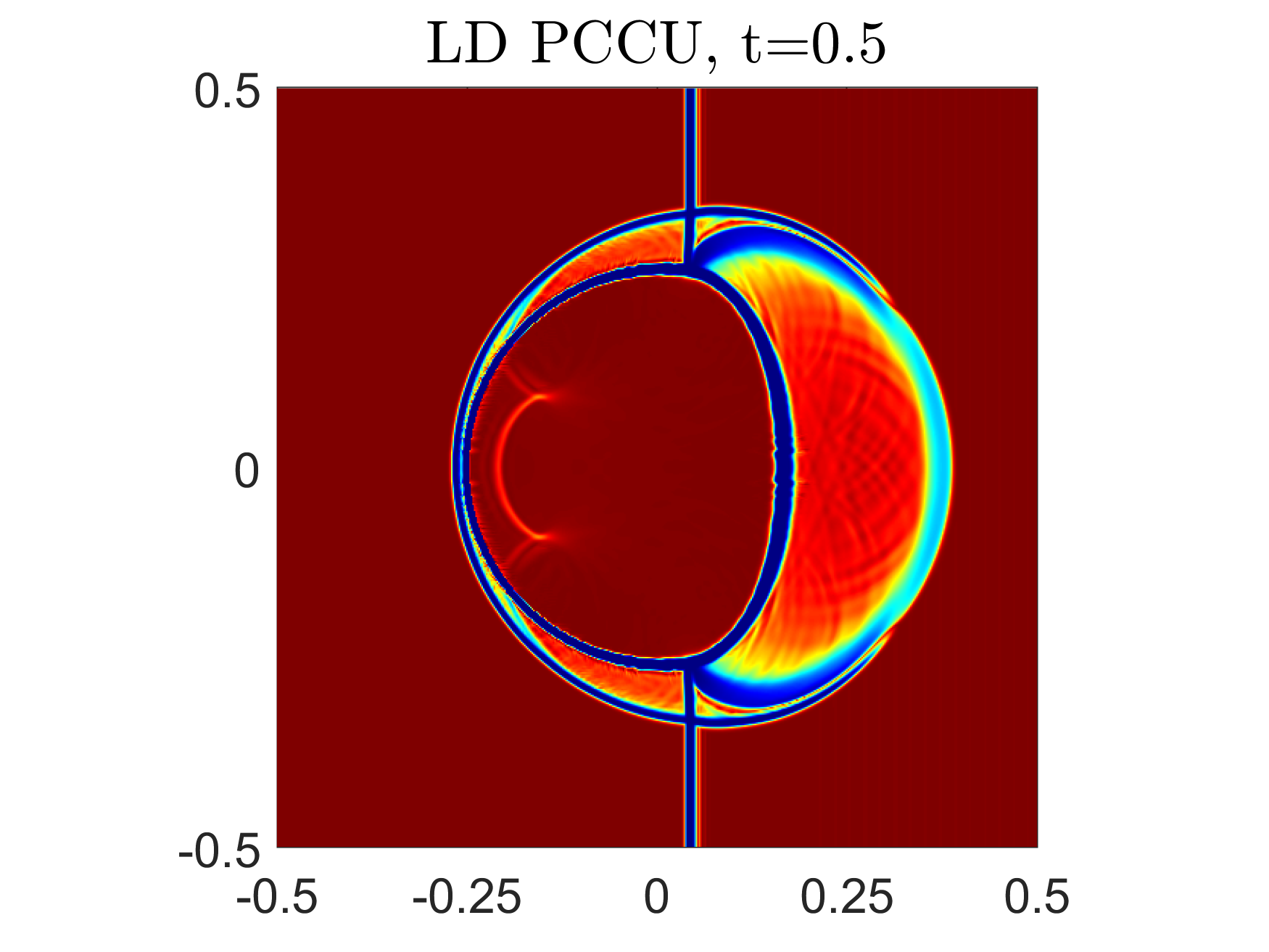}\hspace*{0.5cm}
            \includegraphics[trim=2.1cm 0.4cm 2.1cm 0.2cm, clip, width=5.1cm]{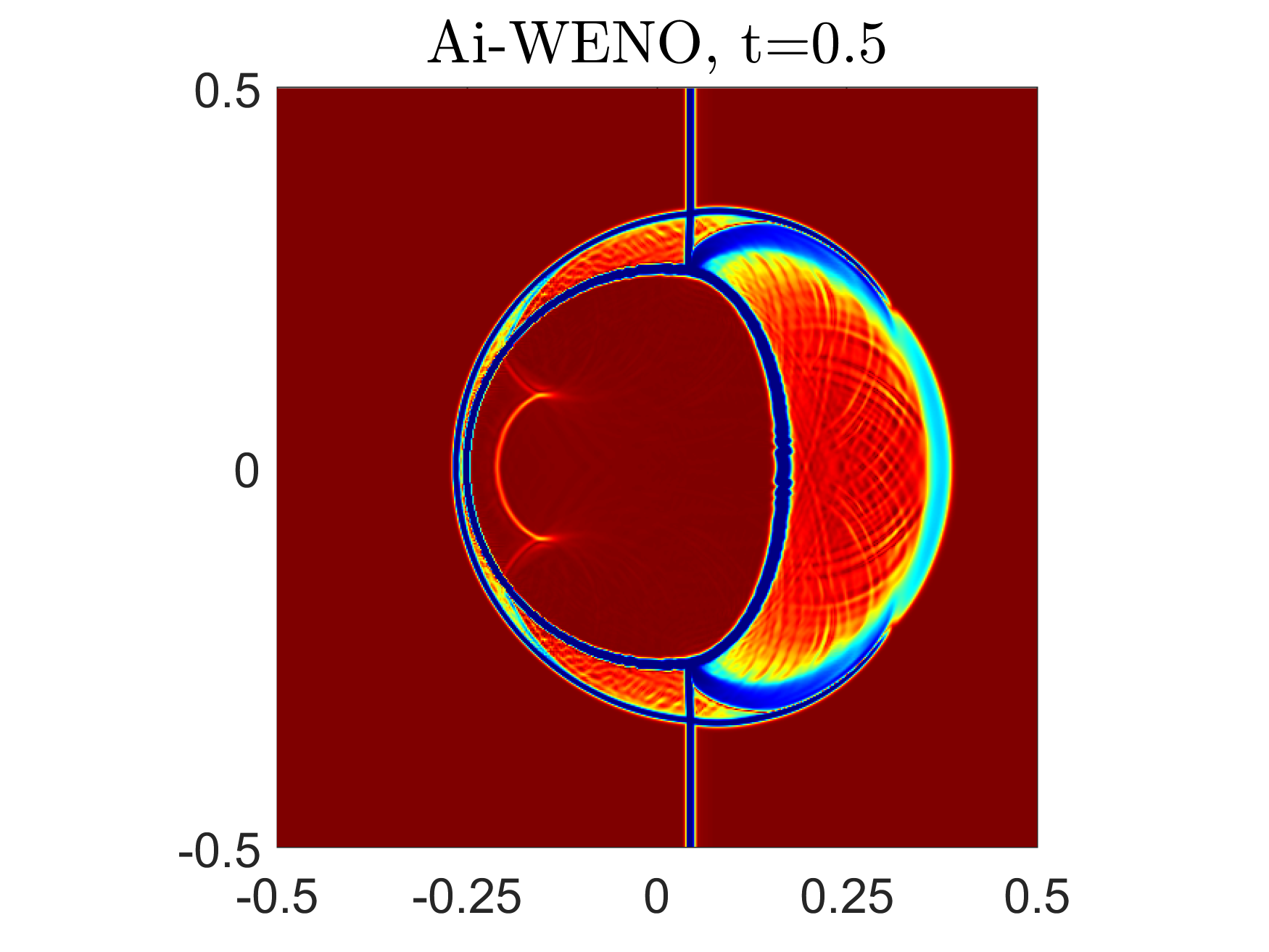}}
\vskip10pt
\centerline{\includegraphics[trim=2.1cm 0.4cm 2.1cm 0.2cm, clip, width=5.1cm]{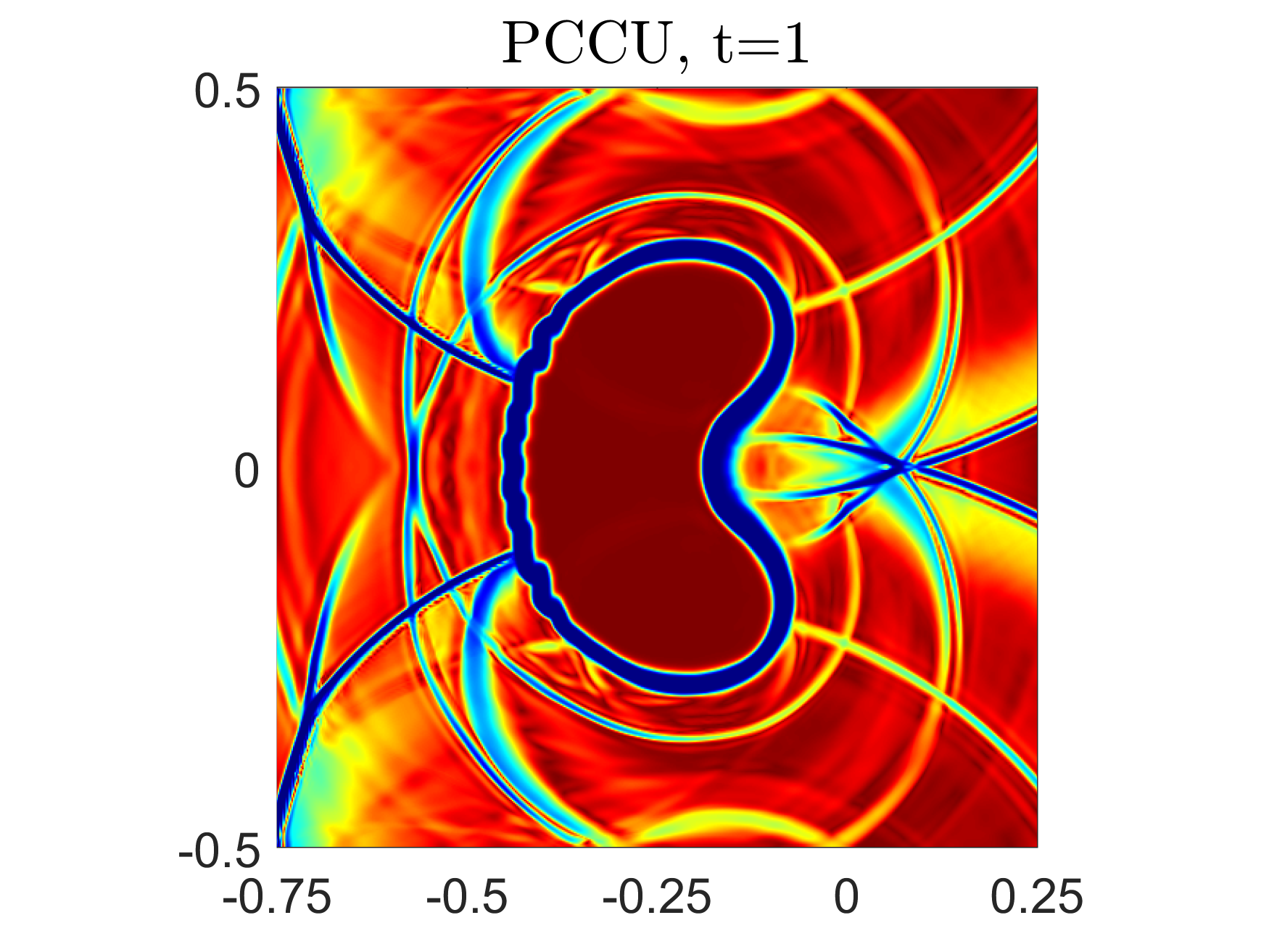}\hspace*{0.5cm}
            \includegraphics[trim=2.1cm 0.4cm 2.1cm 0.2cm, clip, width=5.1cm]{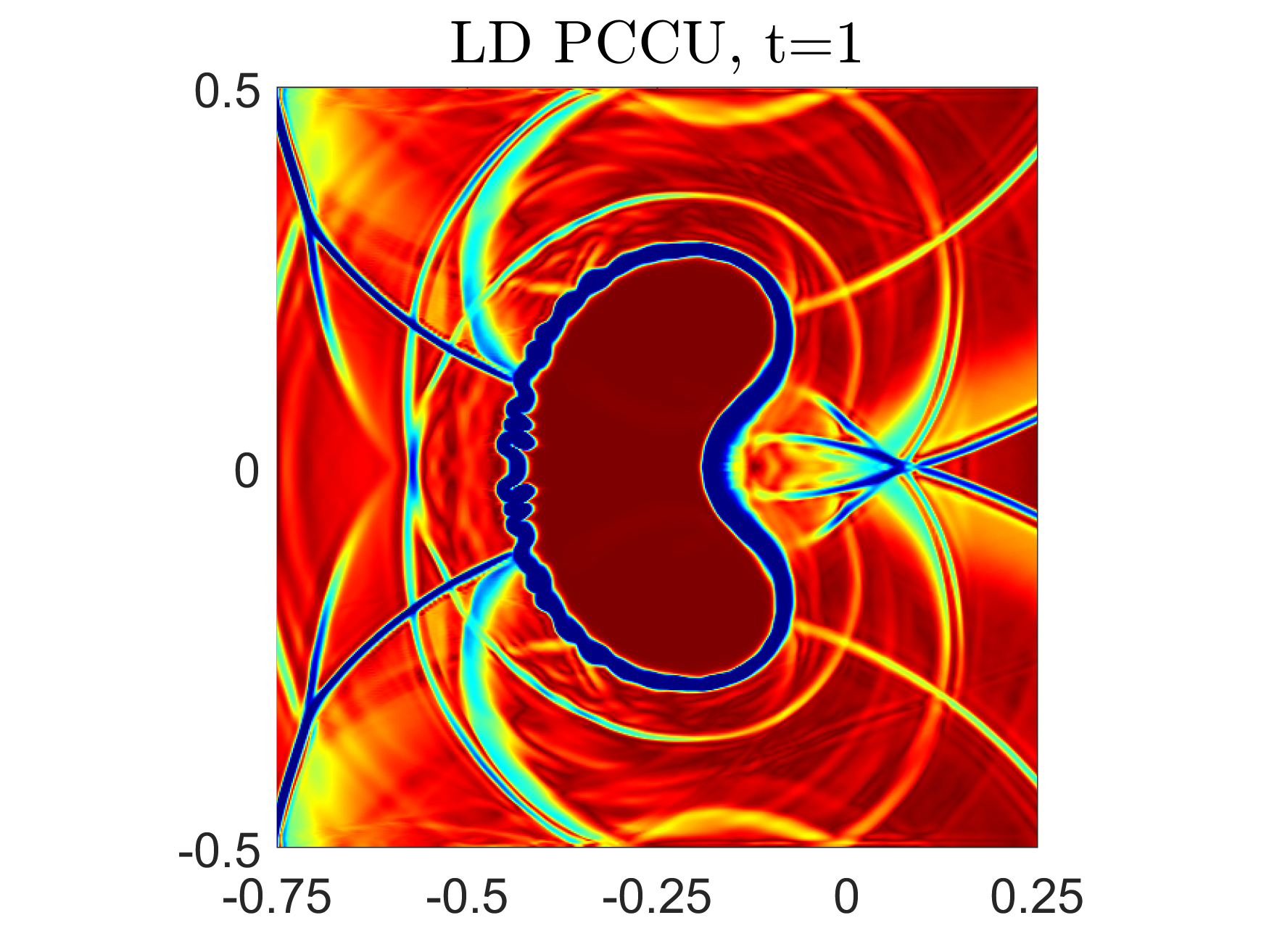}\hspace*{0.5cm}
            \includegraphics[trim=2.1cm 0.4cm 2.1cm 0.2cm, clip, width=5.1cm]{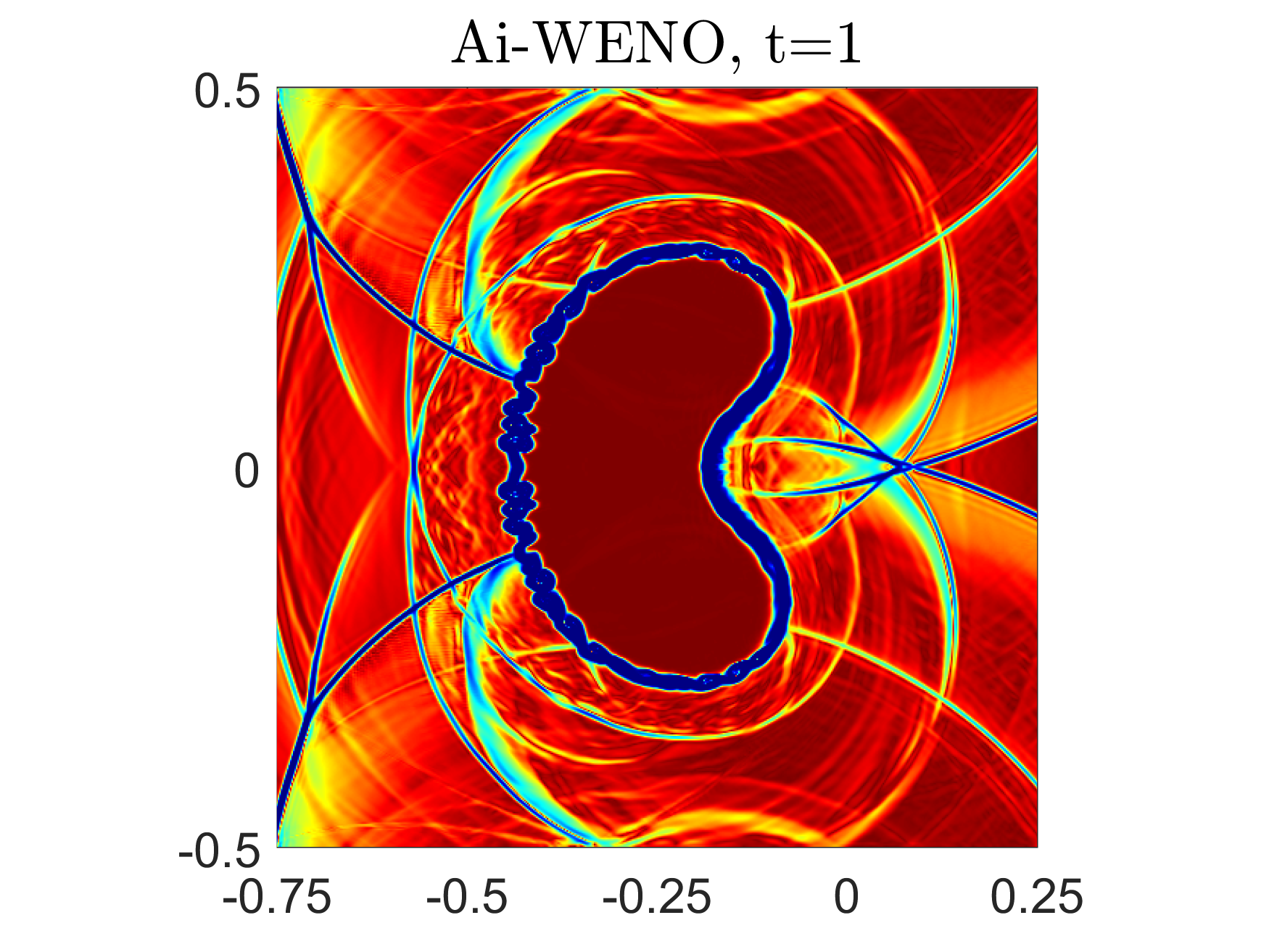}}
\vskip10pt
\centerline{\includegraphics[trim=2.1cm 0.4cm 2.1cm 0.2cm, clip, width=5.1cm]{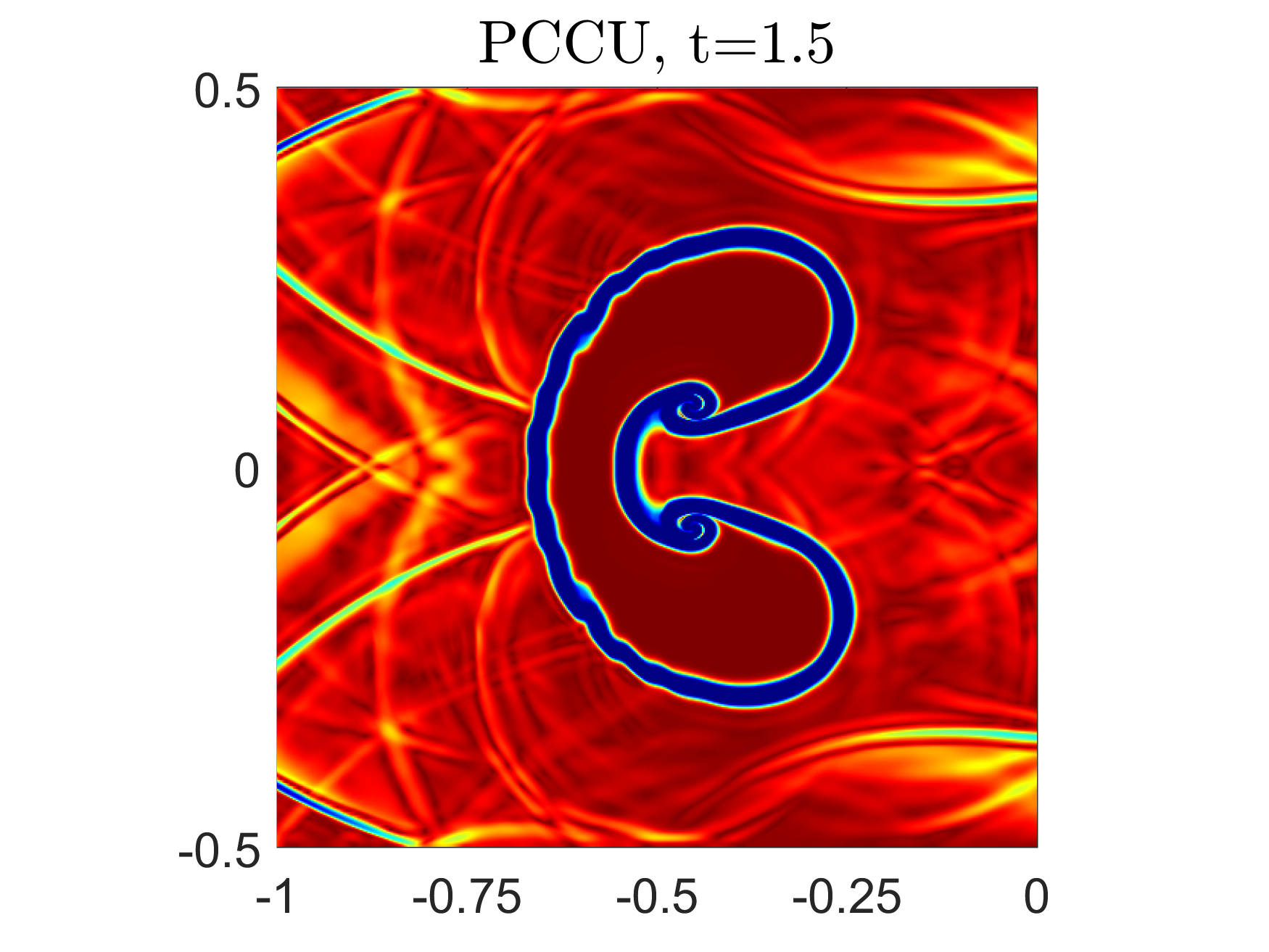}\hspace*{0.5cm}
            \includegraphics[trim=2.1cm 0.4cm 2.1cm 0.2cm, clip, width=5.1cm]{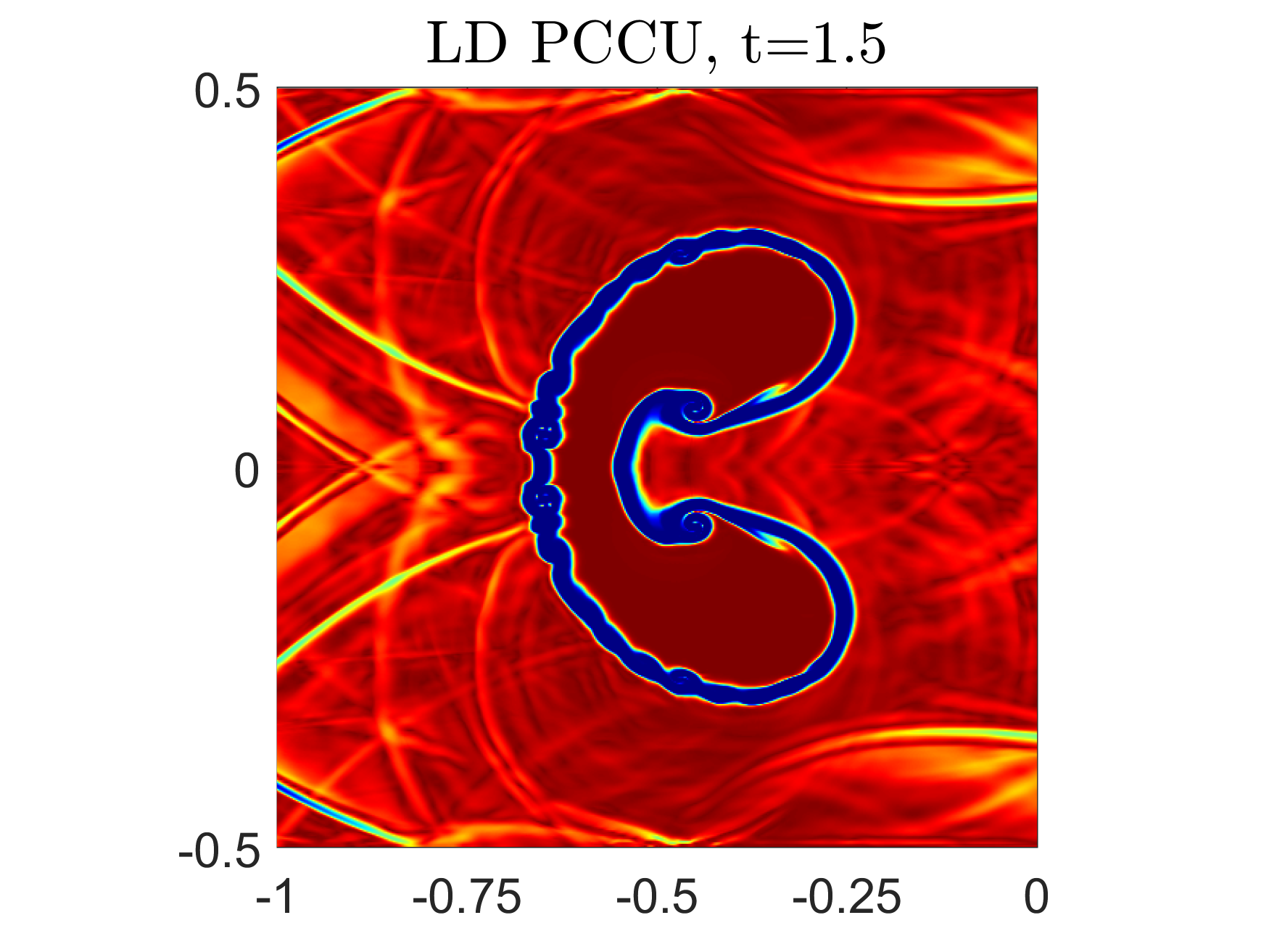}\hspace*{0.5cm}
            \includegraphics[trim=2.1cm 0.4cm 2.1cm 0.2cm, clip, width=5.1cm]{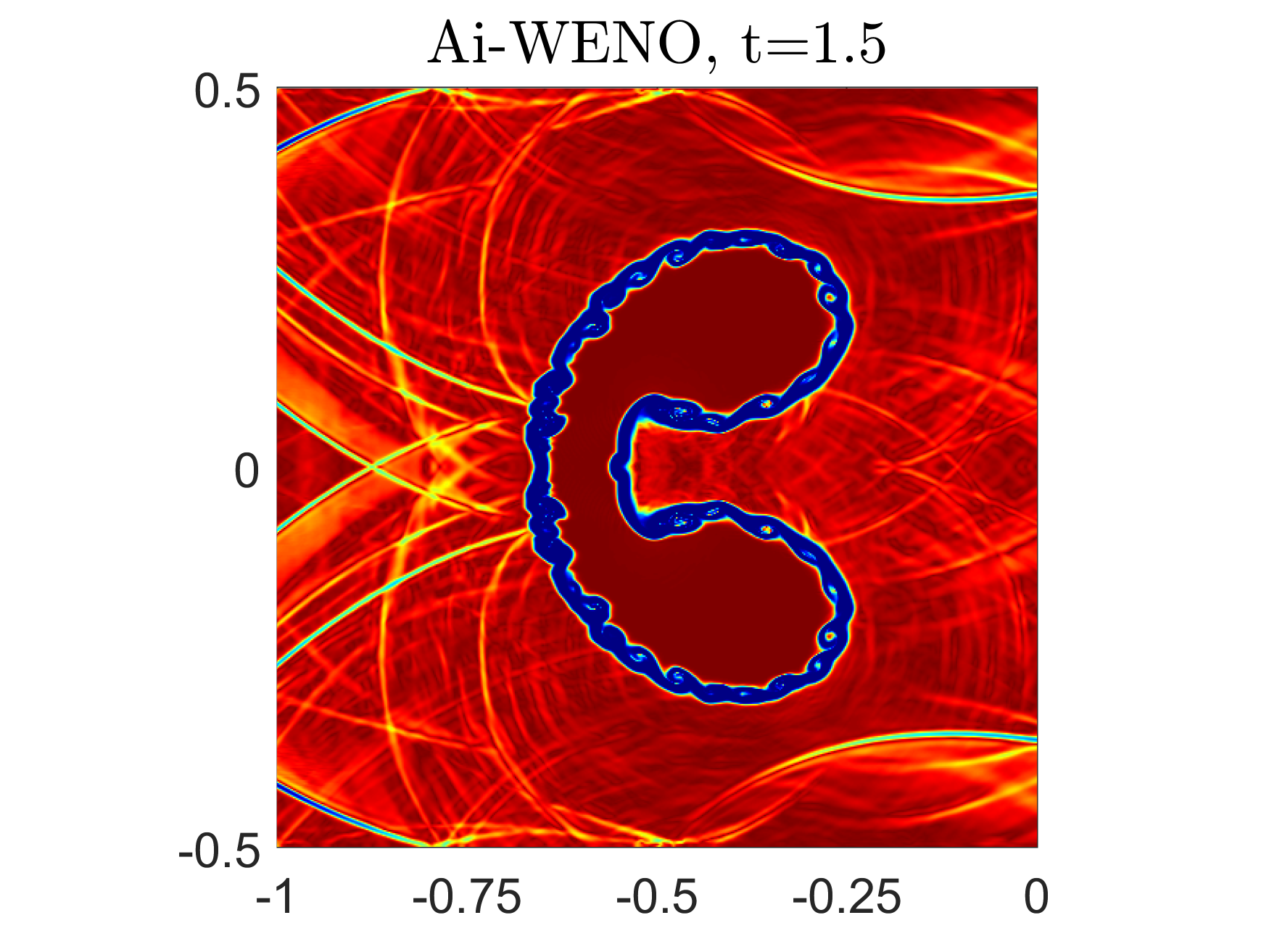}}
\caption{\sf Example 4: Shock-helium bubble interaction by the PCCU (left column), LD PCCU (middle column), and Ai-WENO (right column)
schemes at times $t=0.5$, 1, and 1.5.\label{fig43}}
\end{figure}

\begin{figure}[ht!]
\centerline{\includegraphics[trim=2.1cm 0.4cm 2.1cm 0.2cm, clip, width=5.1cm]{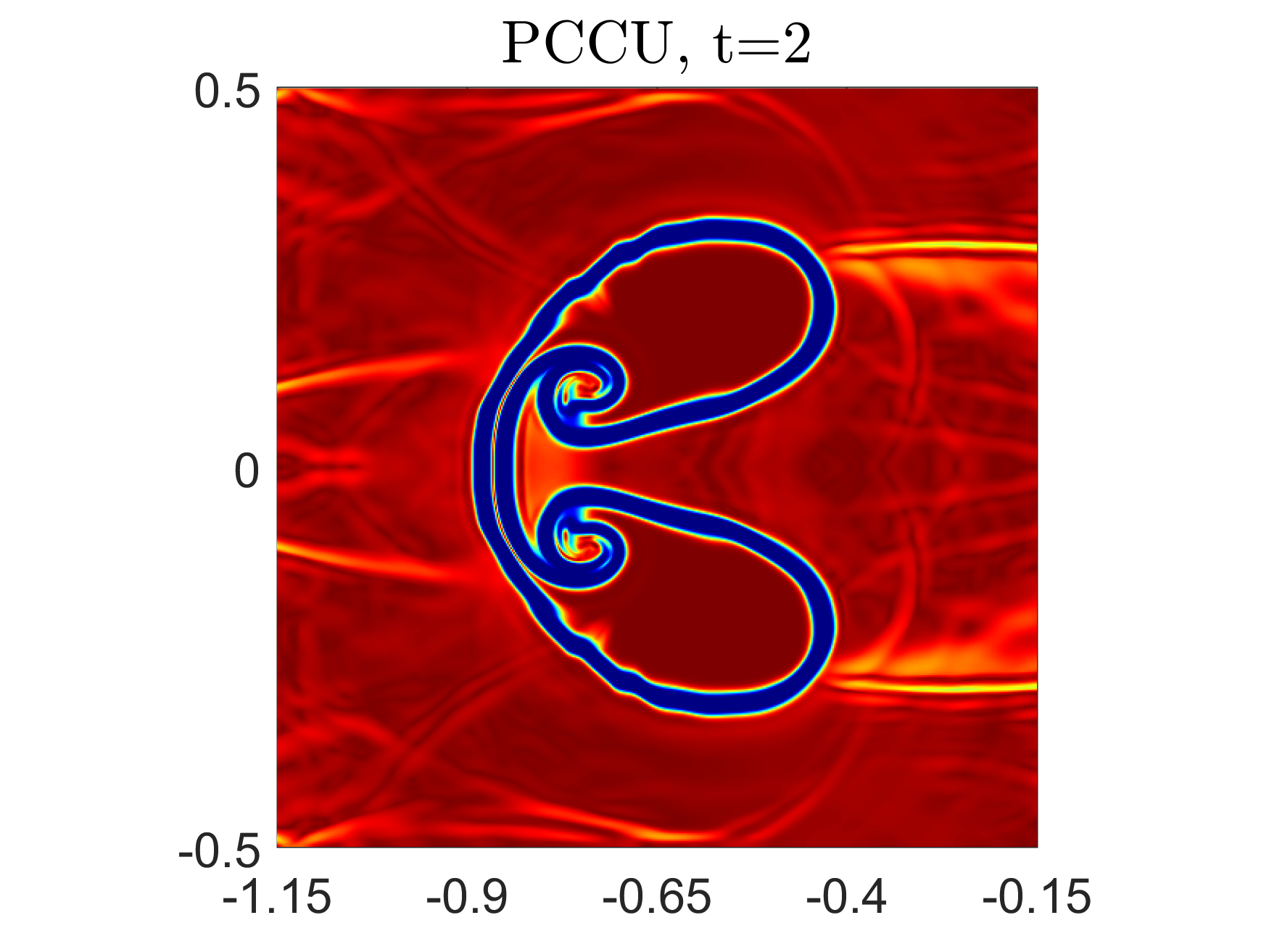}\hspace*{0.5cm}
            \includegraphics[trim=2.1cm 0.4cm 2.1cm 0.2cm, clip, width=5.1cm]{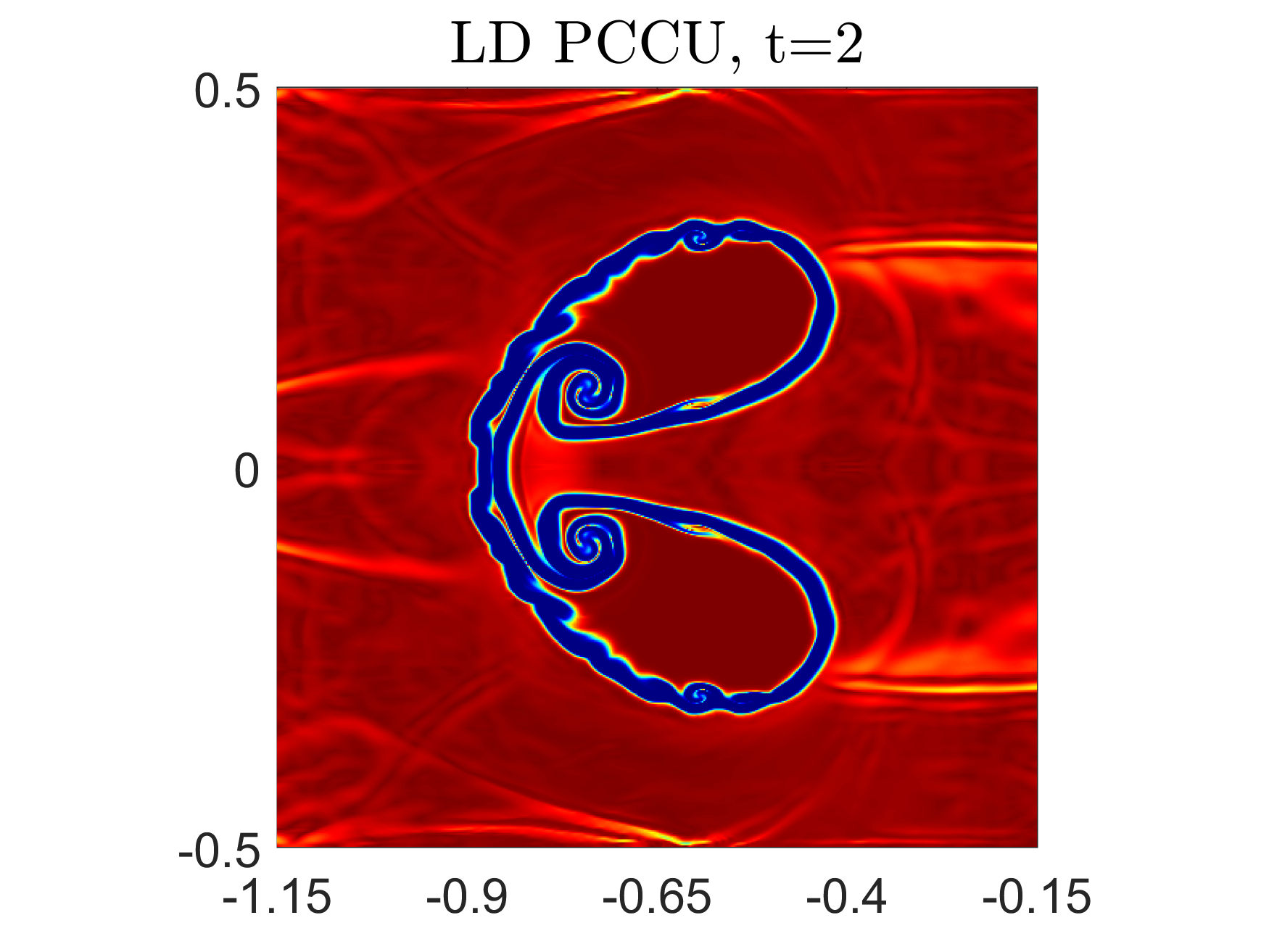}\hspace*{0.5cm}
            \includegraphics[trim=2.1cm 0.4cm 2.1cm 0.2cm, clip, width=5.1cm]{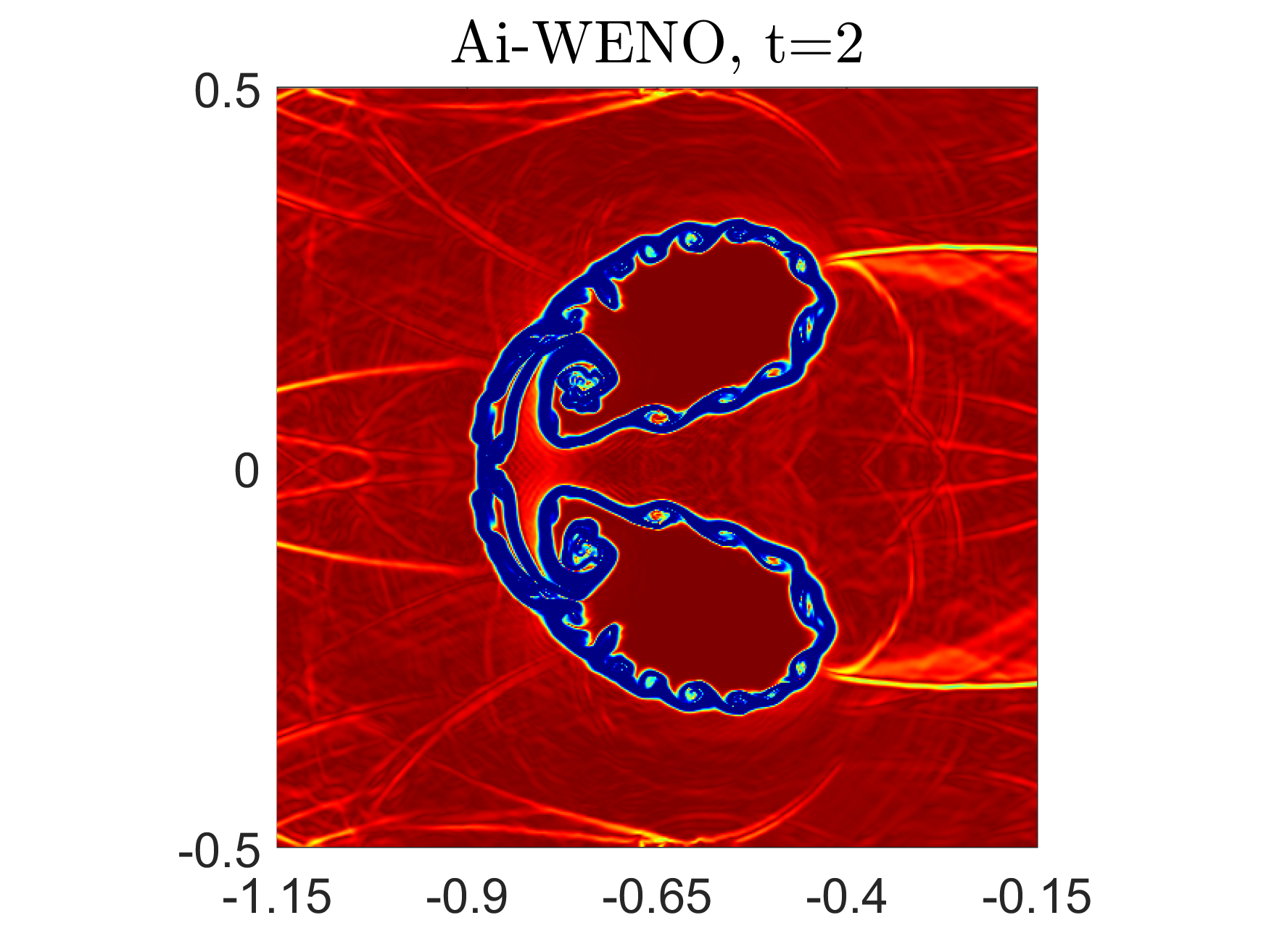}}
\vskip10pt
\centerline{\includegraphics[trim=2.1cm 0.4cm 2.1cm 0.2cm, clip, width=5.1cm]{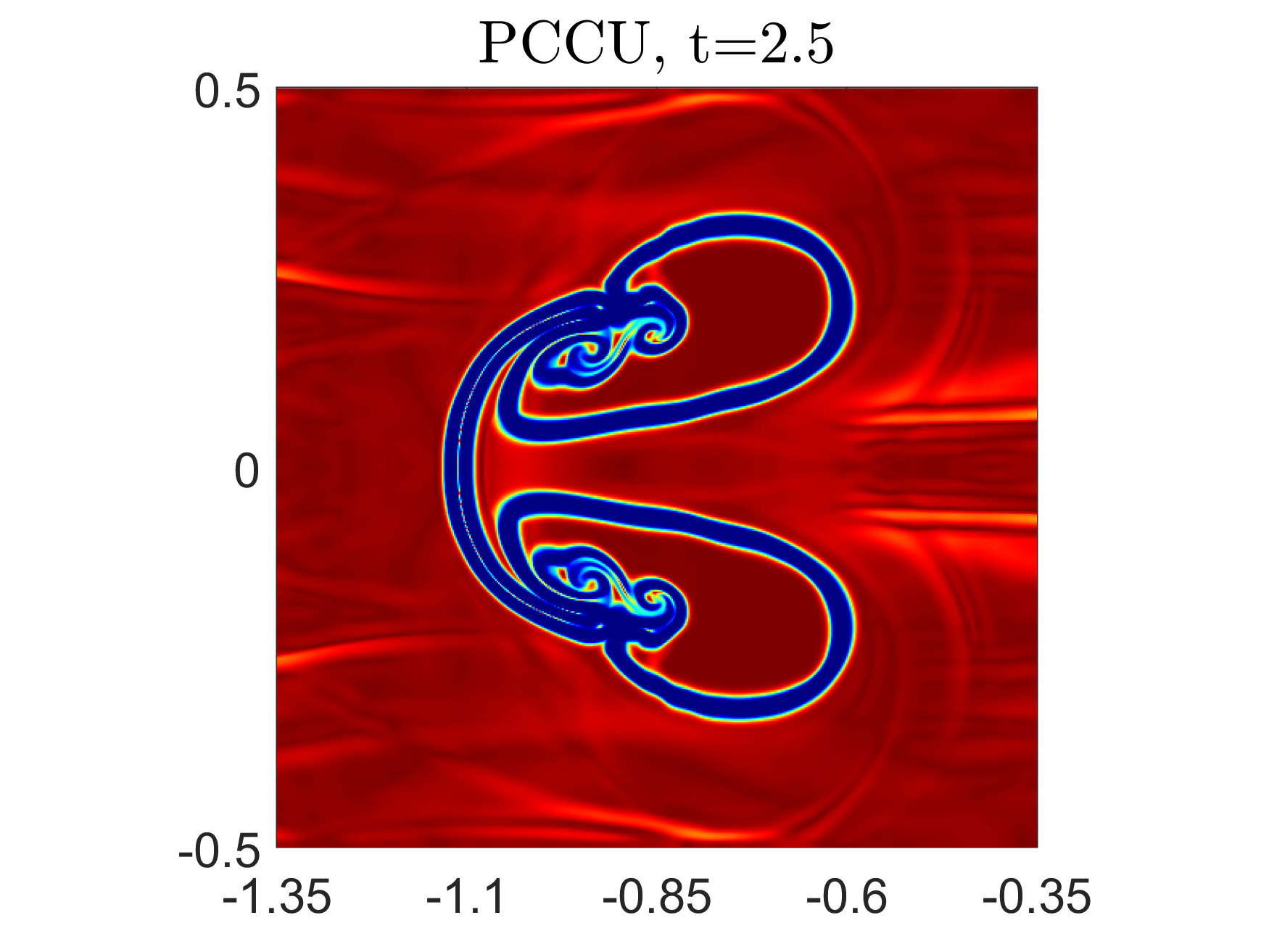}\hspace*{0.5cm}
            \includegraphics[trim=2.1cm 0.4cm 2.1cm 0.2cm, clip, width=5.1cm]{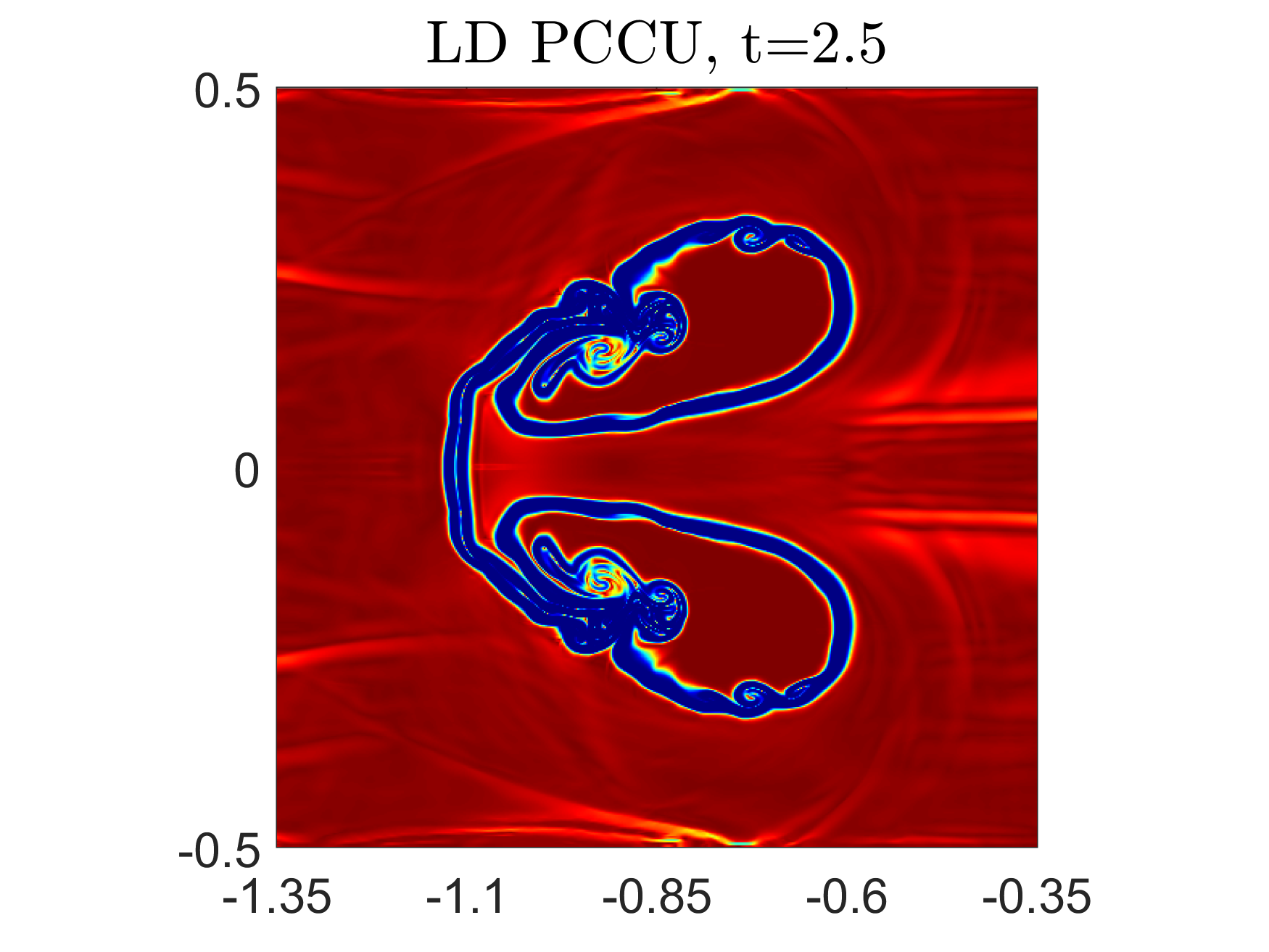}\hspace*{0.5cm}
            \includegraphics[trim=2.1cm 0.4cm 2.1cm 0.2cm, clip, width=5.1cm]{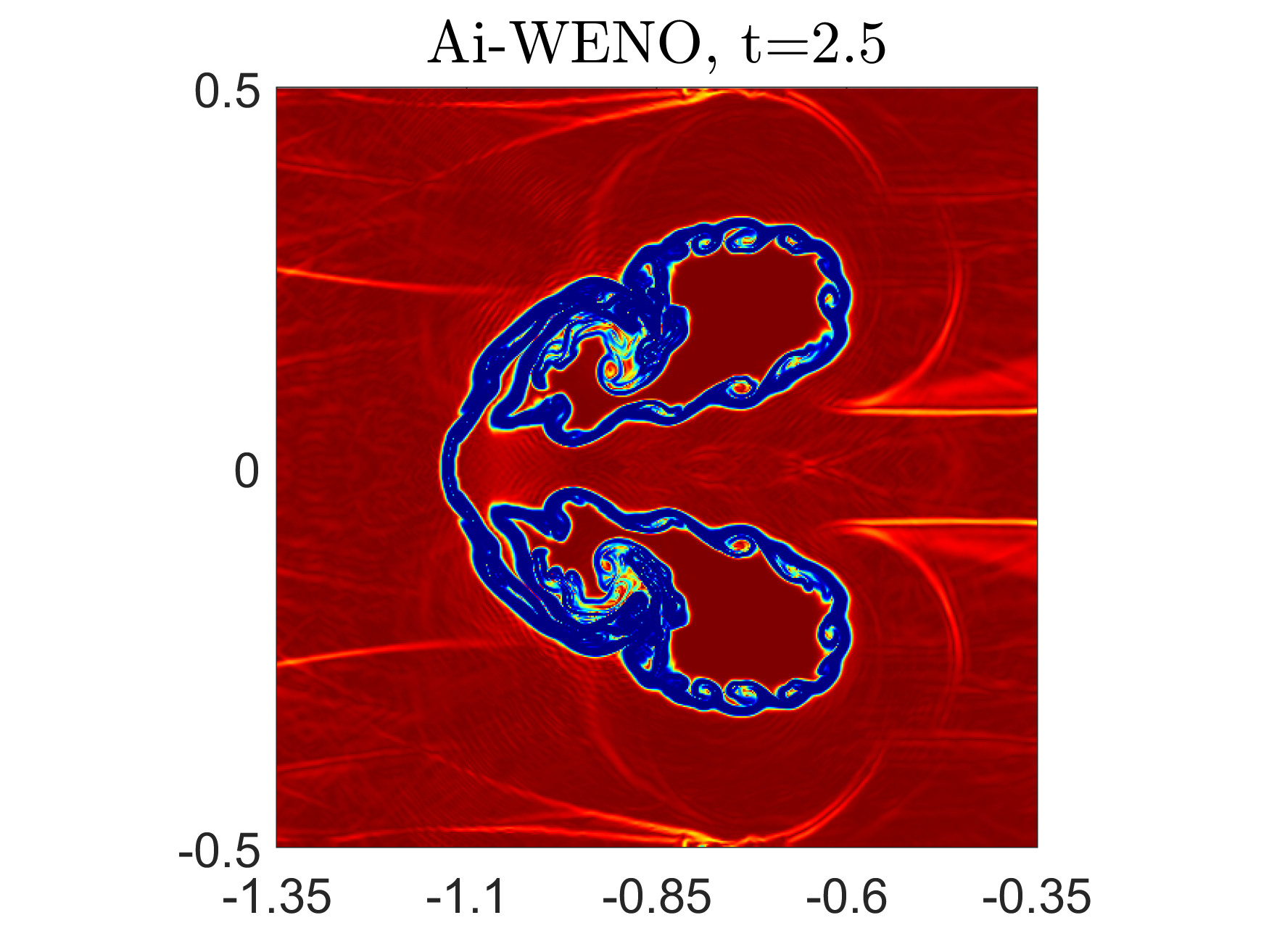}}
\vskip10pt
\centerline{\includegraphics[trim=2.1cm 0.4cm 2.1cm 0.2cm, clip, width=5.1cm]{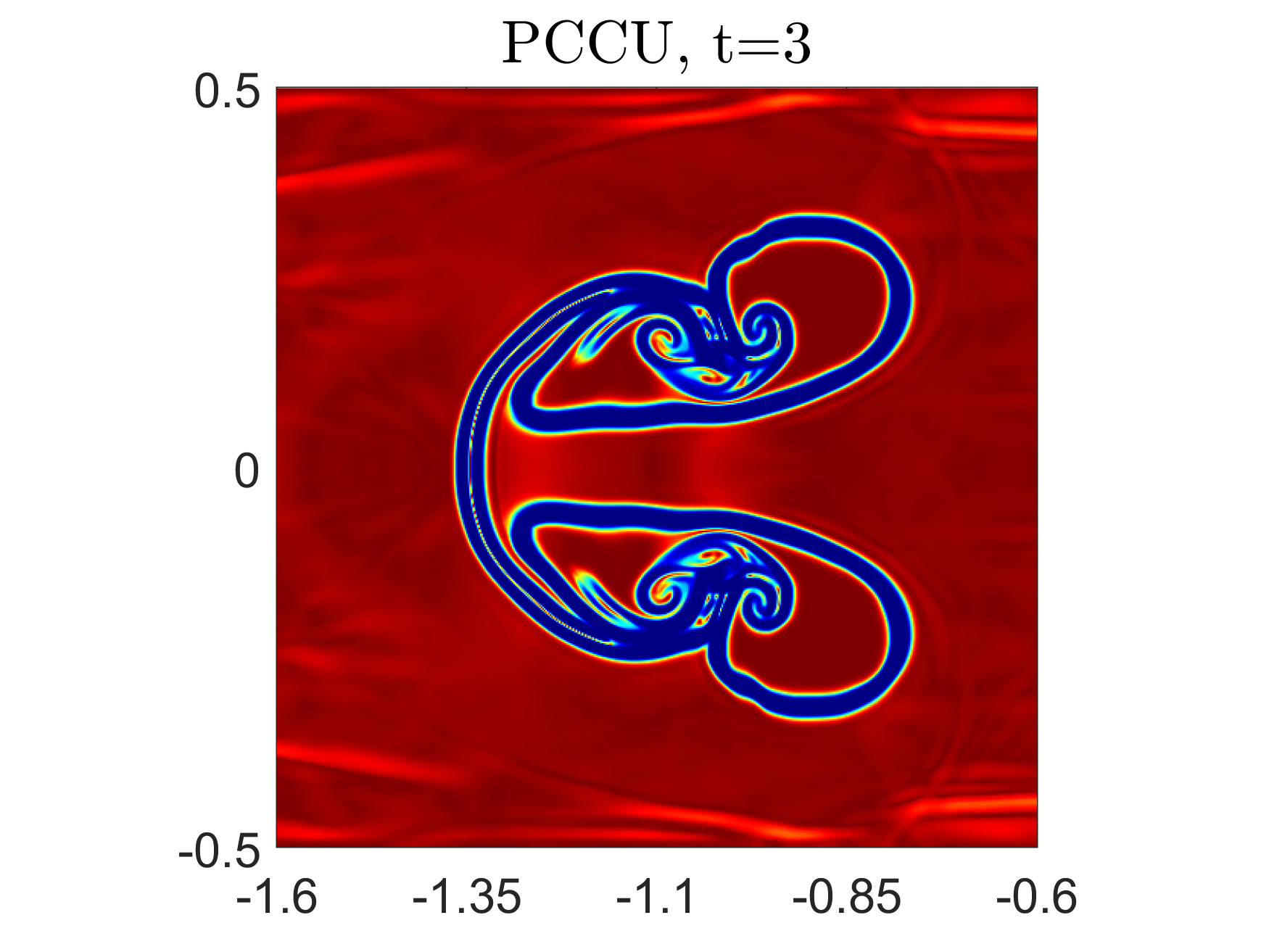}\hspace*{0.5cm}
            \includegraphics[trim=2.1cm 0.4cm 2.1cm 0.2cm, clip, width=5.1cm]{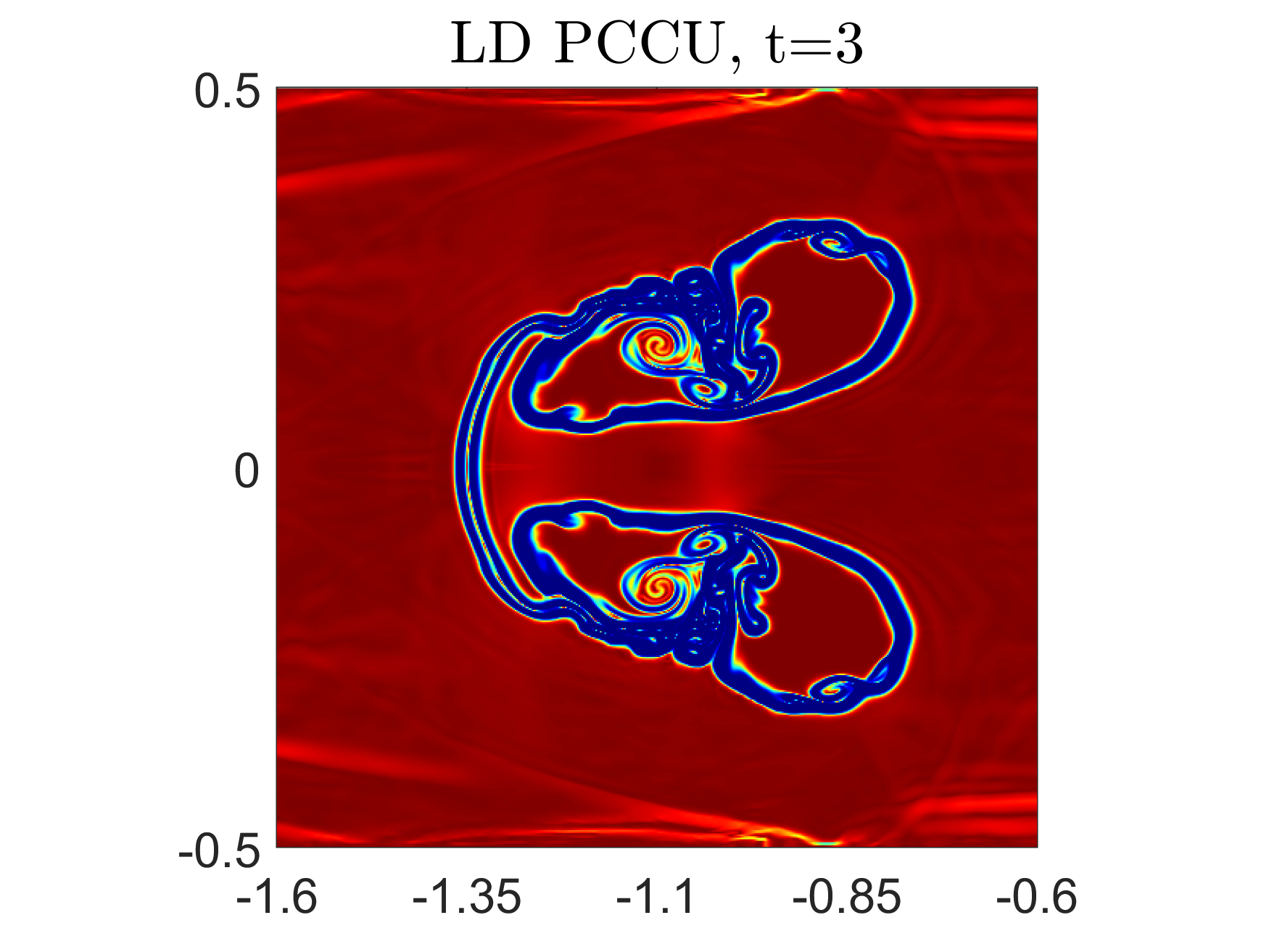}\hspace*{0.5cm}
            \includegraphics[trim=2.1cm 0.4cm 2.1cm 0.2cm, clip, width=5.1cm]{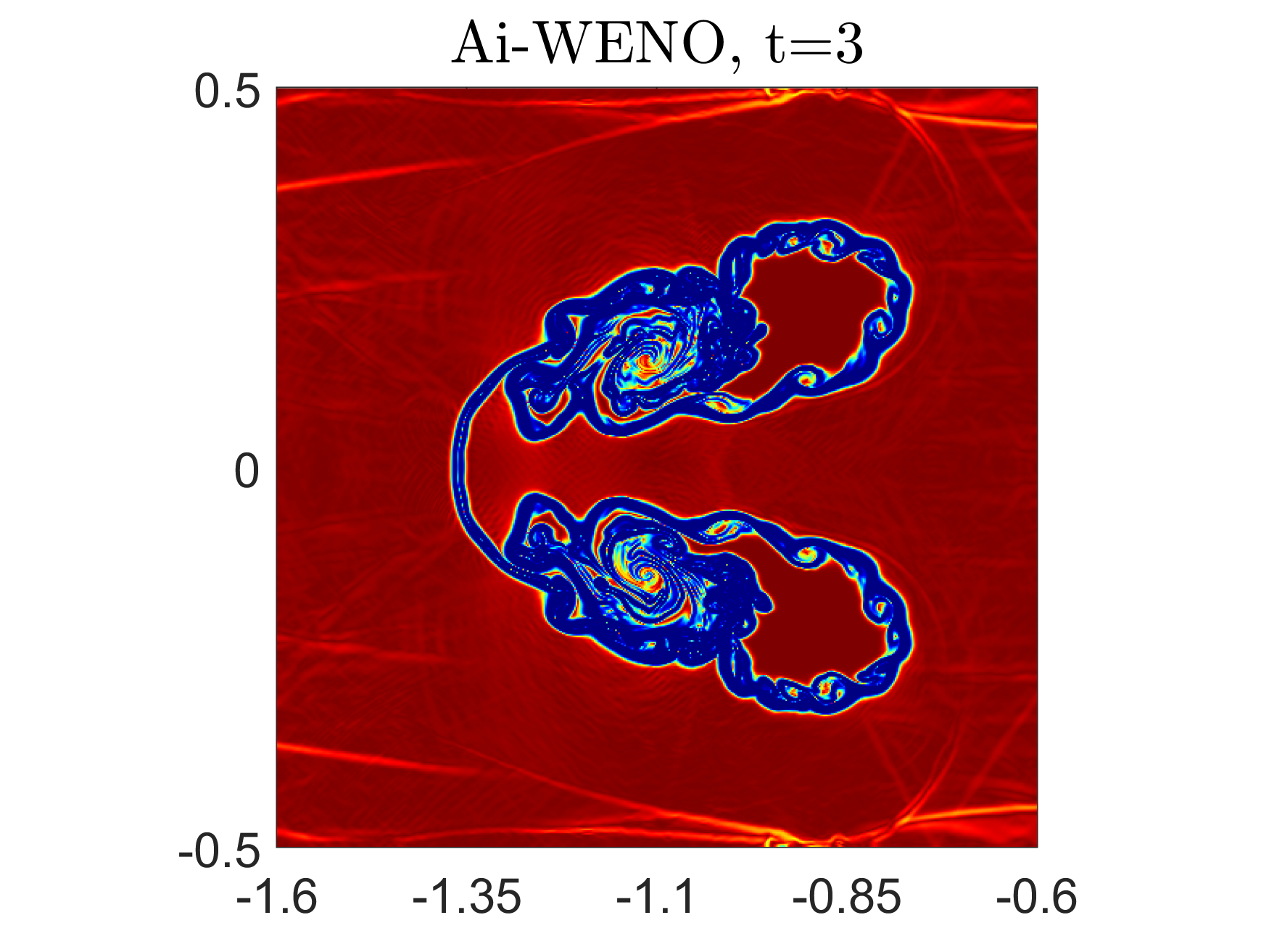}}
\caption{\sf Same as in Figure \ref{fig43}, but at larger times $t=2$, 2.5, and 3.\label{fig44}}
\end{figure}

We compute the numerical solutions until the final time $t=3$ on a uniform mesh with $\dx=\dy=1/500$. In Figures \ref{fig43} and
\ref{fig44}, we present different stages of the shock-bubble interaction computed by the PCCU, LD PCCU, and Ai-WENO schemes. Notice that the
bubble changes its shape and propagates to the left, but in order to focus on the details of the bubble structure, we only zoom at
$[\sigma,\sigma+1]\times[-0.5,0.5]$ square area containing the bubble ($\sigma$ is decreasing in time from -0.5 to -1.6). As one can observe
from the numerical results, the bubble interface develops very complex structures after the bubble is hit by the shock, especially at large
times. The obtained results are in a good agreement with the experimental findings presented in \cite{Haas1987} and the numerical results
reported in \cite{Quirk1996,chertock7,CCK_21}. At the same time, one can see that at a small time $t=0.5$, the resolution of the bubble
interface is significantly improved by the use of the LD PCCU and especially the Ai-WENO schemes. At larger times, the interface develops
instabilities which are smeared by a more dissipative PCCU scheme. The differences in the achieved resolution of the small solution
structures become even more pronounced at the larger times $t=2$, 2.5, and especially at $t=3$. This clearly indicates that the LD PCCU
scheme outperforms the PCCU one, and the further improvement in the Ai-WENO results is much more obvious than in the 1-D examples.

\subsubsection*{Example 5---Shock-R22 Bubble Interaction}
In the second 2-D example also taken from \cite{Quirk1996,chertock7}, a shock wave in the air hits the heavy resting bubble which contains
R22. The initial conditions are
\begin{equation*}
(\rho,u,v,p;\gamma,\pi_\infty)=\begin{cases}(3.1538,0,0,1;1.249,0),&\mbox{in region A},\\(1,0,0,1;1.4,0),&\mbox{in region B},\\
(4/3,-0.3535,0,1.5;1.4,0),&\mbox{in region C}.\end{cases}
\end{equation*}
The regions A, B, and C are the same as in Example 4 and they are specified in Figure \ref{fig44}. In this example, we impose the same
boundary conditions and use the same computational domain as in Example 4.

\begin{figure}[ht!]
\centerline{\includegraphics[trim=2.1cm 0.4cm 2.1cm 0.2cm, clip, width=5.1cm]{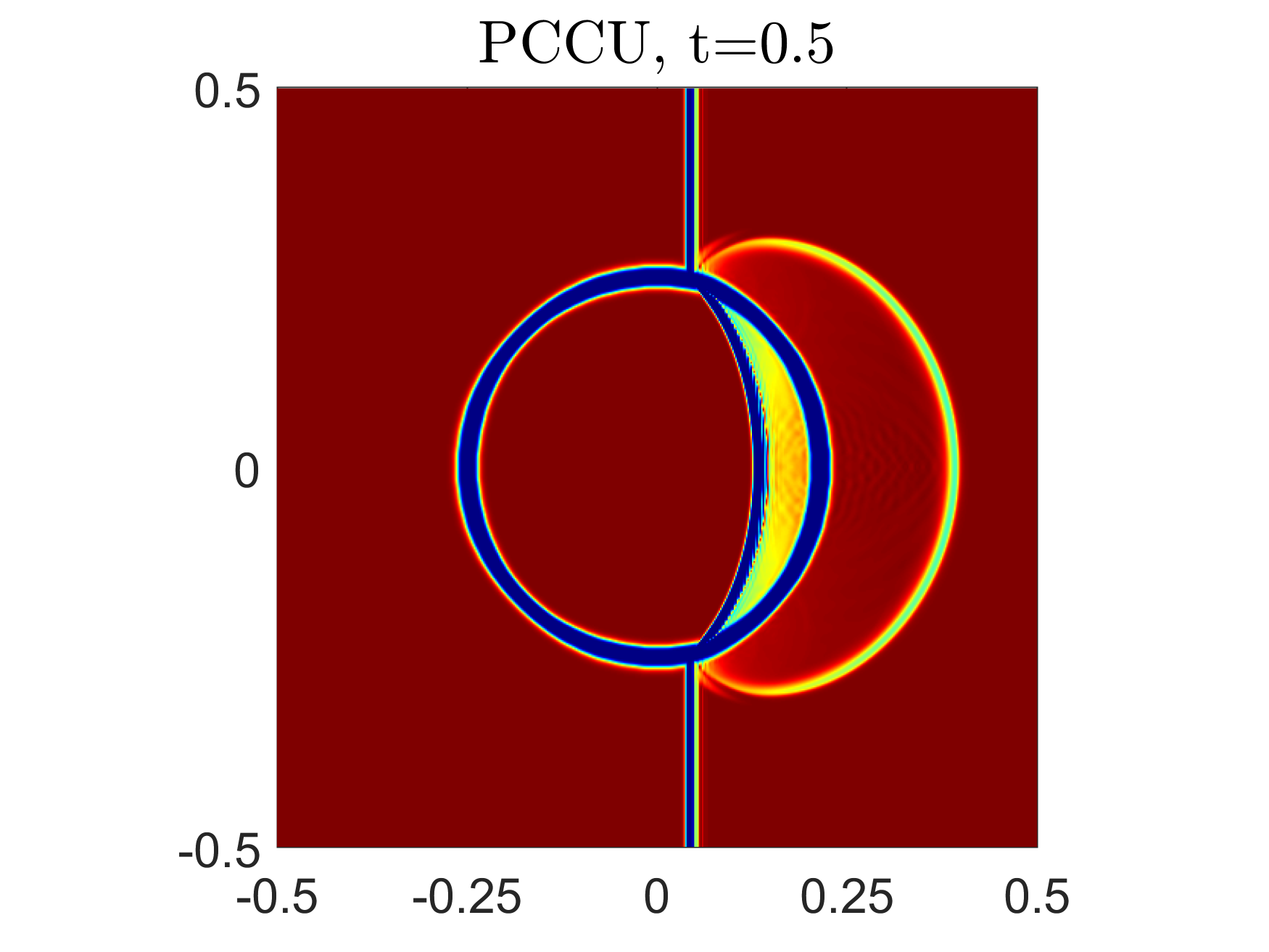}\hspace*{0.5cm}
            \includegraphics[trim=2.1cm 0.4cm 2.1cm 0.2cm, clip, width=5.1cm]{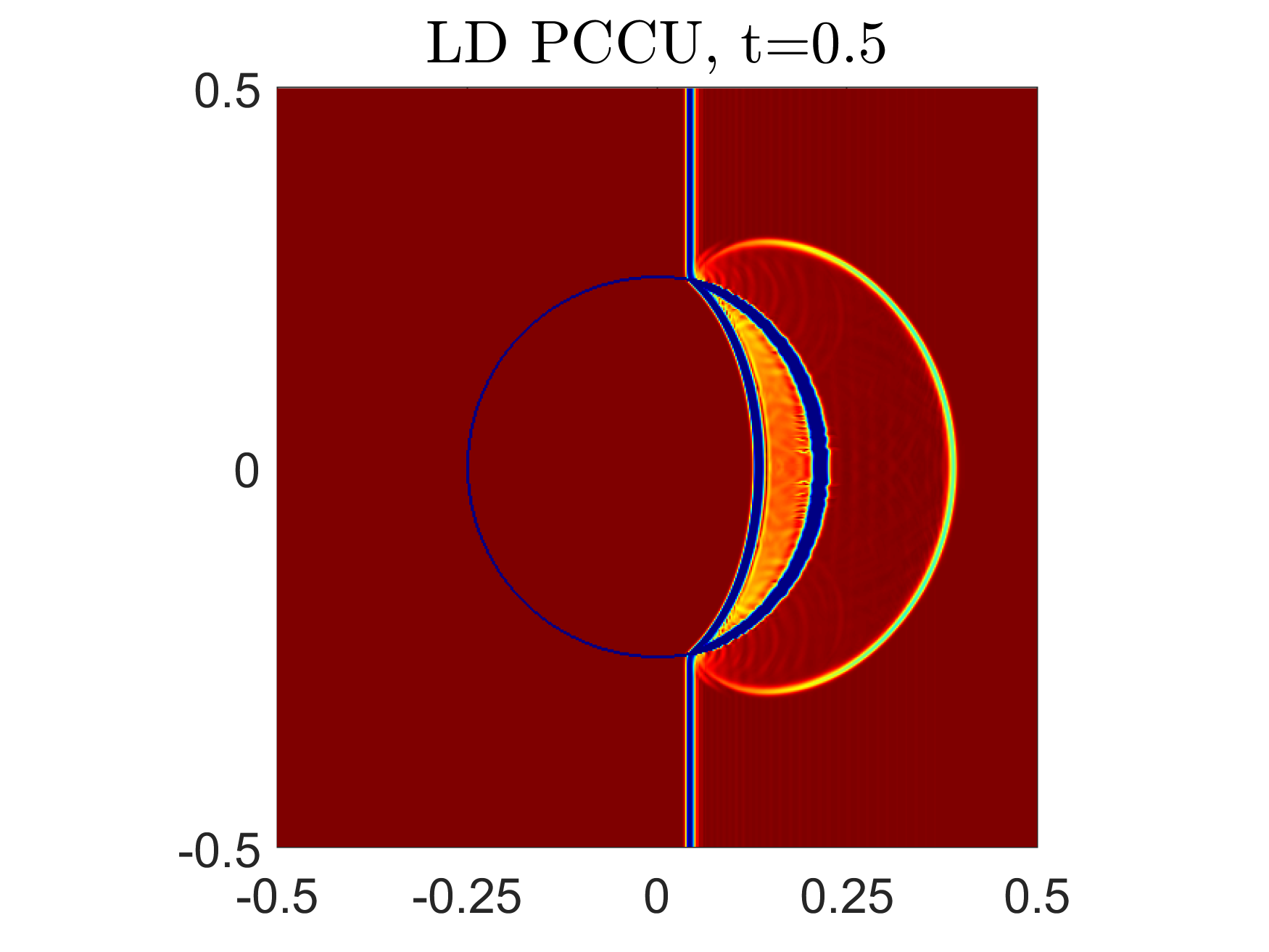}\hspace*{0.5cm}
            \includegraphics[trim=2.1cm 0.4cm 2.1cm 0.2cm, clip, width=5.1cm]{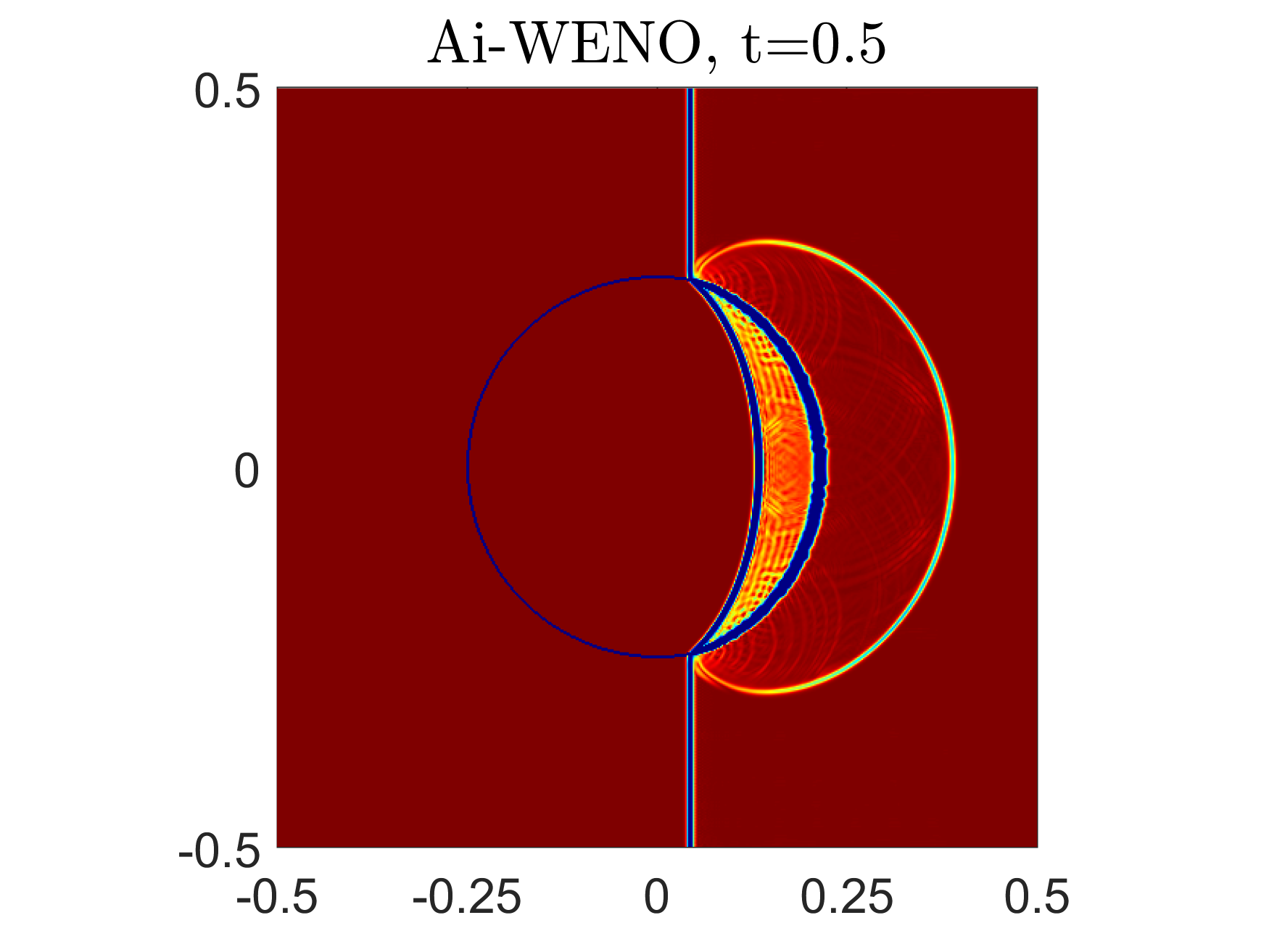}}
\vskip10pt
\centerline{\includegraphics[trim=2.1cm 0.4cm 2.1cm 0.2cm, clip, width=5.1cm]{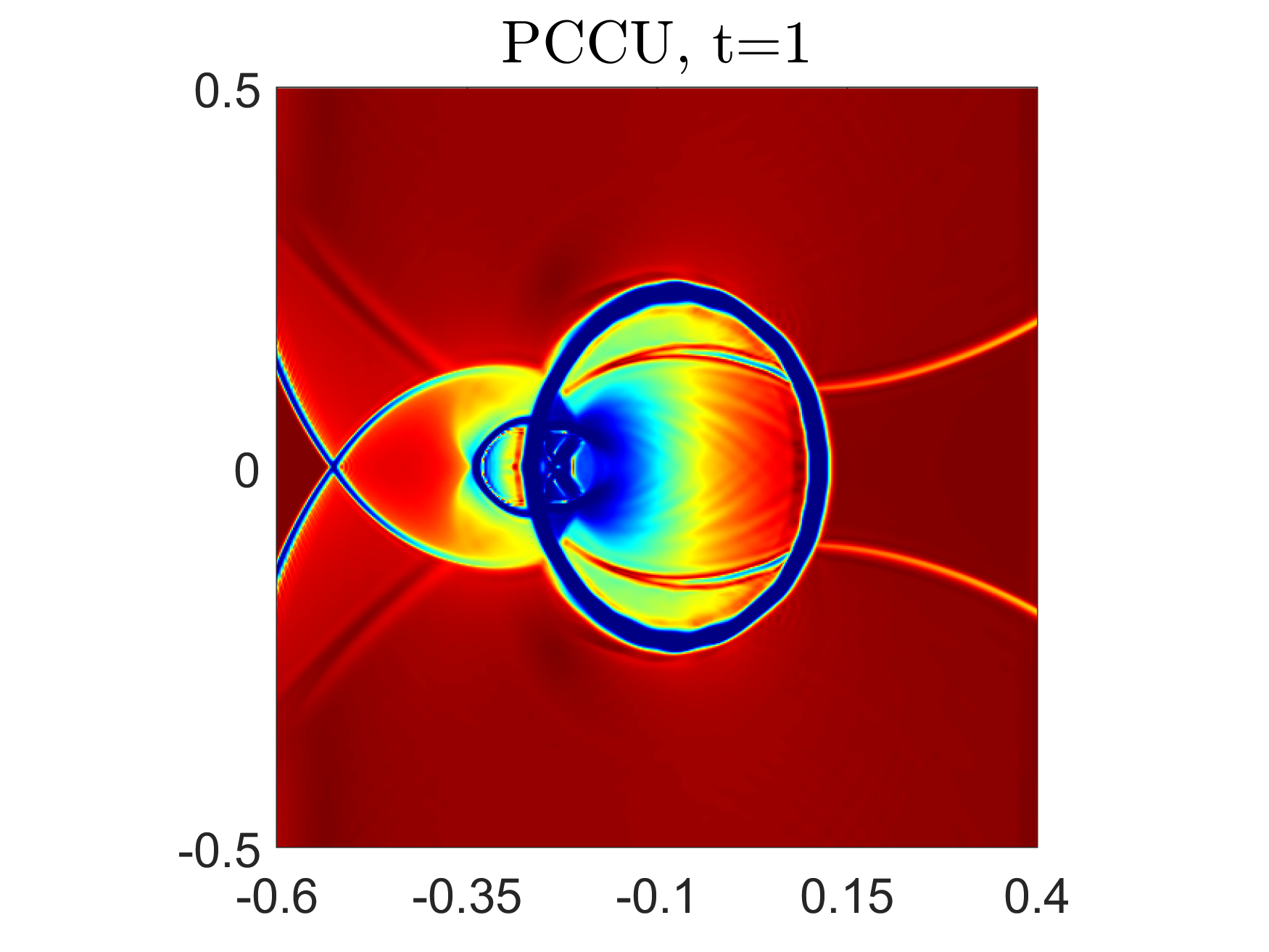}\hspace*{0.5cm}
            \includegraphics[trim=2.1cm 0.4cm 2.1cm 0.2cm, clip, width=5.1cm]{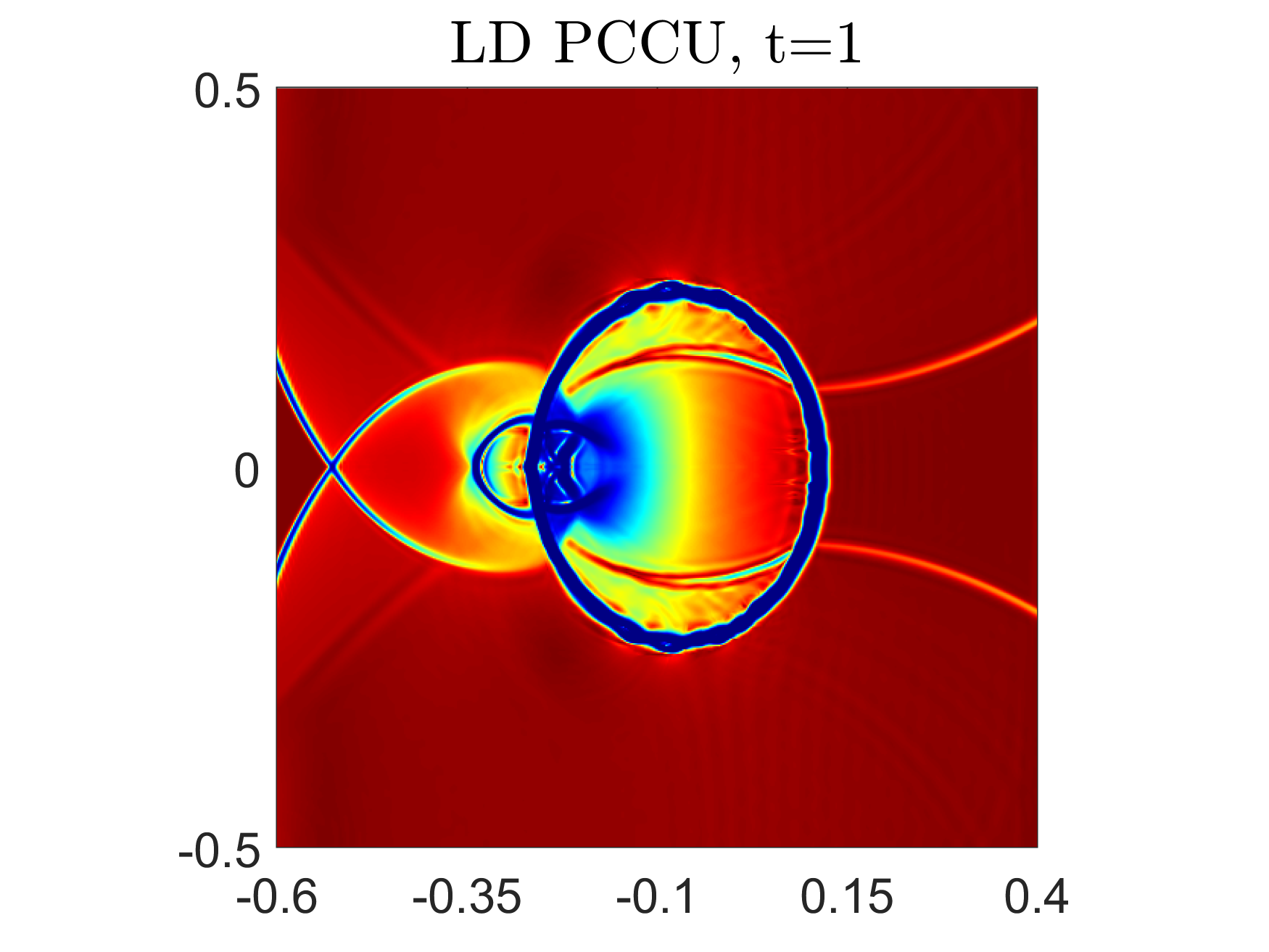}\hspace*{0.5cm}
            \includegraphics[trim=2.1cm 0.4cm 2.1cm 0.2cm, clip, width=5.1cm]{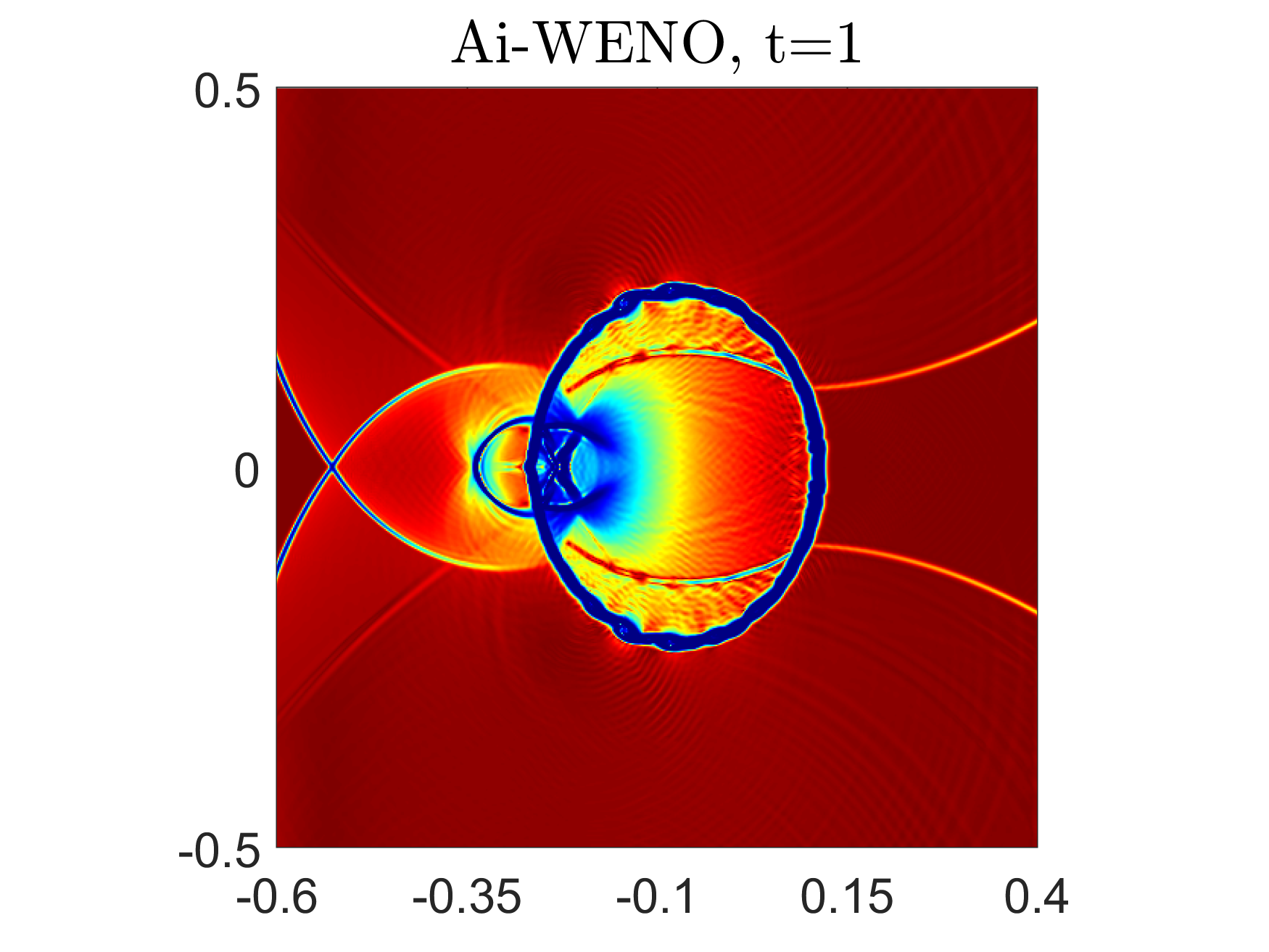}}
\caption{\sf Example 5: Shock-R22 bubble interaction by the PCCU (left column), LD PCCU (middle column), and Ai-WENO (right column) schemes
at times $t=0.5$ and 1.\label{fig45}}
\end{figure}

\begin{figure}[ht!]
\centerline{\includegraphics[trim=2.1cm 0.4cm 2.1cm 0.2cm, clip, width=5.1cm]{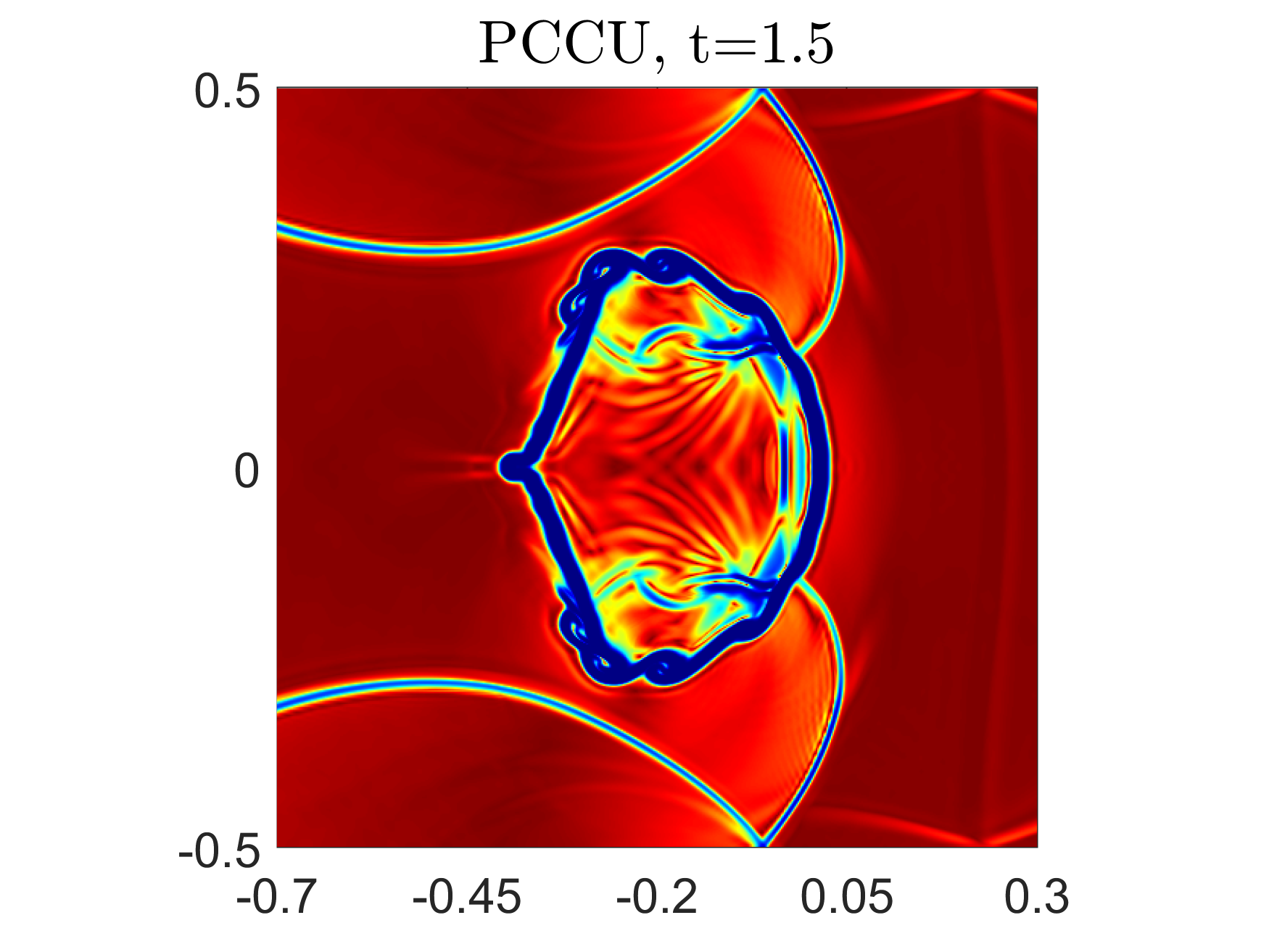}\hspace*{0.5cm}
            \includegraphics[trim=2.1cm 0.4cm 2.1cm 0.2cm, clip, width=5.1cm]{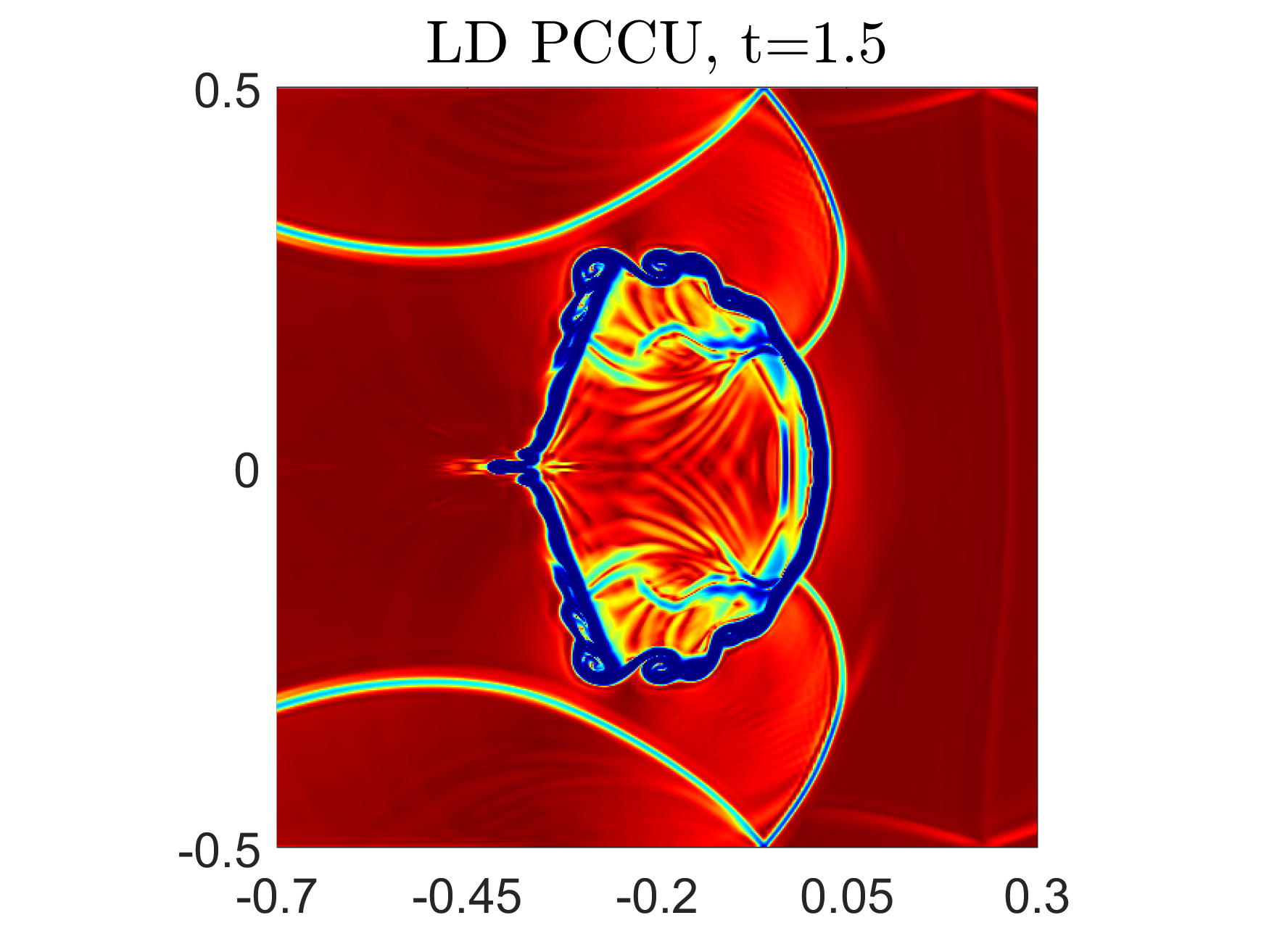}\hspace*{0.5cm}
            \includegraphics[trim=2.1cm 0.4cm 2.1cm 0.2cm, clip, width=5.1cm]{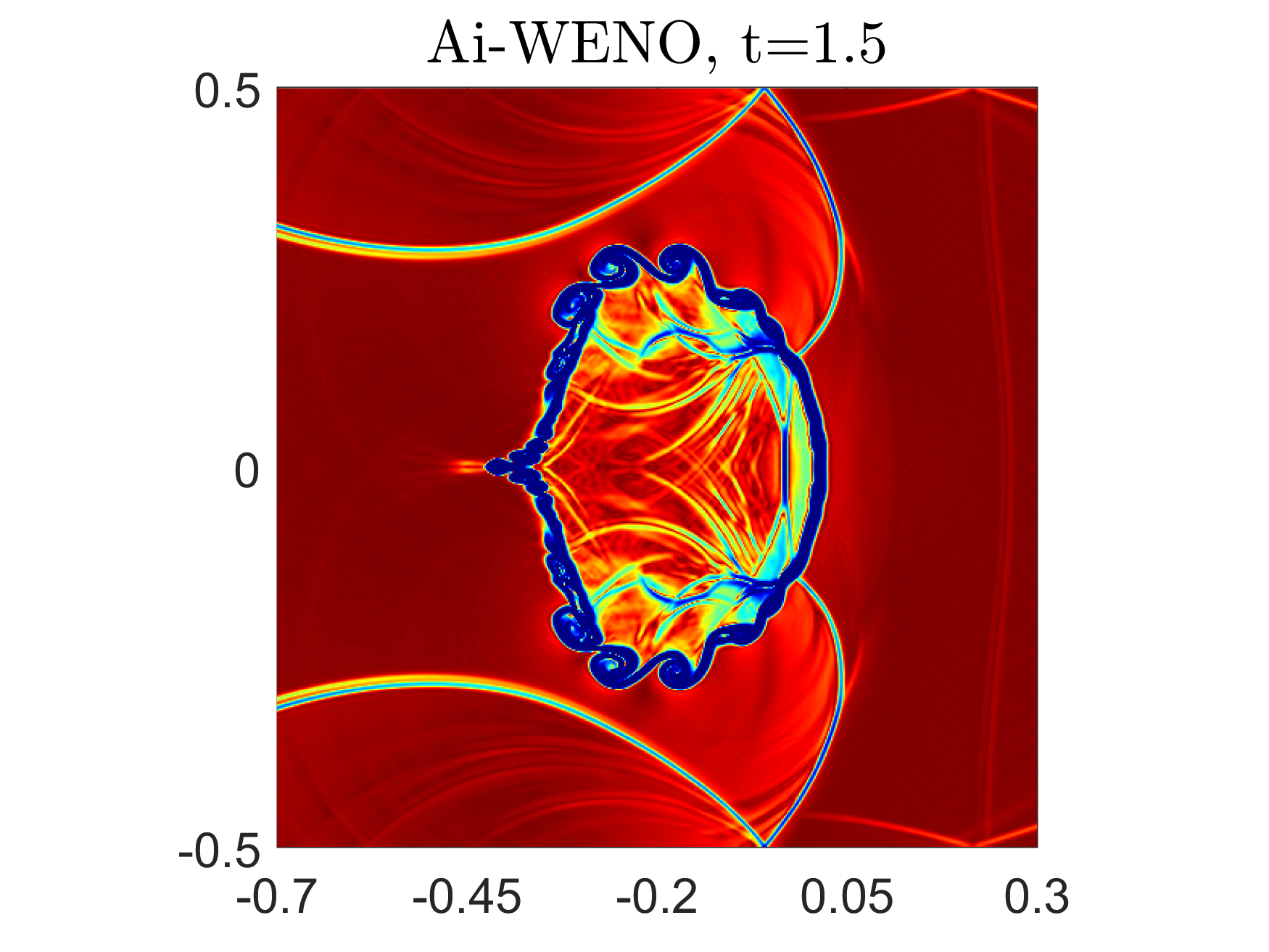}}
\vskip10pt
\centerline{\includegraphics[trim=2.1cm 0.4cm 2.1cm 0.2cm, clip, width=5.1cm]{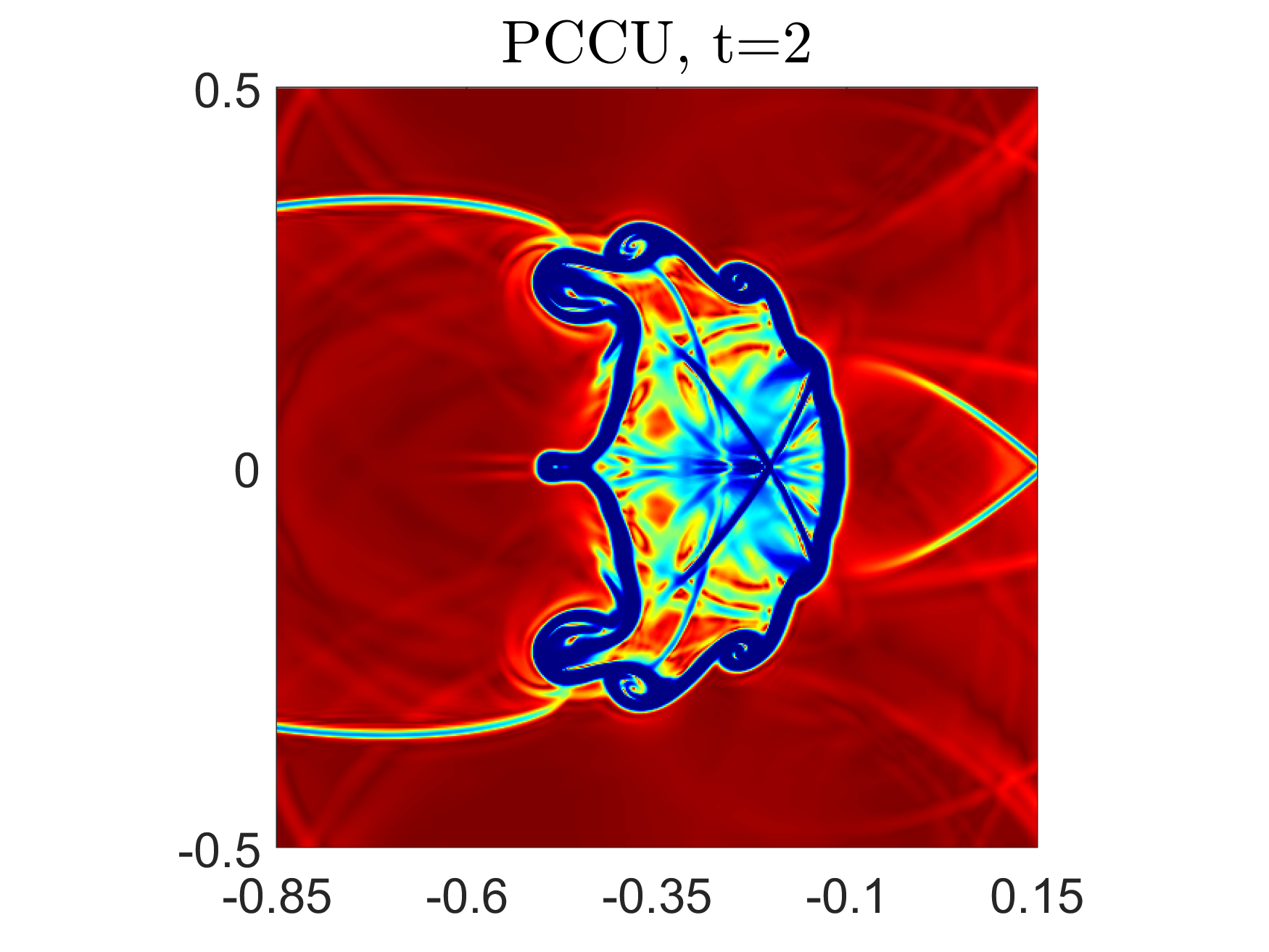}\hspace*{0.5cm}
            \includegraphics[trim=2.1cm 0.4cm 2.1cm 0.2cm, clip, width=5.1cm]{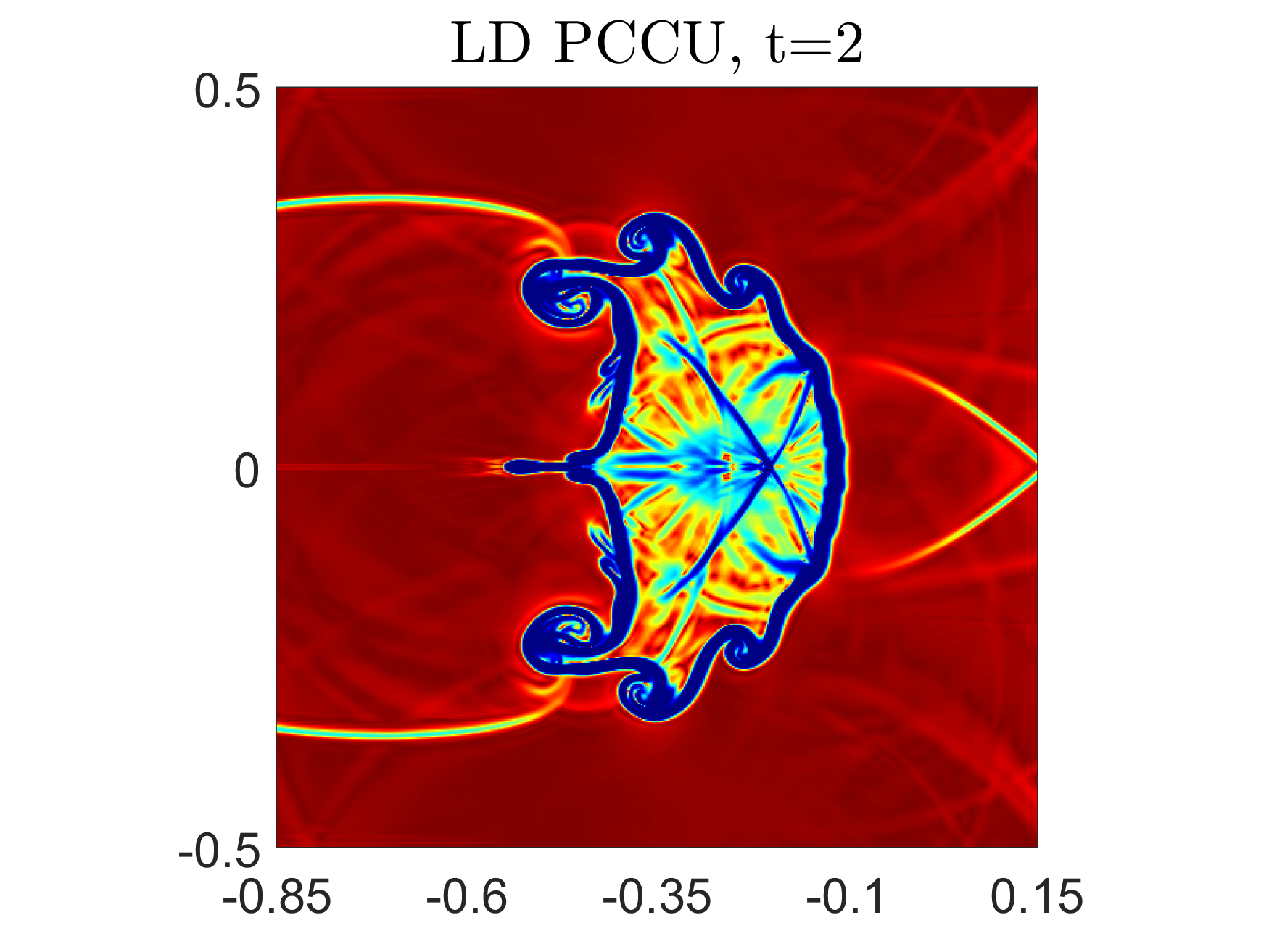}\hspace*{0.5cm}
            \includegraphics[trim=2.1cm 0.4cm 2.1cm 0.2cm, clip, width=5.1cm]{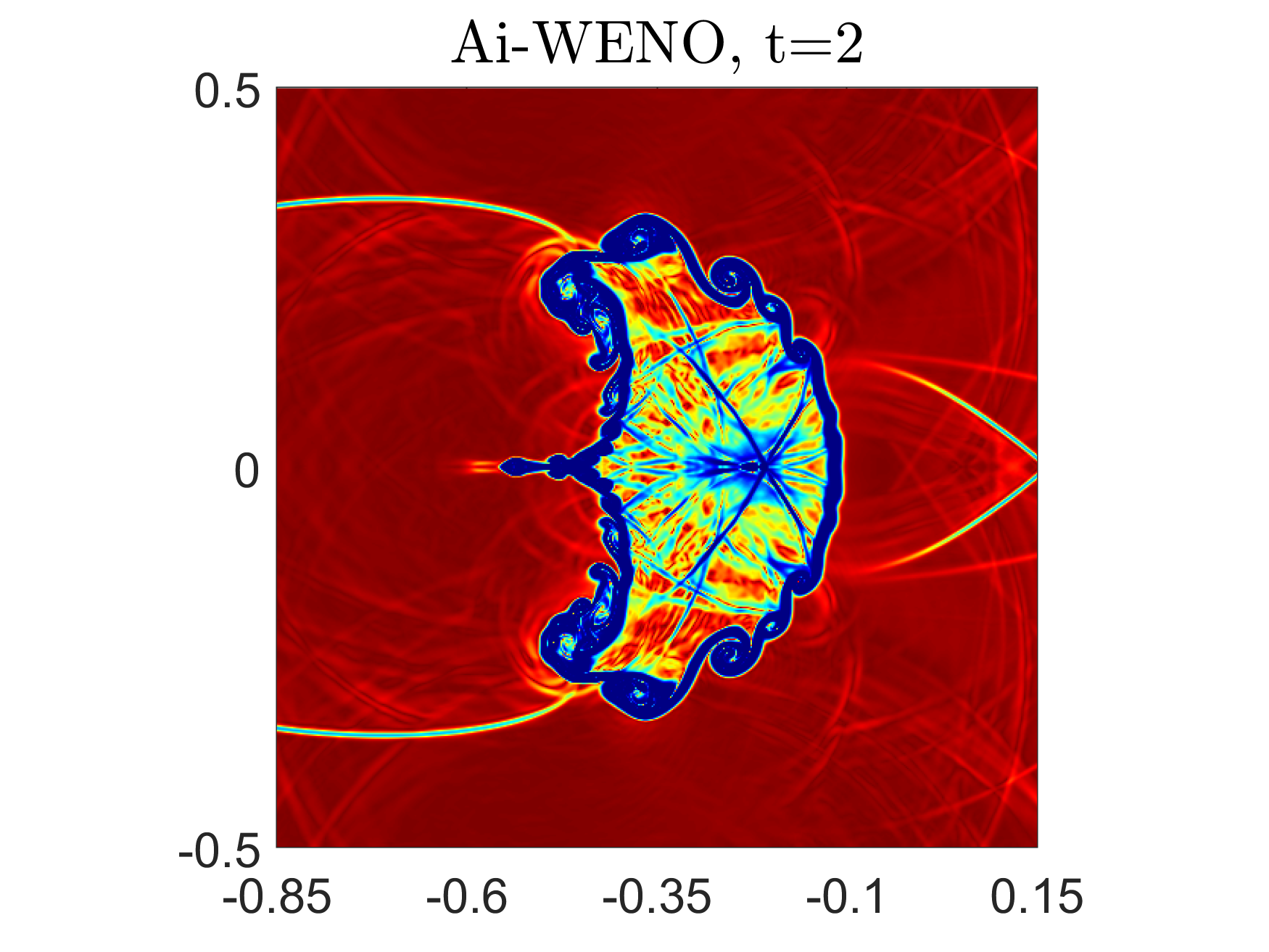}}
\vskip10pt
\centerline{\includegraphics[trim=2.1cm 0.4cm 2.1cm 0.2cm, clip, width=5.1cm]{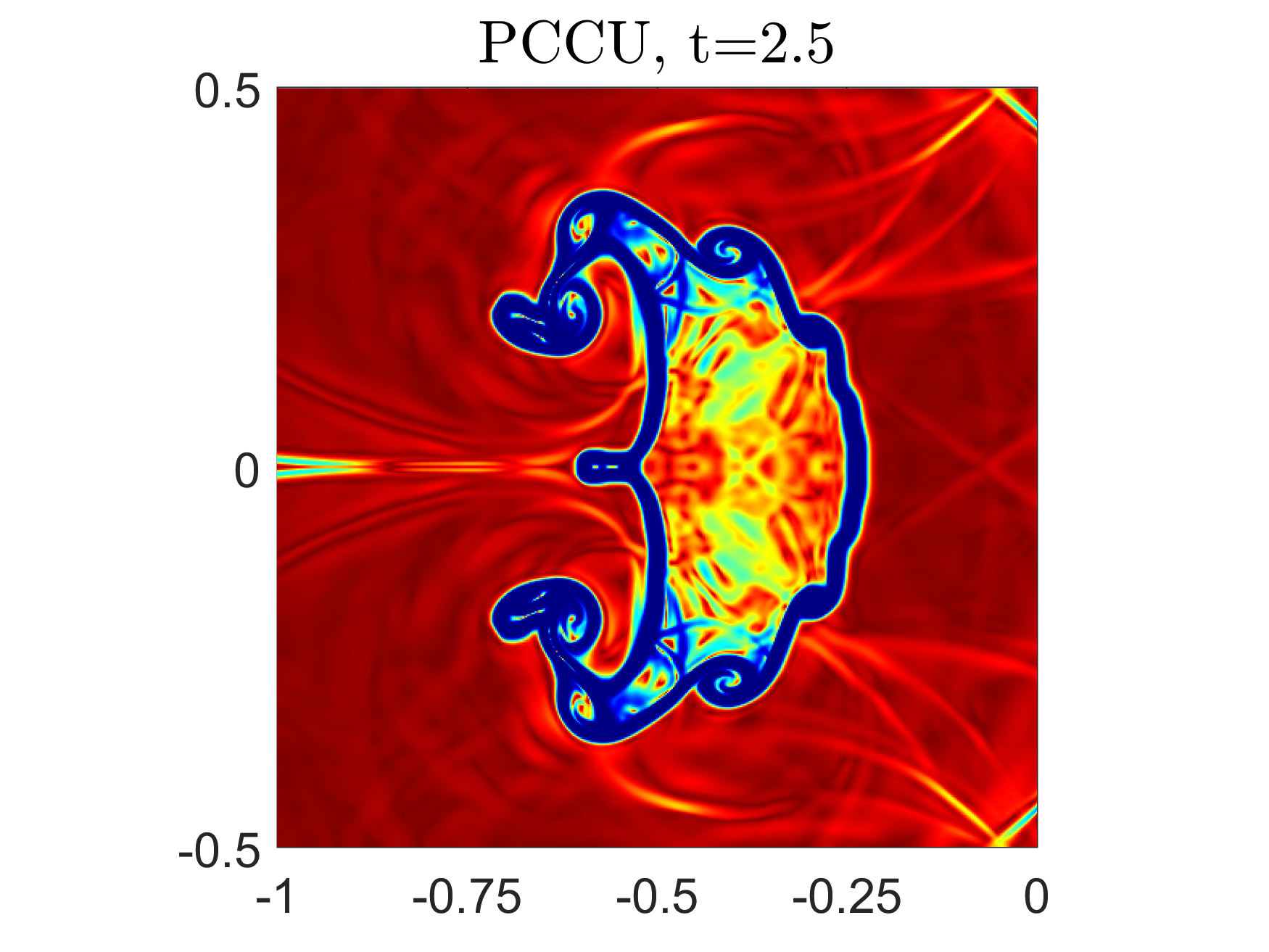}\hspace*{0.5cm}
            \includegraphics[trim=2.1cm 0.4cm 2.1cm 0.2cm, clip, width=5.1cm]{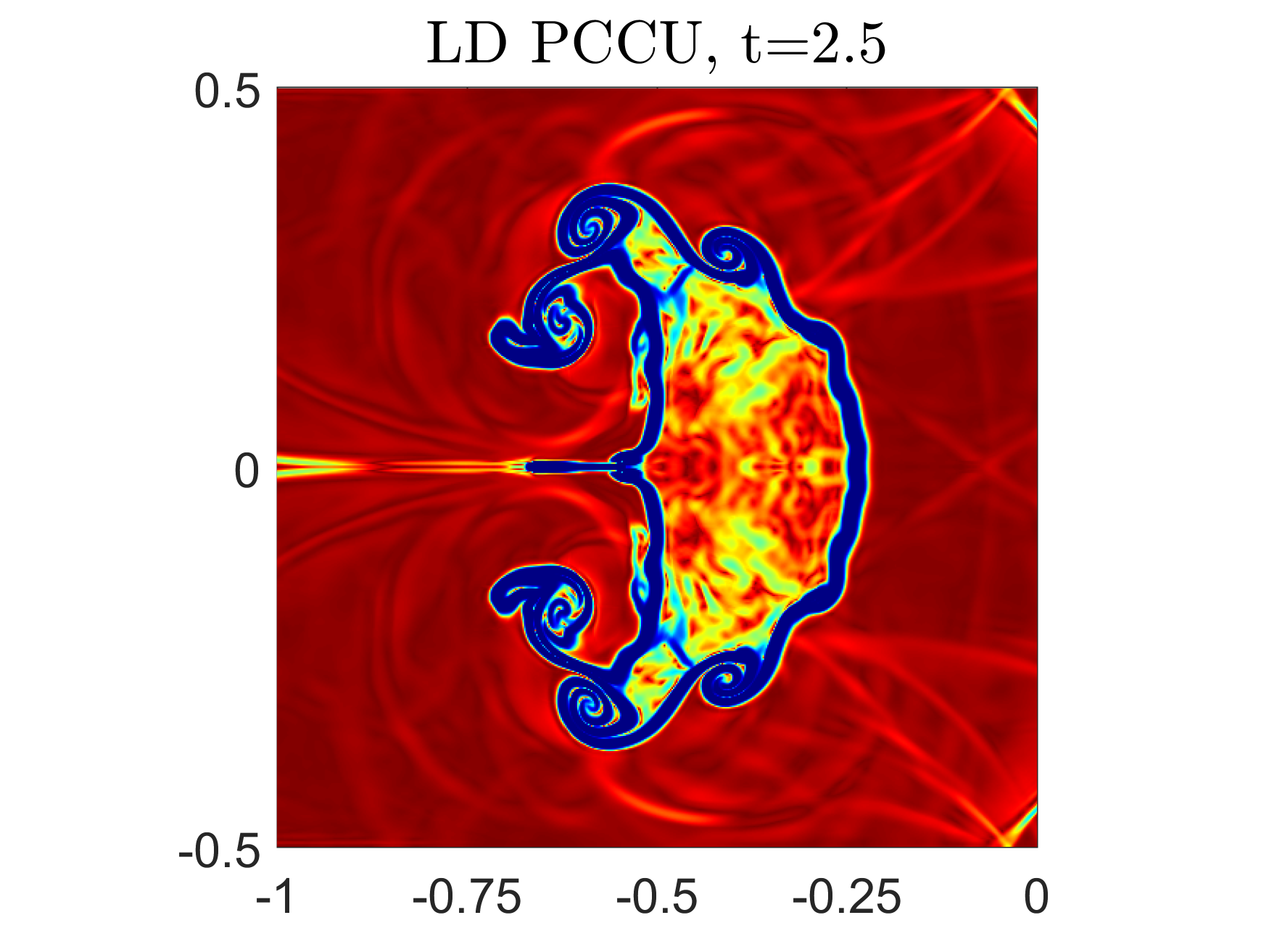}\hspace*{0.5cm}
            \includegraphics[trim=2.1cm 0.4cm 2.1cm 0.2cm, clip, width=5.1cm]{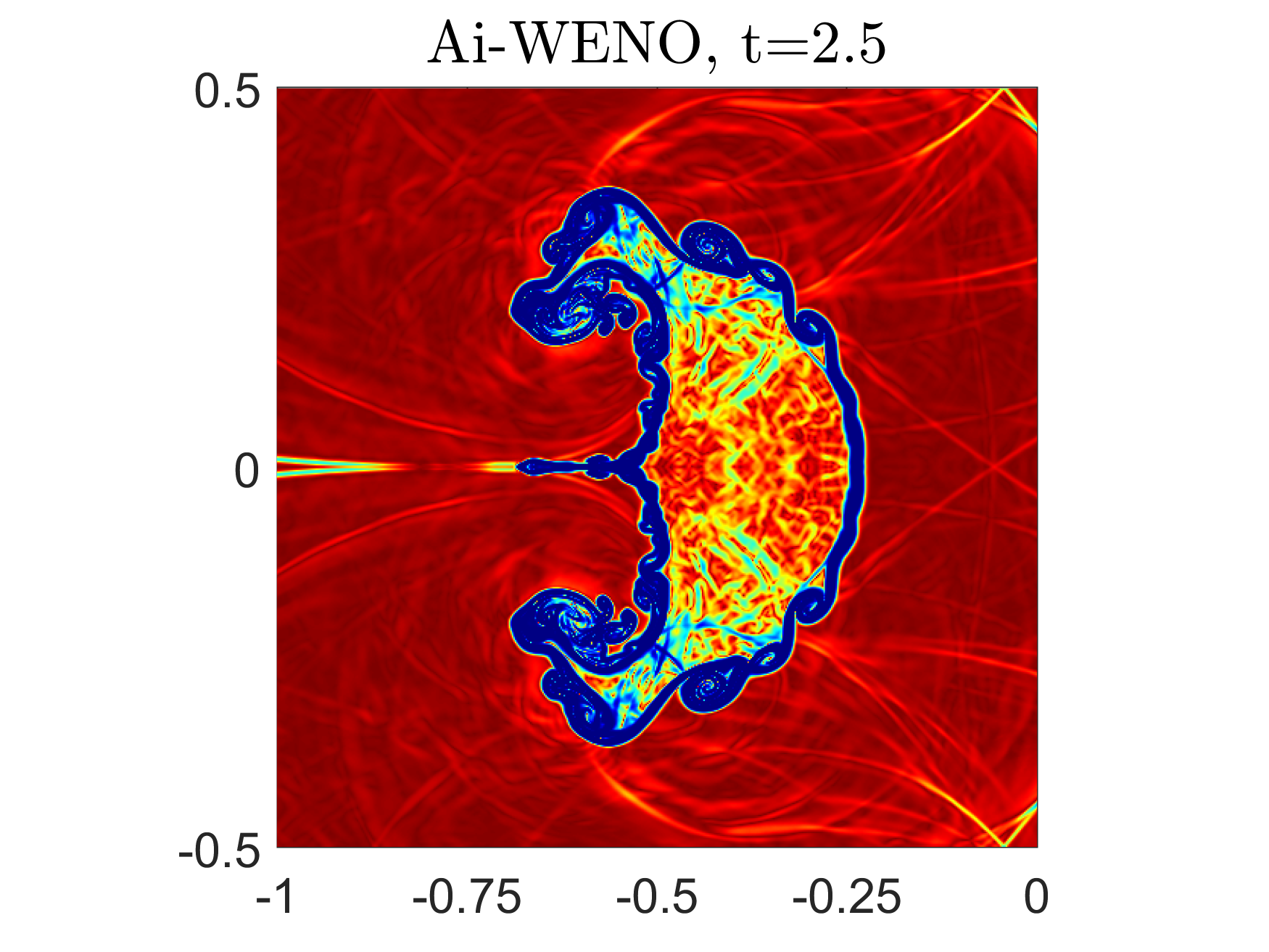}}
\vskip10pt
\centerline{\includegraphics[trim=2.1cm 0.4cm 2.1cm 0.2cm, clip, width=5.1cm]{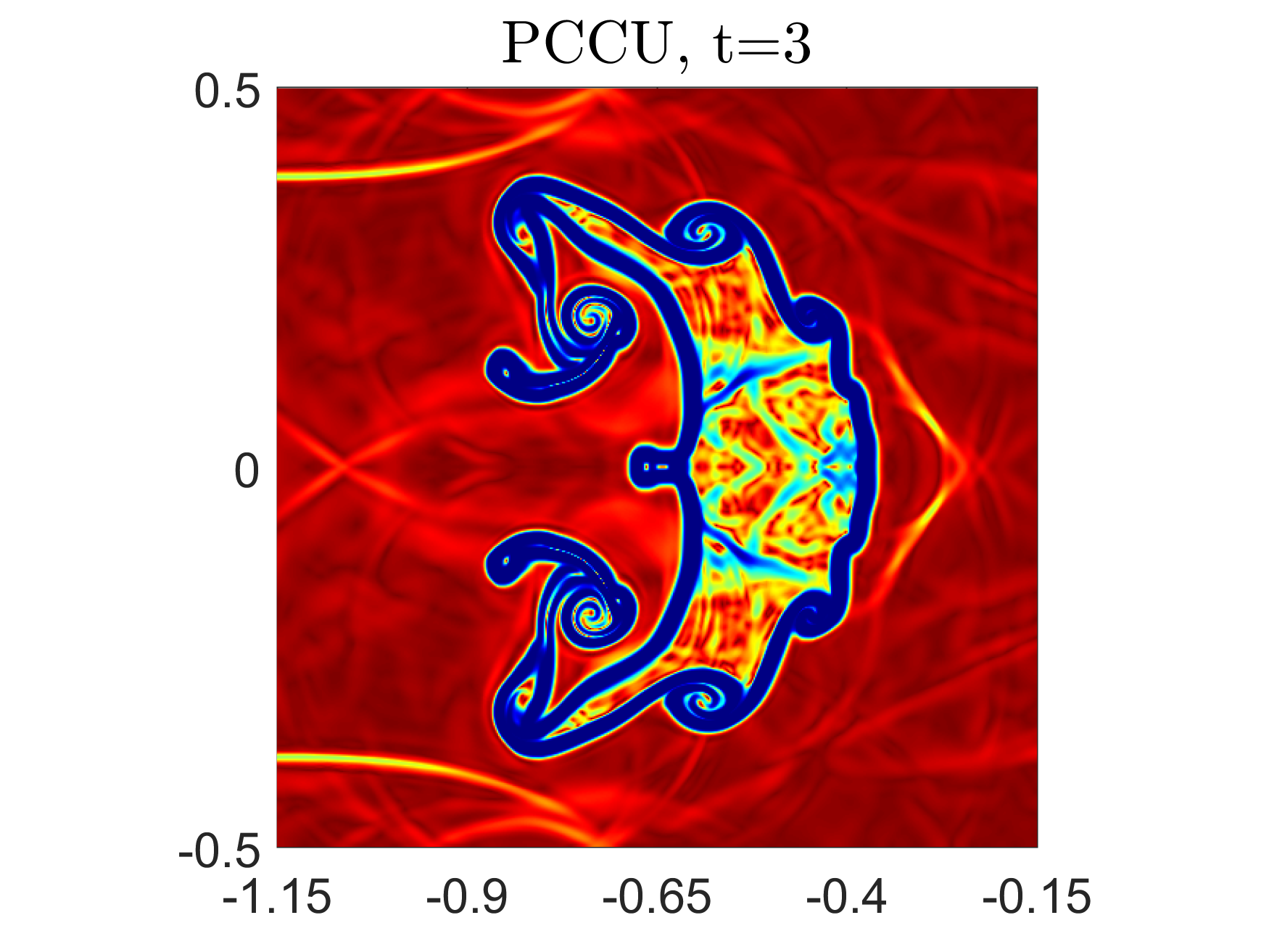}\hspace*{0.5cm}
            \includegraphics[trim=2.1cm 0.4cm 2.1cm 0.2cm, clip, width=5.1cm]{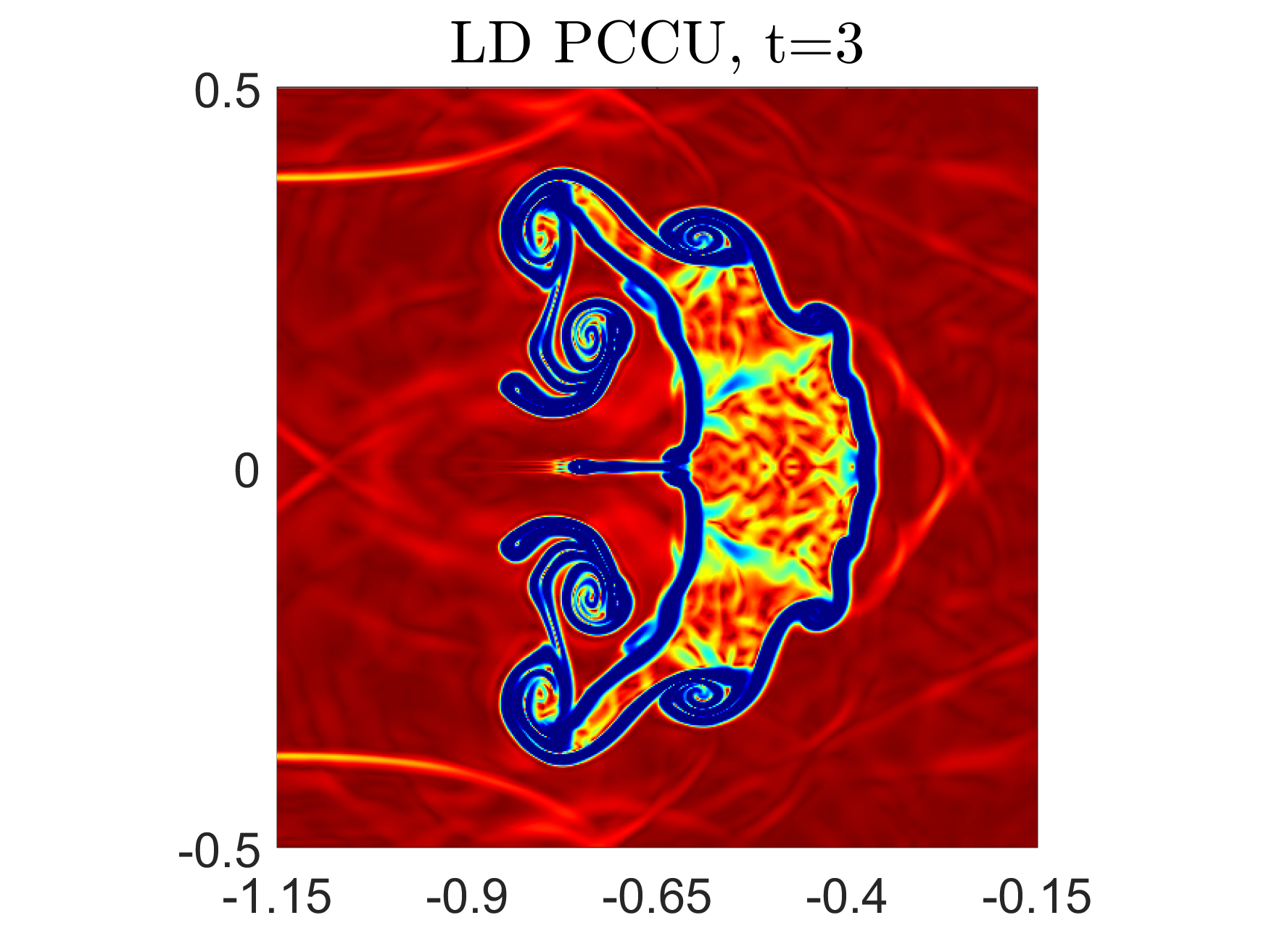}\hspace*{0.5cm}
            \includegraphics[trim=2.1cm 0.4cm 2.1cm 0.2cm, clip, width=5.1cm]{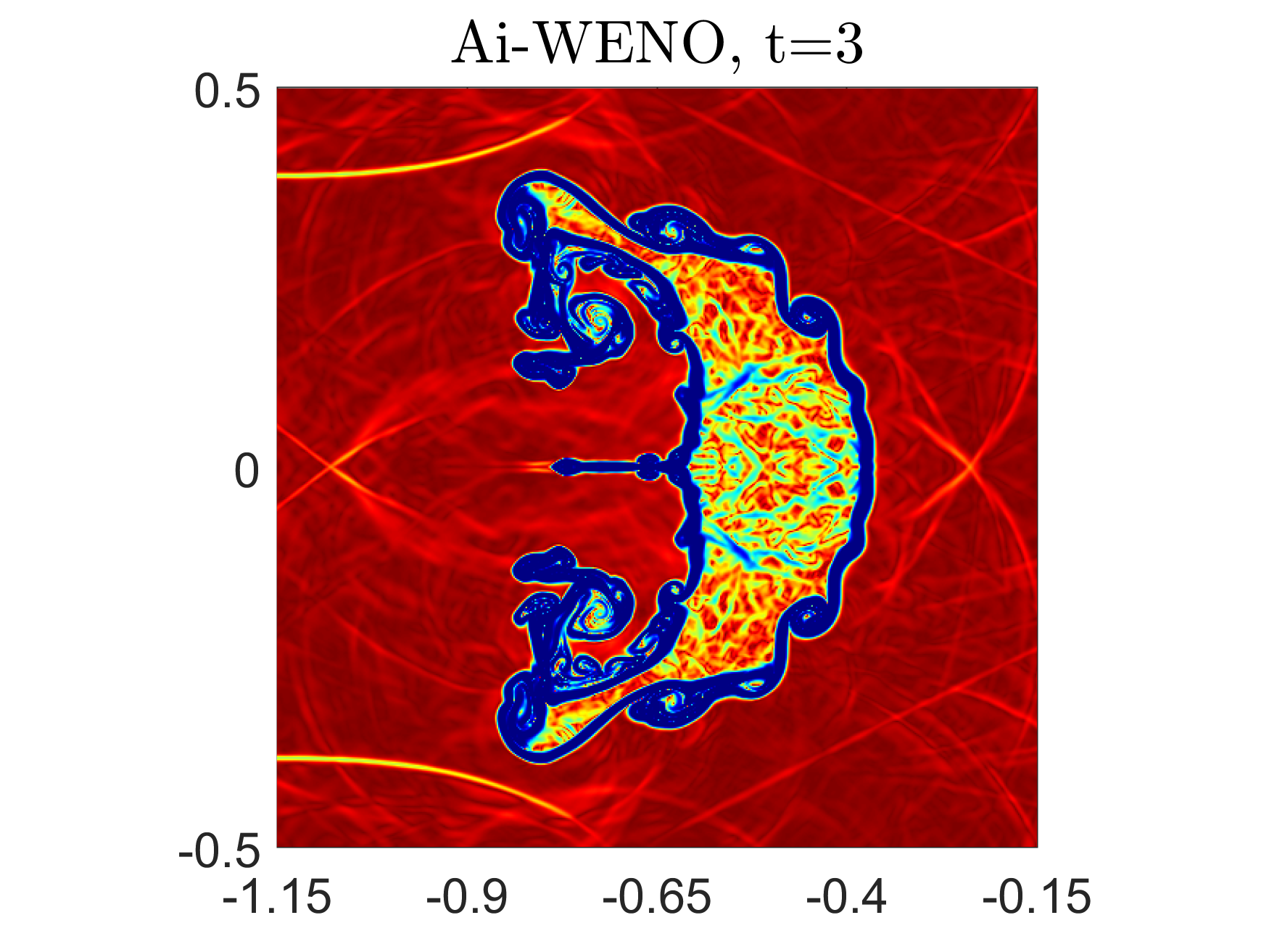}}
\caption{\sf Same as in Figure \ref{fig45}, but at larger times $t=1.5$, 2, 2.5, and 3.\label{fig46}}
\end{figure}

We compute the numerical solutions until the final time $t=3$ on a uniform mesh with $\dx=\dy=1/500$. In Figures \ref{fig45} and
\ref{fig46}, we present different stages of the shock-bubble interaction computed by the PCCU, LD PCCU, and Ai-WENO schemes. As one can see,
the bubble changes its shape and propagates to the left, and in order to focus on the details of the bubble structure, we only zoom at
$[\sigma,\sigma+1]\times[-0.5,0.5]$ square area containing the bubble ($\sigma$ is decreasing in time from -0.5 to -1.15). Compared with
Example 4, the bubble moves to the left a little slower and develops totally different structures as the R22 is heavier than the Helium. The
obtained results are in a good agreement with the numerical results reported in \cite{Quirk1996,chertock7,CCK_21}. Similar to Example 4, at
a small time $t=0.5$, the resolution of the bubble interface is significantly improved by the use of either the LD PCCU or Ai-WENO schemes,
and the improvement in this example is even more pronounced. By the time $t=1$, the interface develop instabilities which are smeared by a
more dissipative PCCU scheme. As time progresses, the solutions develop very complex small structures, which are better resolved by the
schemes containing smaller amount of numerical dissipation, namely, by the LD PCCU and Ai-WENO schemes.

\subsubsection*{Example 6---A Cylindrical Explosion Problem}
In this 2-D example, which is a modification of an example from \cite{XL17}, we consider the case where a cylindrical explosive source is
located between an air-water interface and an impermeable wall. The initial conditions are given by
\begin{equation*}
(\rho,u,v,p;\gamma,\pi_\infty)=\begin{cases}(1.27,0,0,8290;2,0),&\mbox{$(x-5)^2+(y-2)^2<1$},\\(0.02,0,0,1;1.4,0),&\mbox{$y>4$},\\
(1,0,0,1;7.15,3309),&\mbox{otherwise},\end{cases}
\end{equation*}
the solid wall boundary conditions are imposed at the bottom, and the free boundary conditions are prescribed on the other sides of the
computational domain $[0,10]\times[0,6]$.

Notice that the initial conditions contain three---not two---different fluids, and therefore we need to modify the way the fluid interfaces
are detected as the criteria \eref{3.7a} and \eref{3.7af} are applicable in the two-fluid case only. We first set
$\Gamma_{\rm I}=1/(\gamma_{\rm I}-1)$, $\Gamma_{\rm II}={1}/{(\gamma_{\rm II}-1)}$, and $\Gamma_{III}=1/(\gamma_{\rm III}-1)$, where
$\gamma_{\rm I}=2$, $\gamma_{\rm II}=1.4$, and $\gamma_{\rm III}=7.15$ are the specific heat ratios for the three fluids. We then introduce
$\widehat\Gamma_1:=(\Gamma_{\rm I}+\Gamma_{\rm III})/2$, $\widehat\Gamma_2:=(\Gamma_{\rm II}+\Gamma_{\rm III})/2$, and replace the
conditions \eref{3.7a} and \eref{3.7af} with
\begin{equation*}
(\xbar\Gamma_{j,k}-\widehat\Gamma_1)(\xbar\Gamma_{j+1,k}-\widehat\Gamma_1)<0\quad{\rm or}\quad
(\xbar\Gamma_{j,k}-\widehat\Gamma_2)(\xbar\Gamma_{j+1,k}-\widehat\Gamma_2)<0
\end{equation*}
and
\begin{equation*}
(\xbar\Gamma_{j,k}-\widehat\Gamma_1)(\xbar\Gamma_{j,k+1}-\widehat\Gamma_1)<0\quad{\rm or}\quad
(\xbar\Gamma_{j,k}-\widehat\Gamma_2)(\xbar\Gamma_{j,k+1}-\widehat\Gamma_2)<0
\end{equation*}
respectively.

\begin{figure}[ht!]
\centerline{\includegraphics[trim=1.2cm 1.5cm 1.1cm 1cm, clip, width=5.5cm]{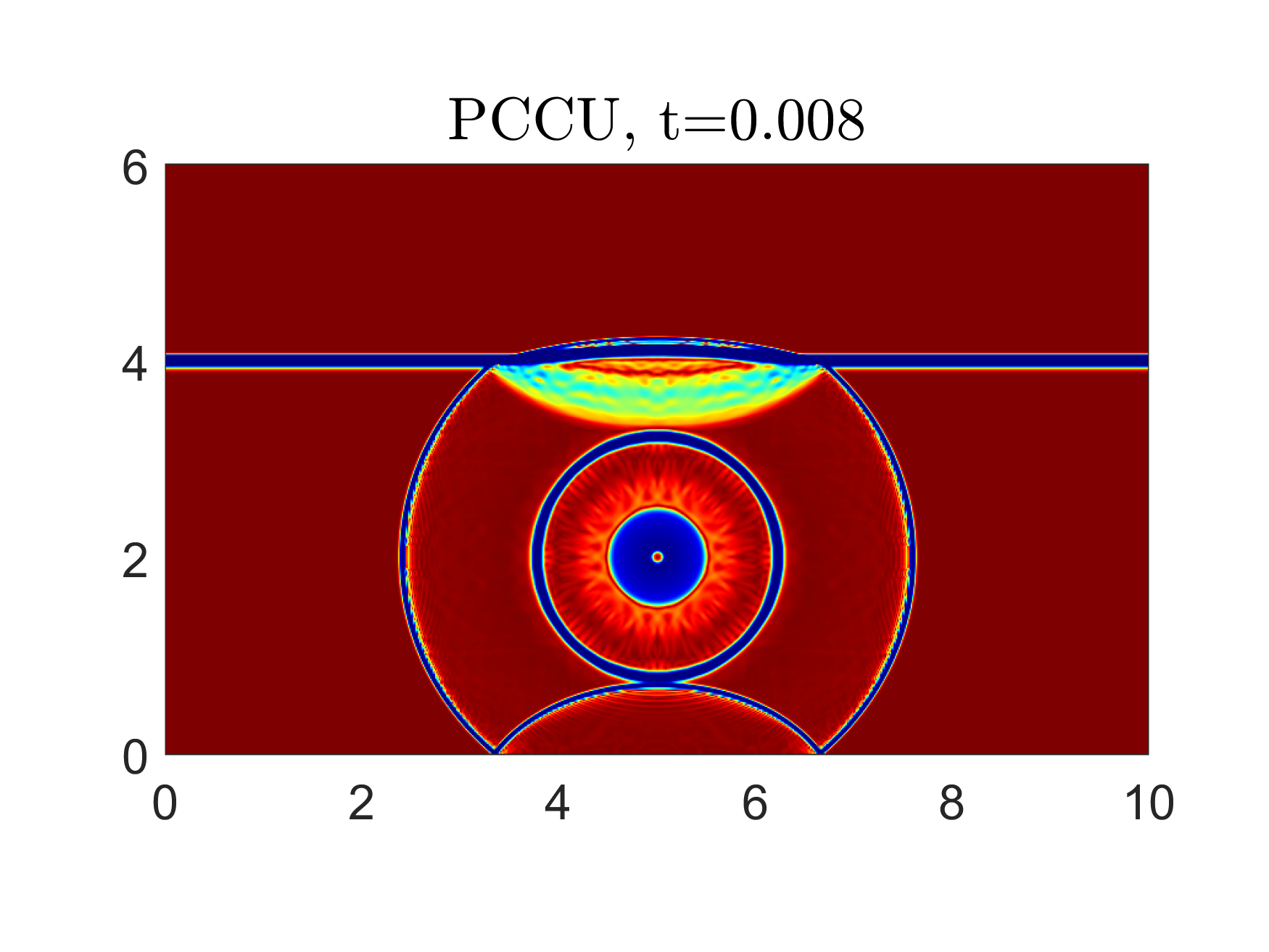}\hspace{0.2cm}
            \includegraphics[trim=1.2cm 1.5cm 1.1cm 1cm, clip, width=5.5cm]{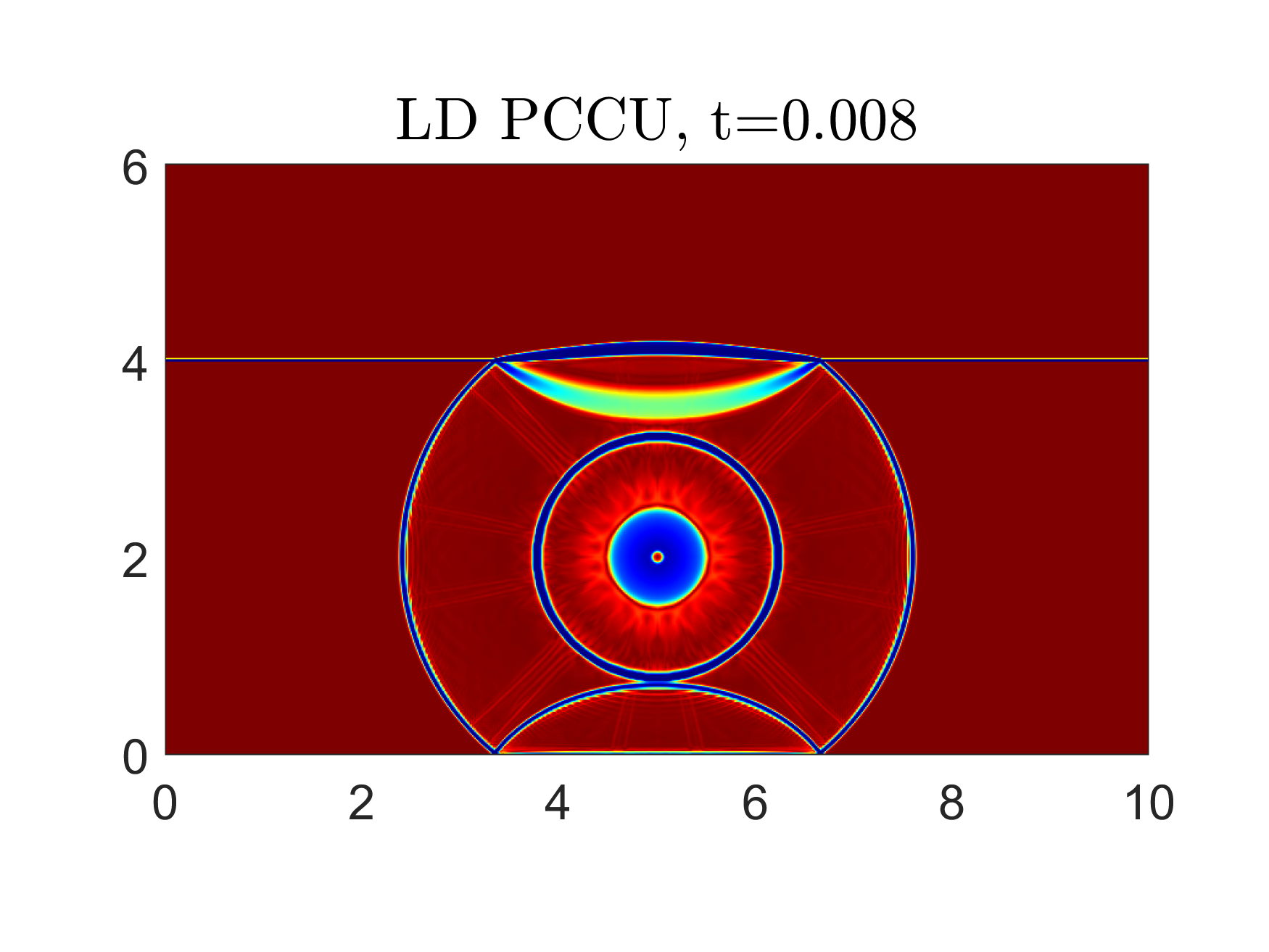}\hspace{0.2cm}
            \includegraphics[trim=1.2cm 1.5cm 1.1cm 1cm, clip, width=5.5cm]{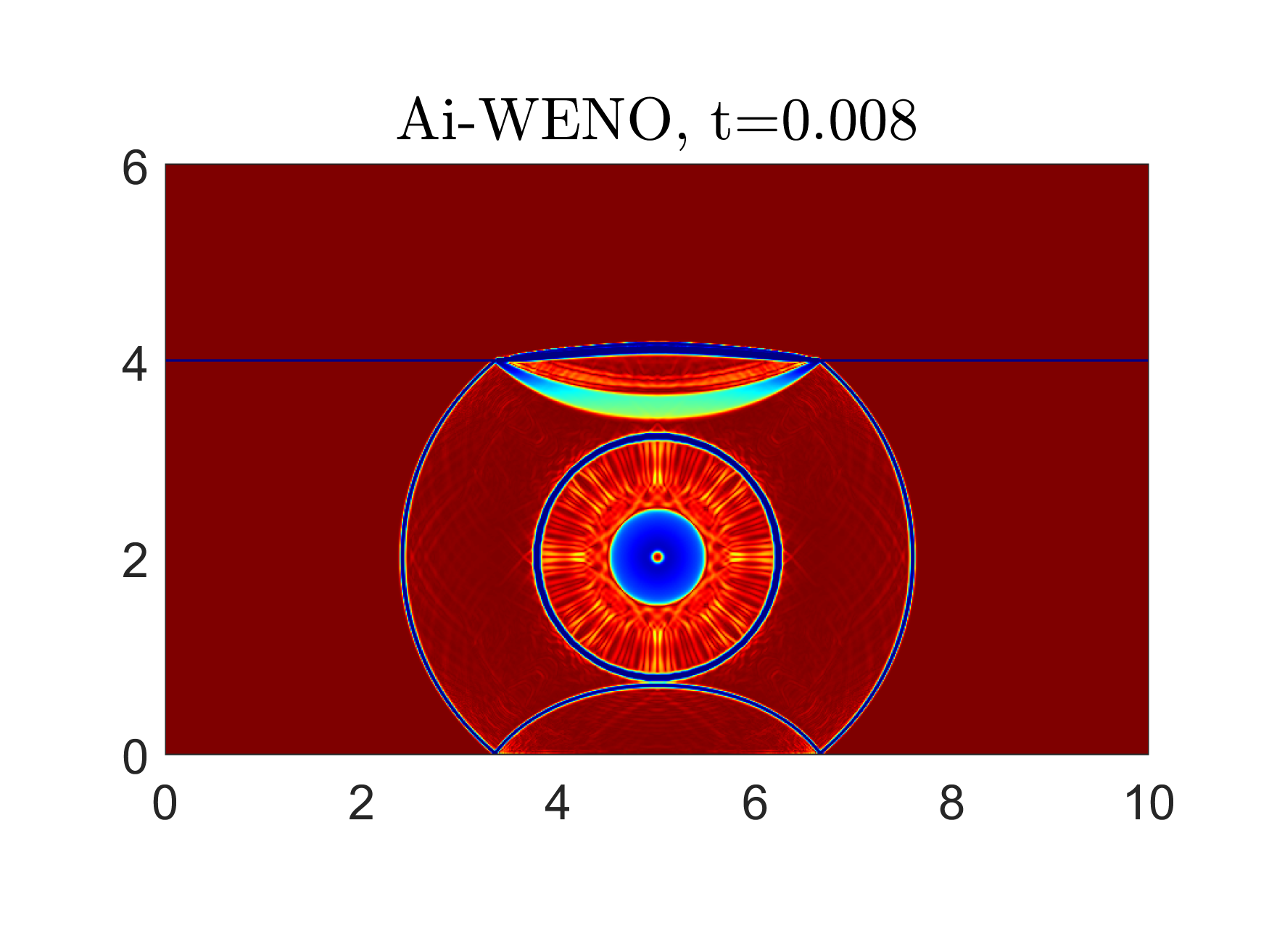}}
\vskip5pt
\centerline{\includegraphics[trim=1.2cm 1.5cm 1.1cm 1cm, clip, width=5.5cm]{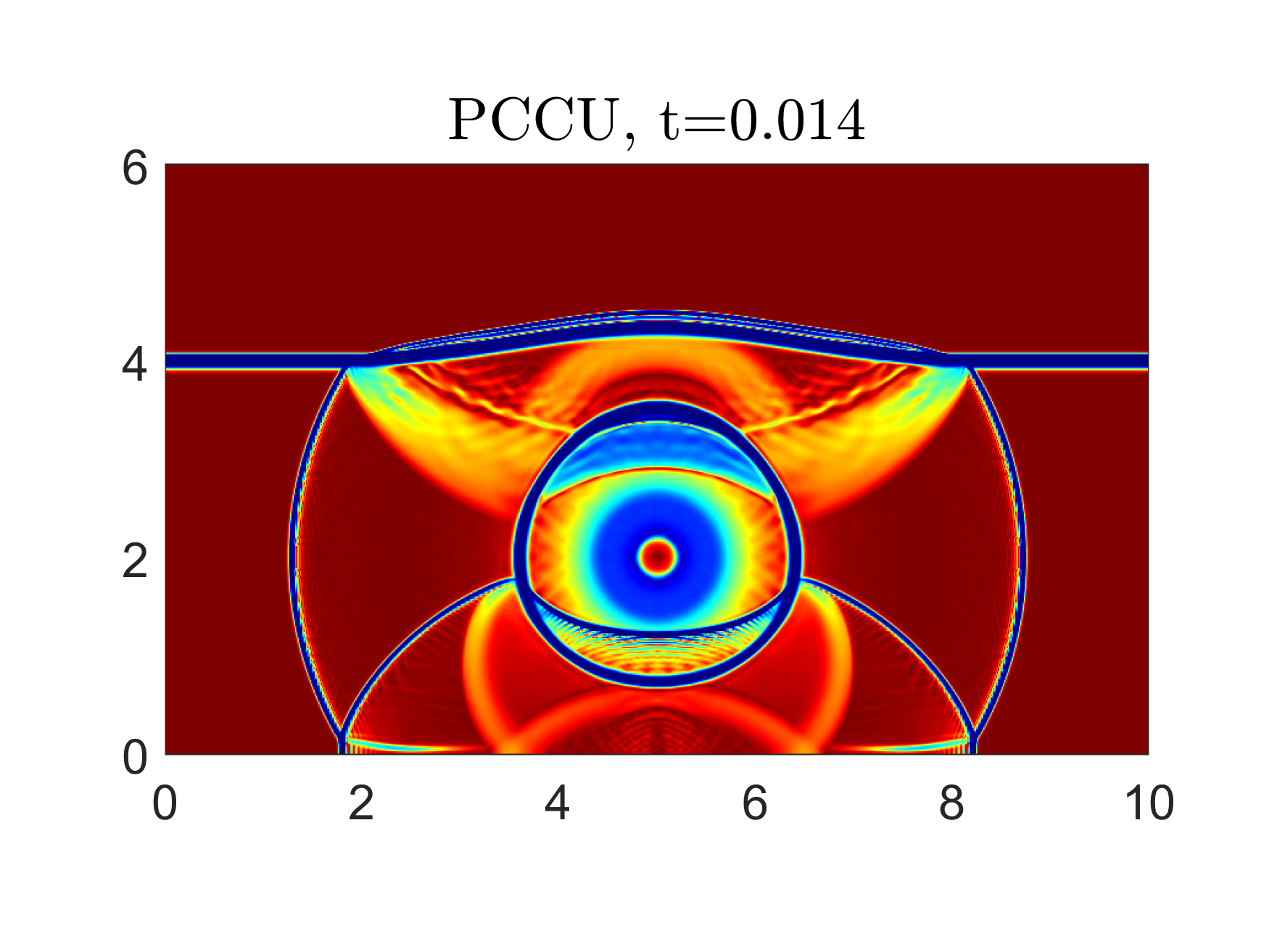}\hspace{0.2cm}
            \includegraphics[trim=1.2cm 1.5cm 1.1cm 1cm, clip, width=5.5cm]{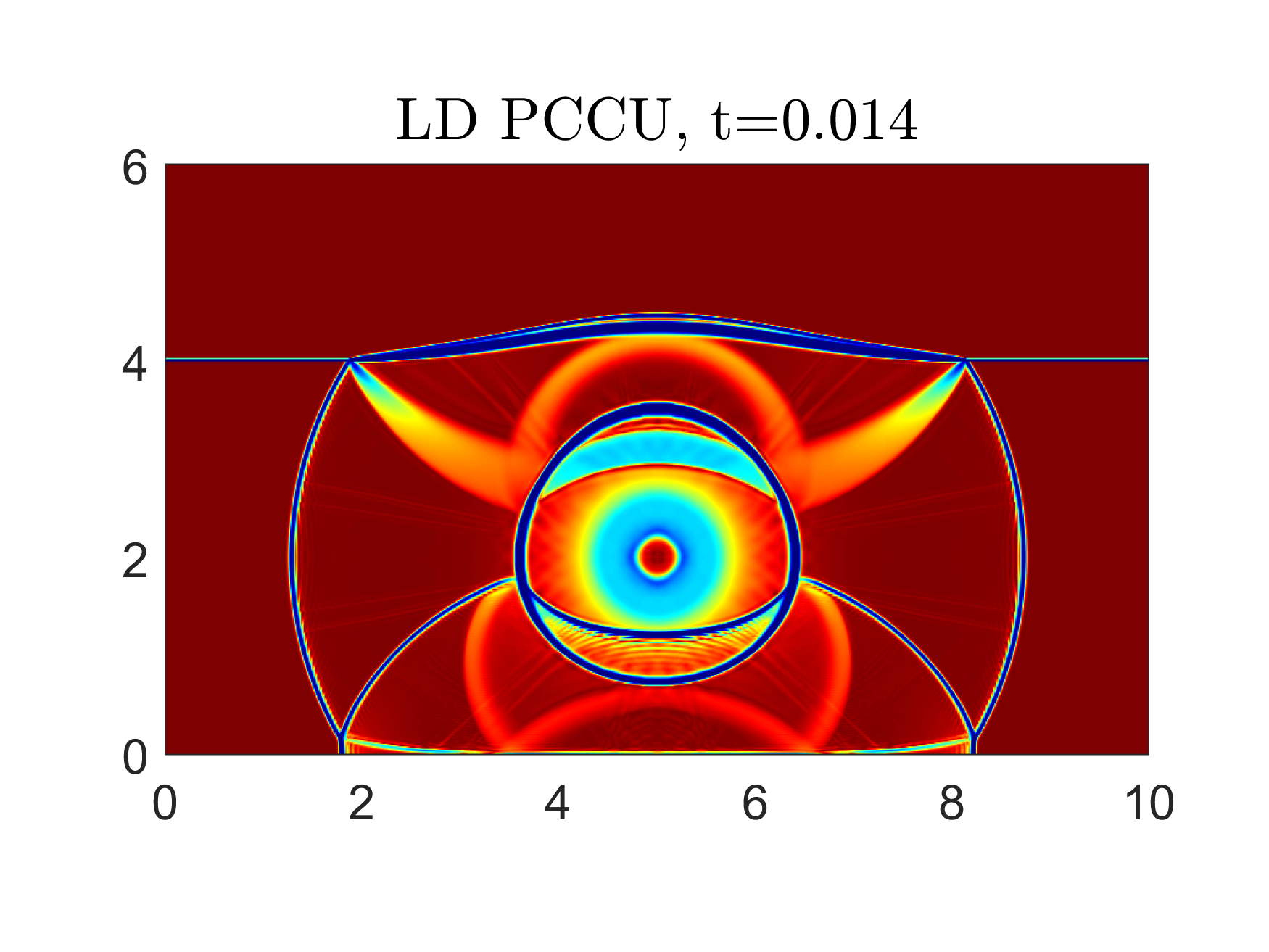}\hspace{0.2cm}
            \includegraphics[trim=1.2cm 1.5cm 1.1cm 1cm, clip, width=5.5cm]{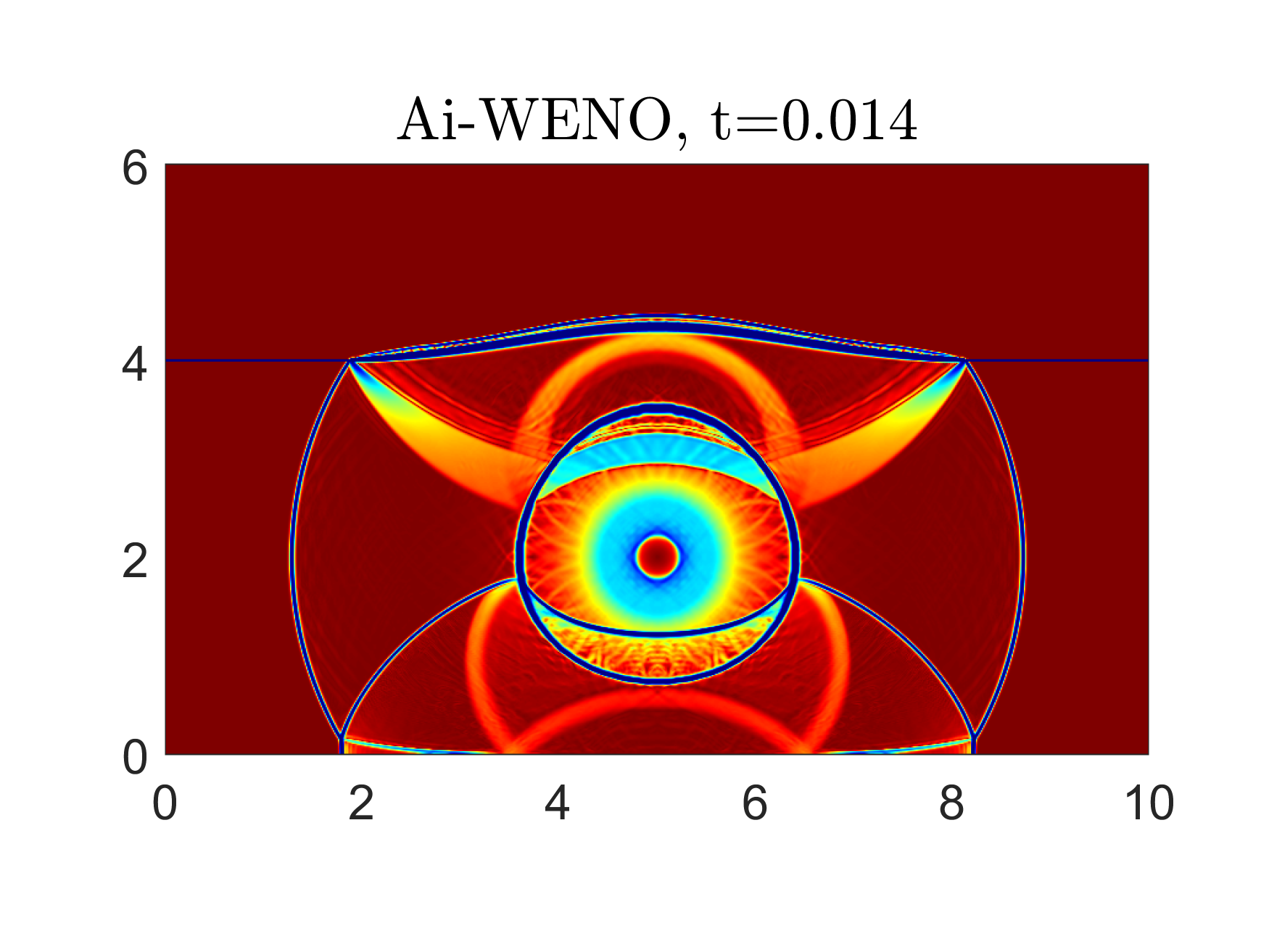}}
\vskip5pt
\centerline{\includegraphics[trim=1.2cm 1.5cm 1.1cm 1cm, clip, width=5.5cm]{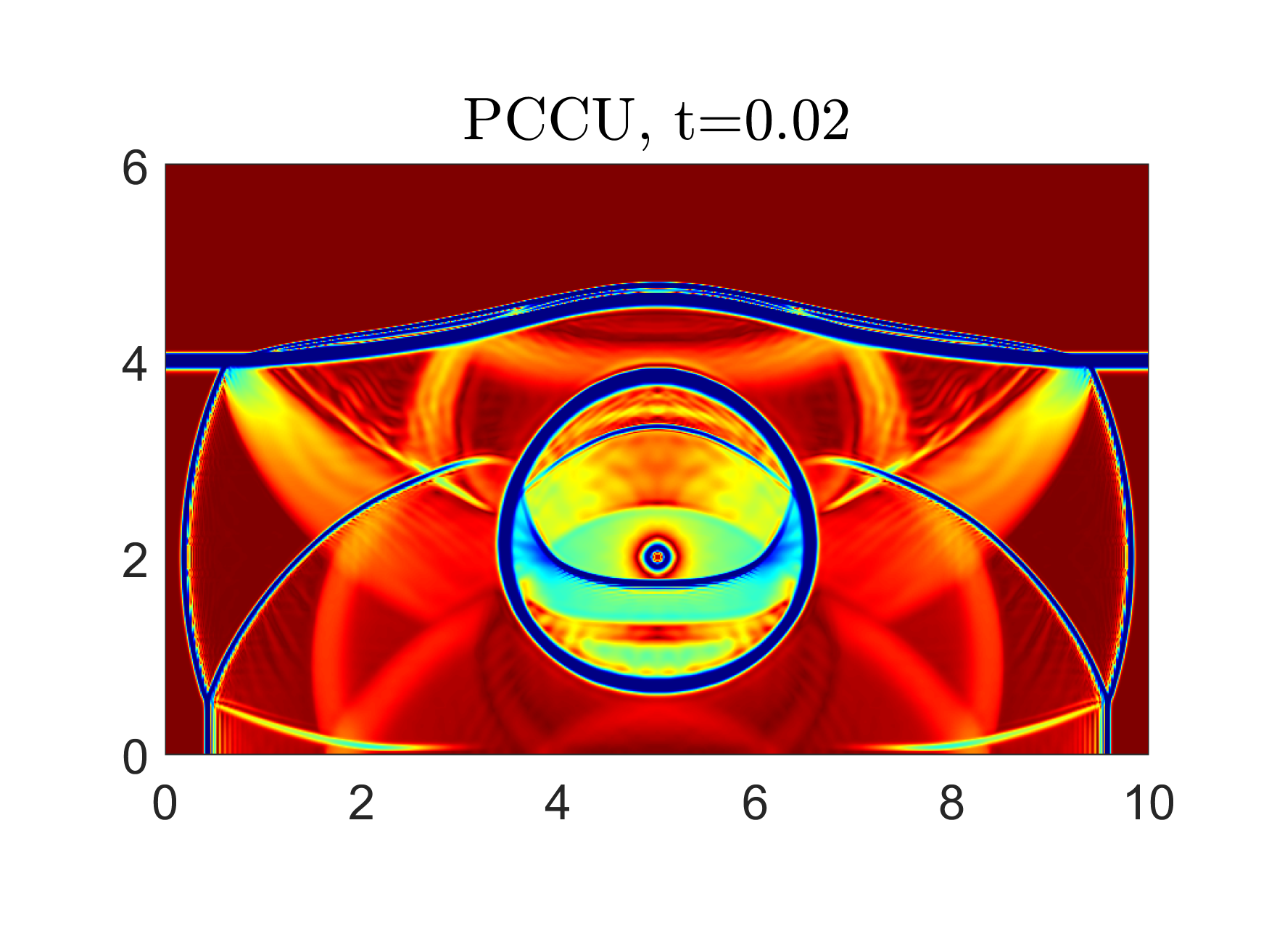}\hspace{0.2cm}
            \includegraphics[trim=1.2cm 1.5cm 1.1cm 1cm, clip, width=5.5cm]{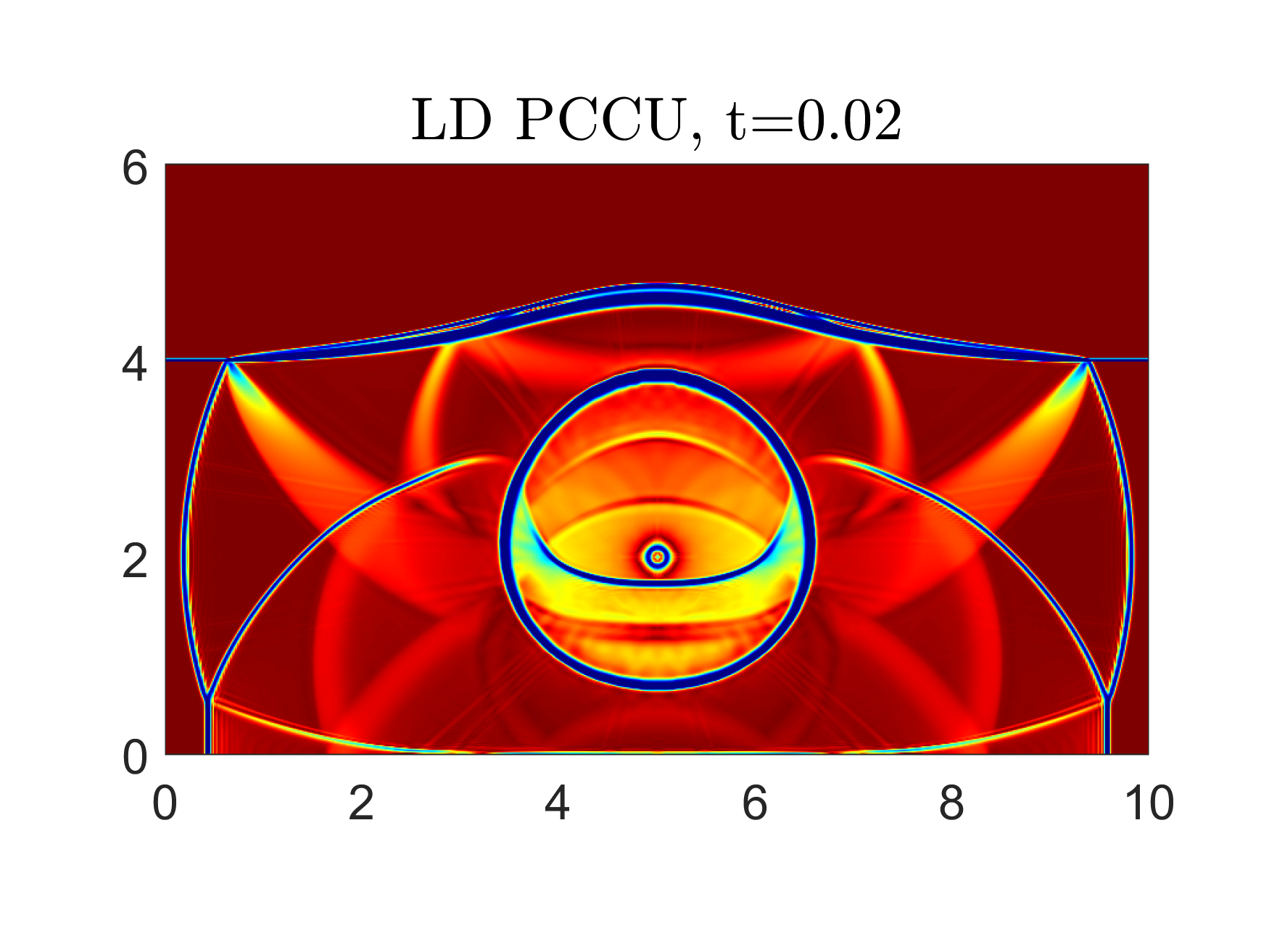}\hspace{0.2cm}
            \includegraphics[trim=1.2cm 1.5cm 1.1cm 1cm, clip, width=5.5cm]{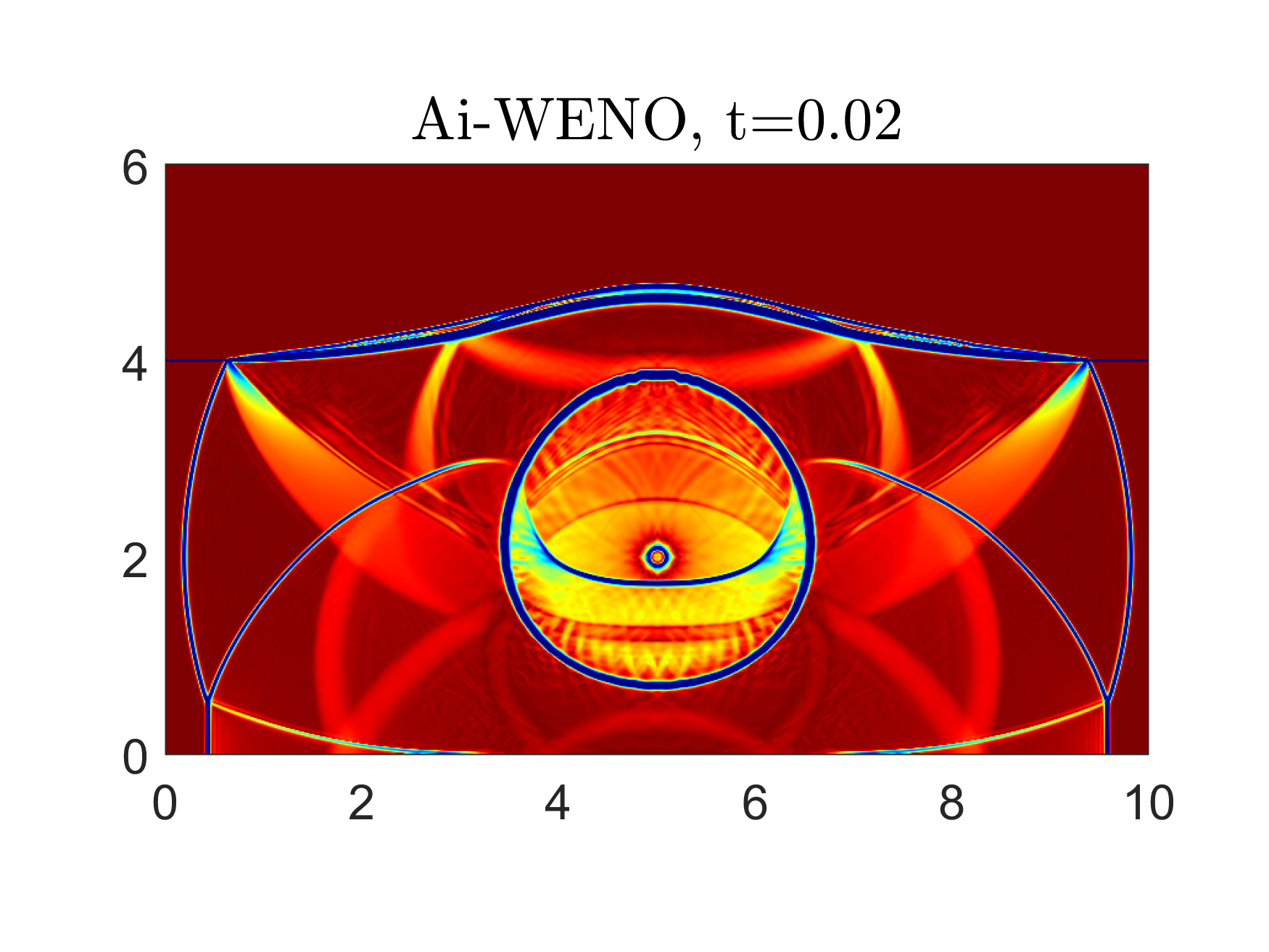}}
\caption{\sf Example 6: Solutions computed by the PCCU (left column), LD PCCU (middle column), and Ai-WENO (right column) schemes at times
$t=0.008$, 0.014, and $0.02$.\label{fig46a}}
\end{figure}

\begin{figure}[ht]
\centerline{\includegraphics[trim=2.1cm 0.4cm 1.9cm 0.2cm, clip, width=4.7cm]{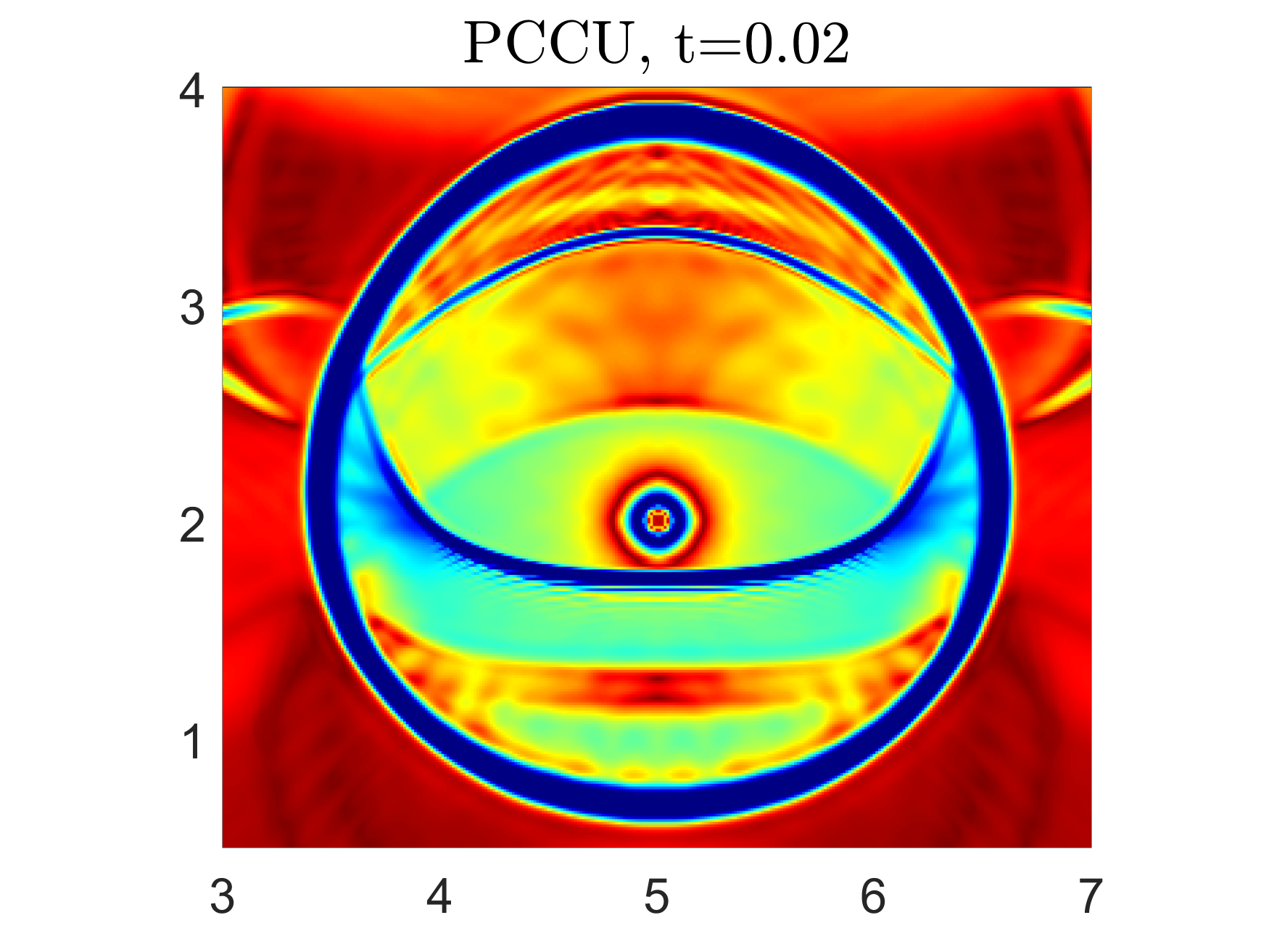}\hspace{0.5cm}
            \includegraphics[trim=2.1cm 0.4cm 1.9cm 0.2cm, clip, width=4.7cm]{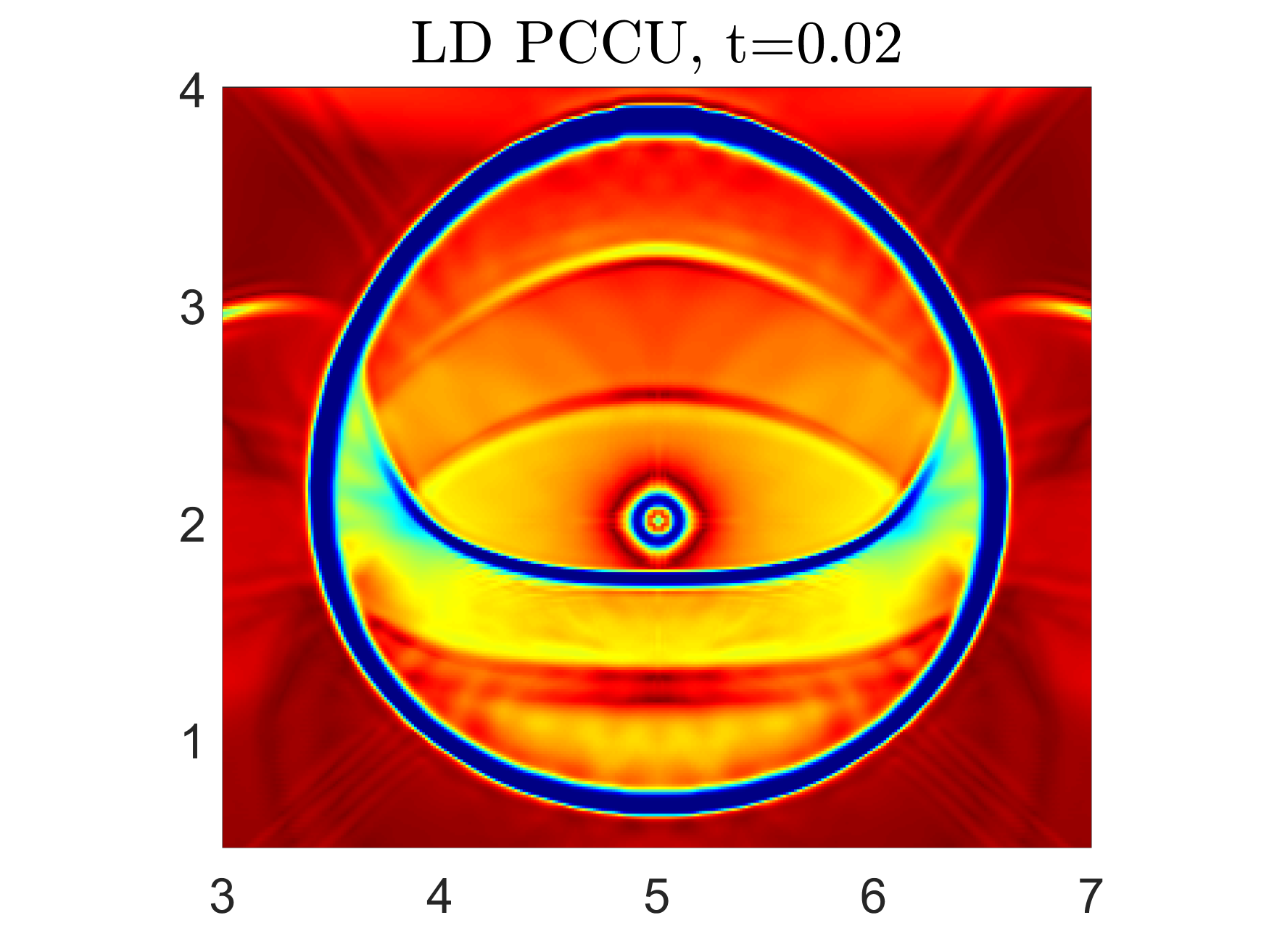}\hspace{0.5cm}
            \includegraphics[trim=2.1cm 0.4cm 1.9cm 0.2cm, clip, width=4.7cm]{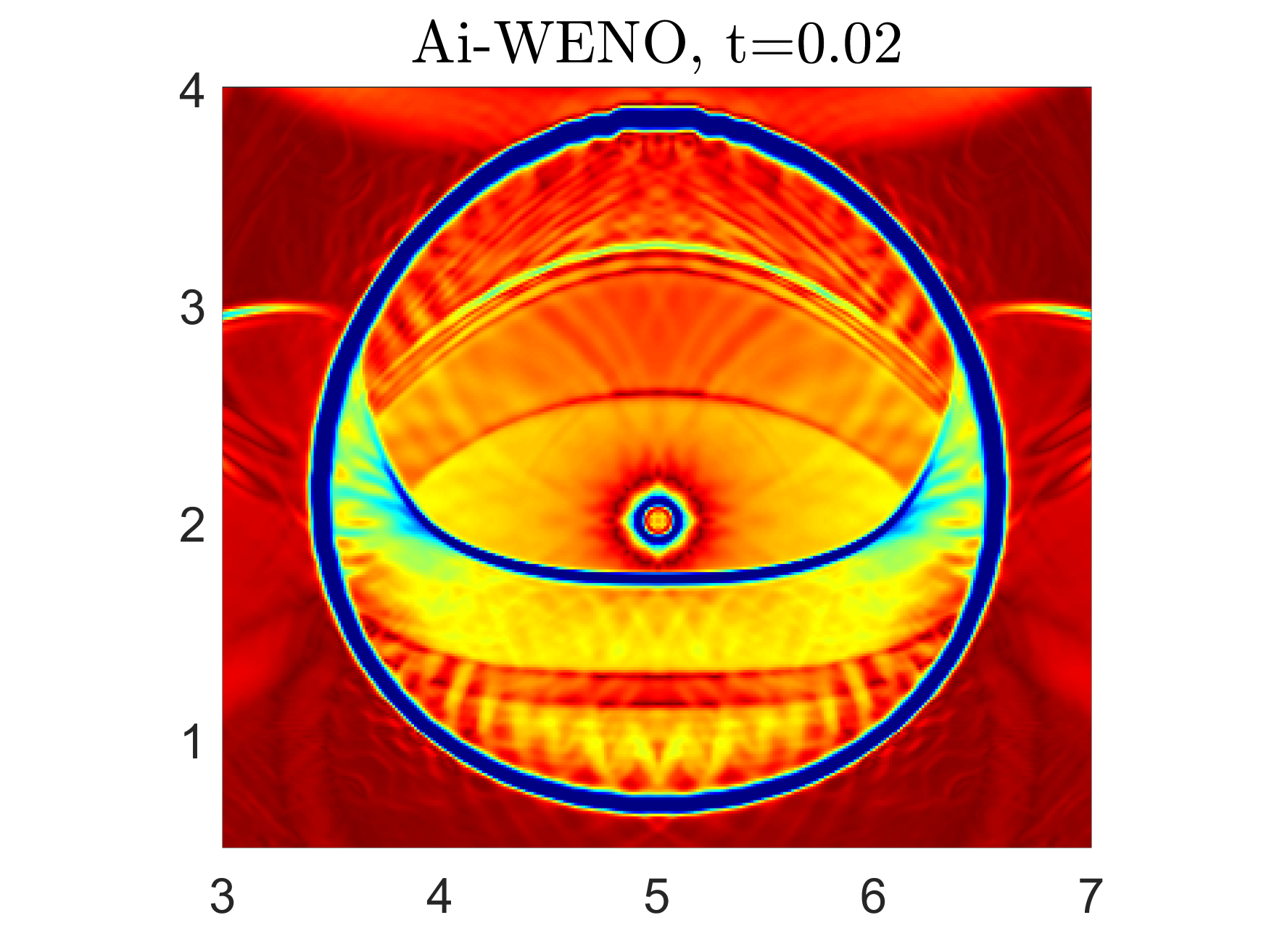}}
\caption{\sf Example 6: Solutions computed by the PCCU (left), LD PCCU (middle), and Ai-WENO (right) schemes at time $t=0.02$; zoom at the
bubble area.\label{fig46b}}
\end{figure}

We compute the numerical solutions until the final time $t=0.02$ on a uniform mesh with $\dx=\dy=1/80$ by the studied PCCU, LD PCCU, and
Ai-WENO schemes. In Figure \ref{fig46a}, we present time snapshots of the obtained results, which qualitatively look very similar to the
numerical results reported in \cite{XL17}. As one can see, the LD PCCU scheme captures both the material interfaces and many of the
developed wave structures substantially sharper than the PCCU scheme and the use of the Ai-WENO scheme further enhances the resolution. In
particular, one can observe more pronounced small structures in the Ai-WENO solution (especially inside the bubble) at the large time
$t=0.02$; see Figure \ref{fig46b}, where we zoom at the bubble area. This example, once again, clearly indicates that the LD PCCU and
Ai-WENO schemes outperform the PCCU one.

\subsubsection*{Example 7---Water-Air Shock-Bubble Interaction}
In the last 2-D example, which is taken from \cite{CZL20}, we consider the interaction of a shock in water with a gas bubble. The initial
conditions
\begin{equation*}
(\rho,u,v,p;\gamma,\pi_\infty)=\begin{cases}(0.0012,0,0,1;1.4,0),&\mbox{$(x-6)^2+(y-6)^2<9$},\\
(1.325,-68.525,0,19153;4.4,6000),&\mbox{$x>11.4$},\\(1,0,0,1;4.4,6000),&\mbox{otherwise},\end{cases}
\end{equation*}
correspond to a cylindrical air bubble impacted by a Mach 1.72 shock initiated in water. In this example, the initial data are prescribed in
the computational domain $[0,12]\times[0,12]$ with the solid wall boundary conditions imposed on the top and bottom and the free boundary
conditions on the left and right edges of the computational domain.

\begin{figure}[ht!]
 \centerline{\includegraphics[trim=2.5cm 0.4cm 2.4cm 0.2cm, clip, width=4.7cm]{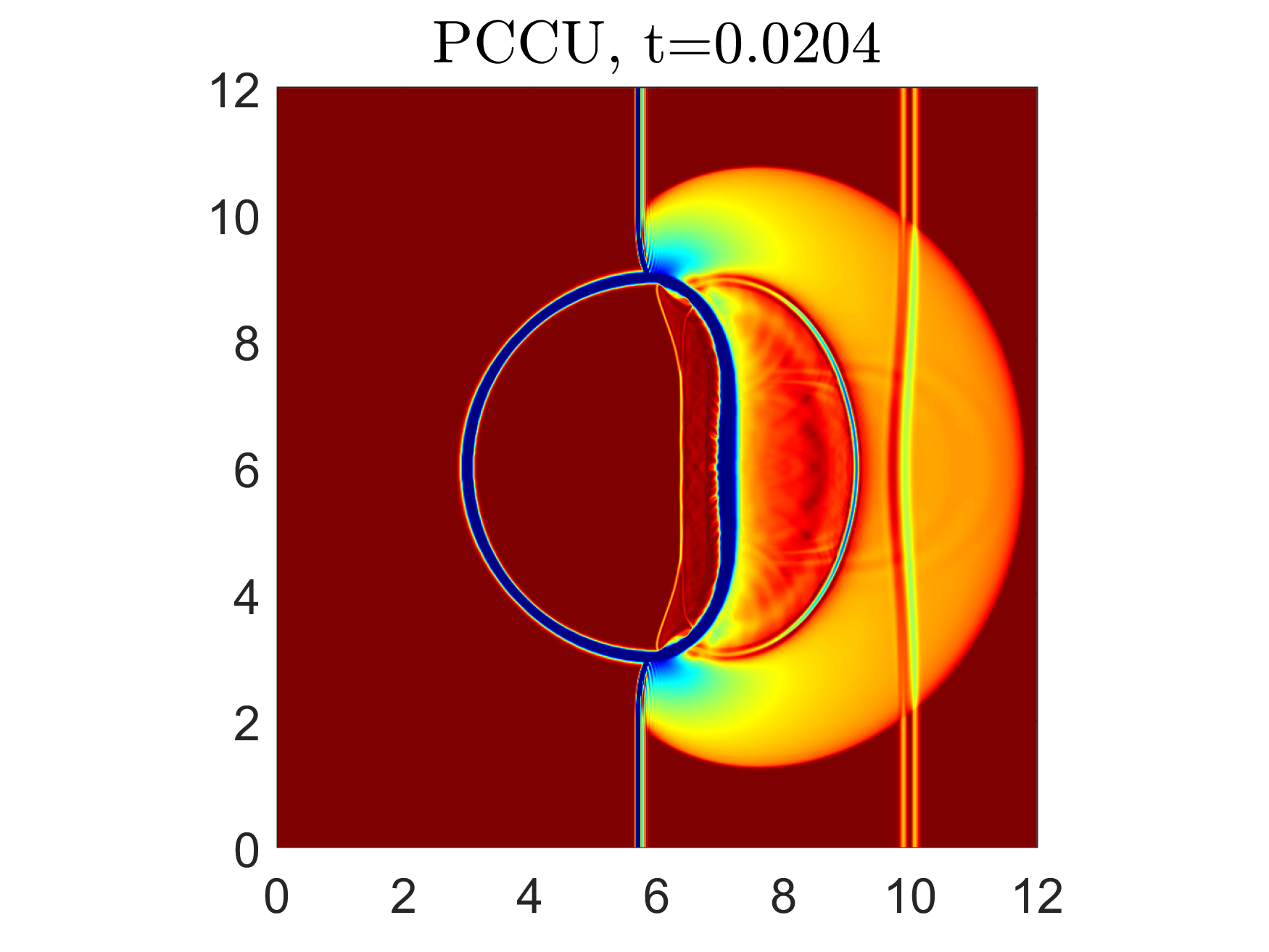}\hspace{0.5cm}
             \includegraphics[trim=2.5cm 0.4cm 2.4cm 0.2cm, clip, width=4.7cm]{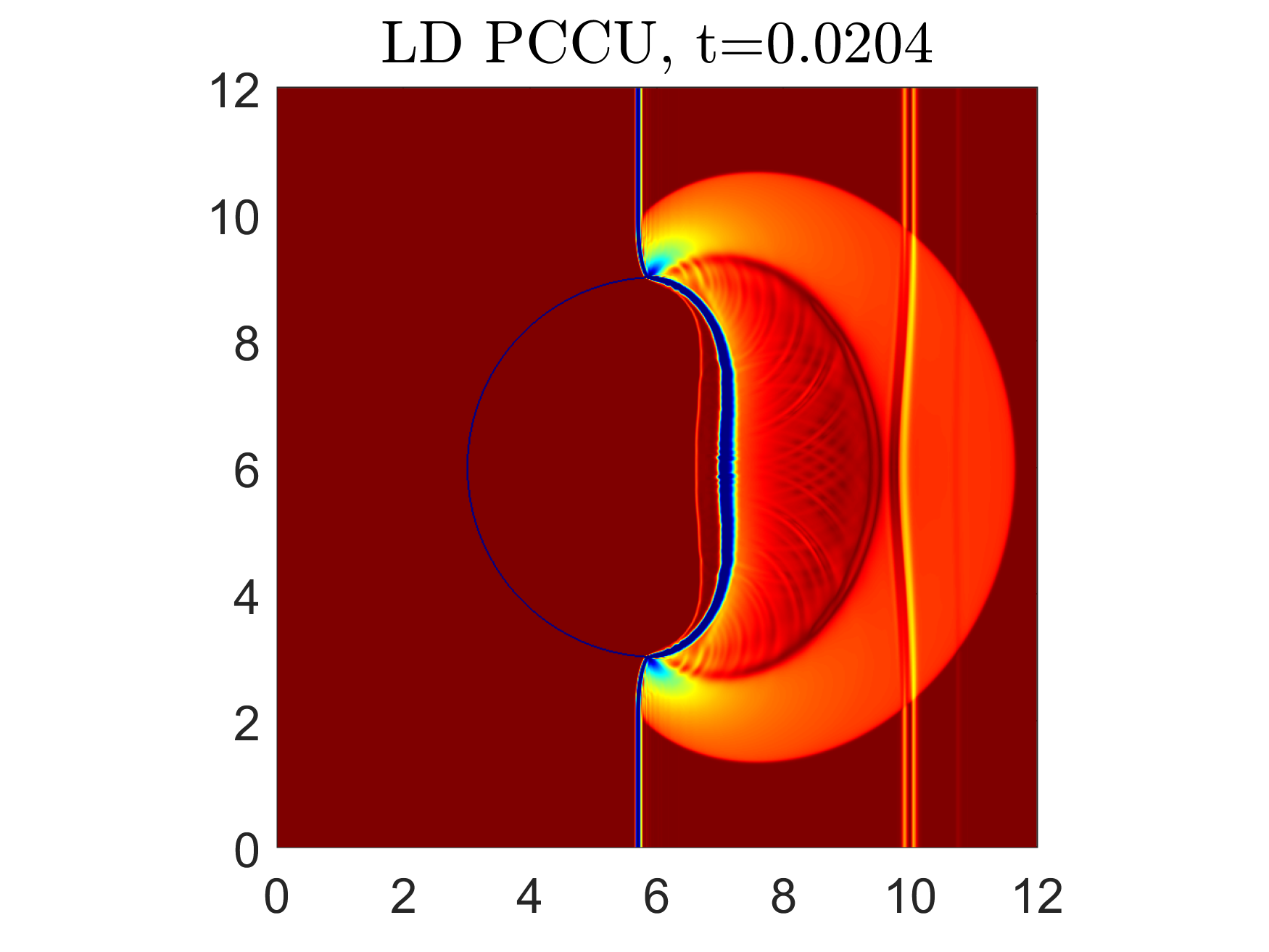}\hspace{0.5cm}
             \includegraphics[trim=2.5cm 0.4cm 2.4cm 0.2cm, clip, width=4.7cm]{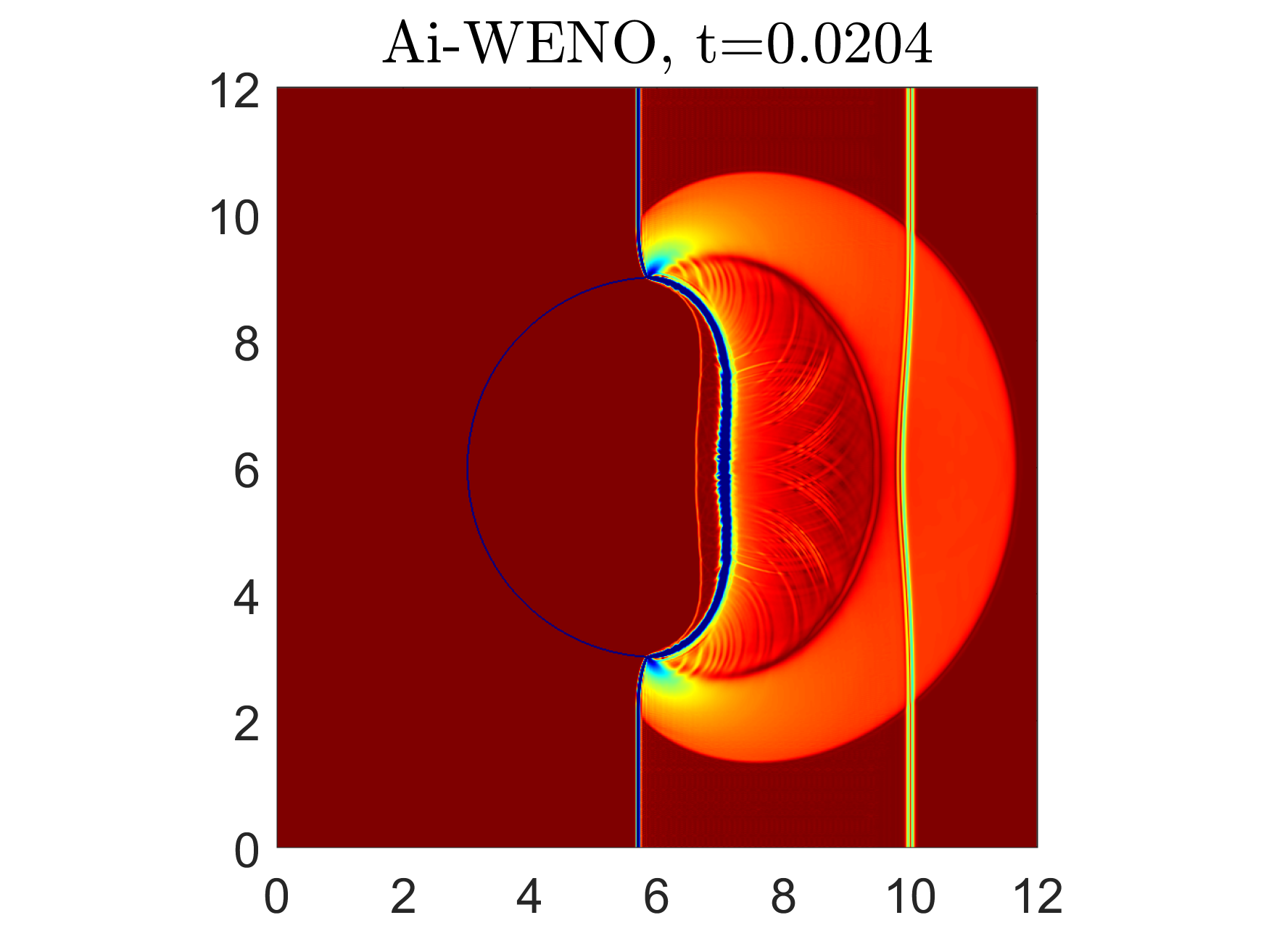}}
\vskip10pt
\centerline{\includegraphics[trim=2.5cm 0.4cm 2.4cm 0.2cm, clip, width=4.7cm]{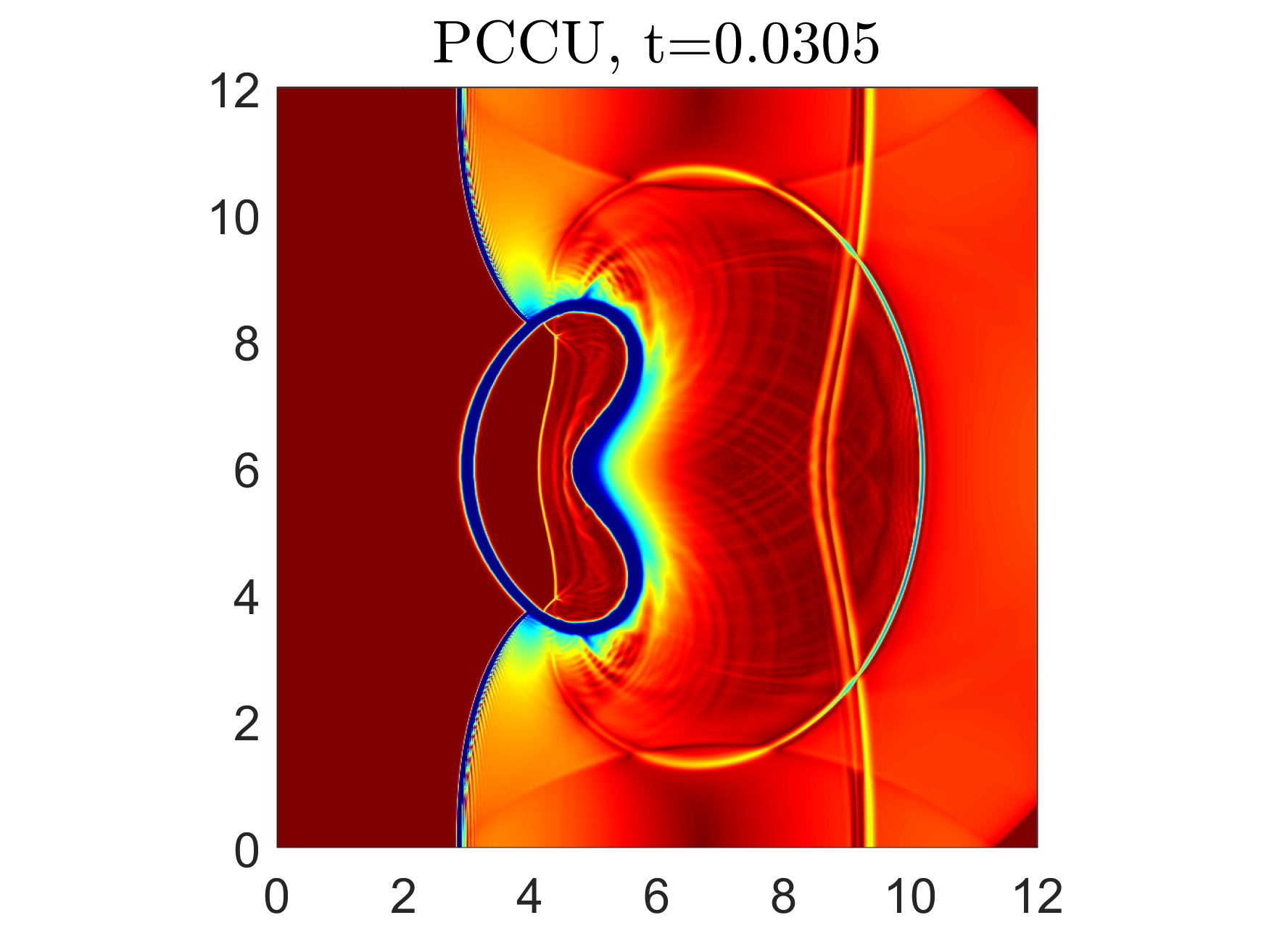}\hspace{0.5cm}
            \includegraphics[trim=2.5cm 0.4cm 2.4cm 0.2cm, clip, width=4.7cm]{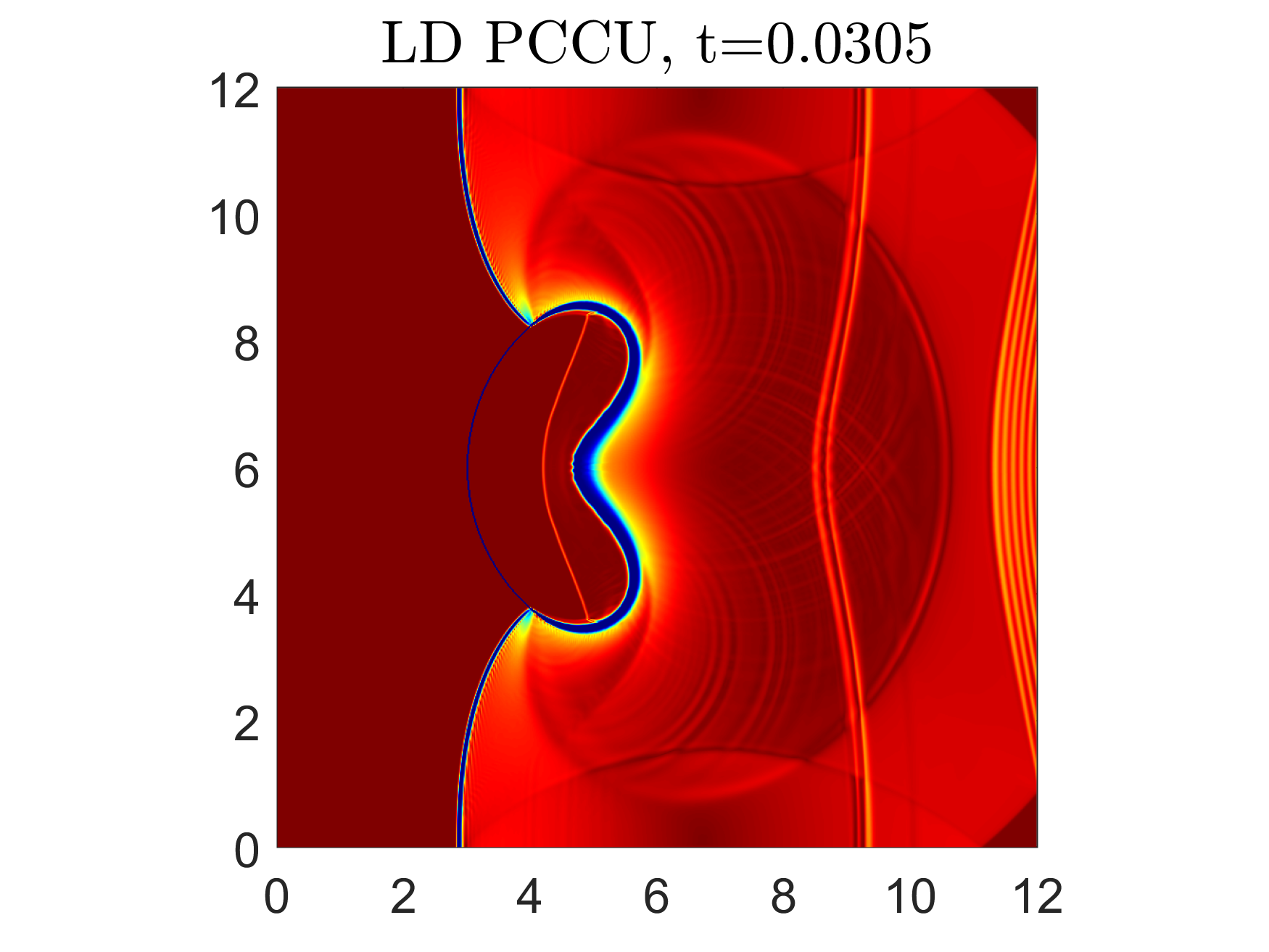}\hspace{0.5cm}
            \includegraphics[trim=2.5cm 0.4cm 2.4cm 0.2cm, clip, width=4.7cm]{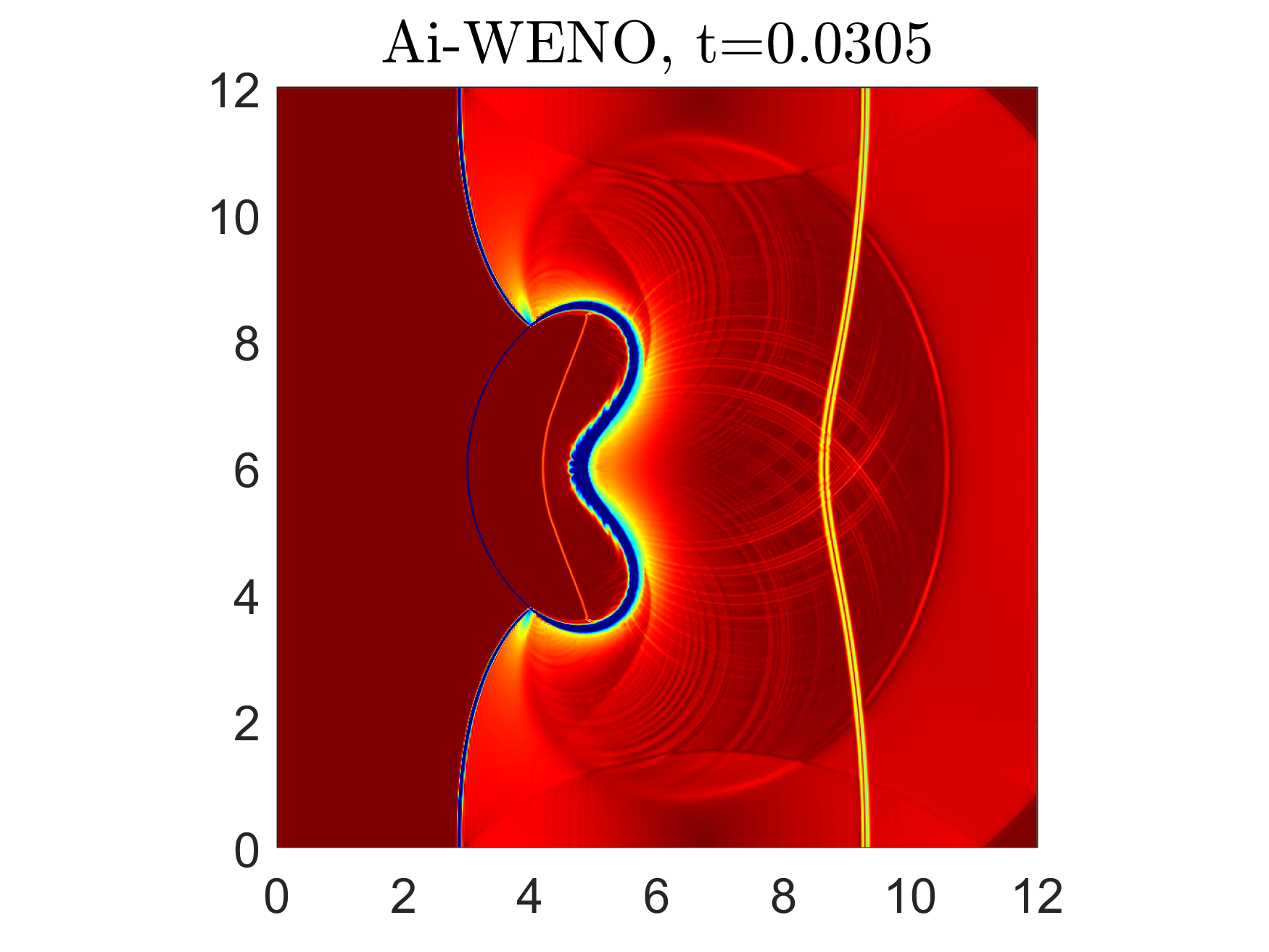}}
\vskip10pt
\centerline{\includegraphics[trim=2.5cm 0.4cm 2.4cm 0.2cm, clip, width=4.7cm]{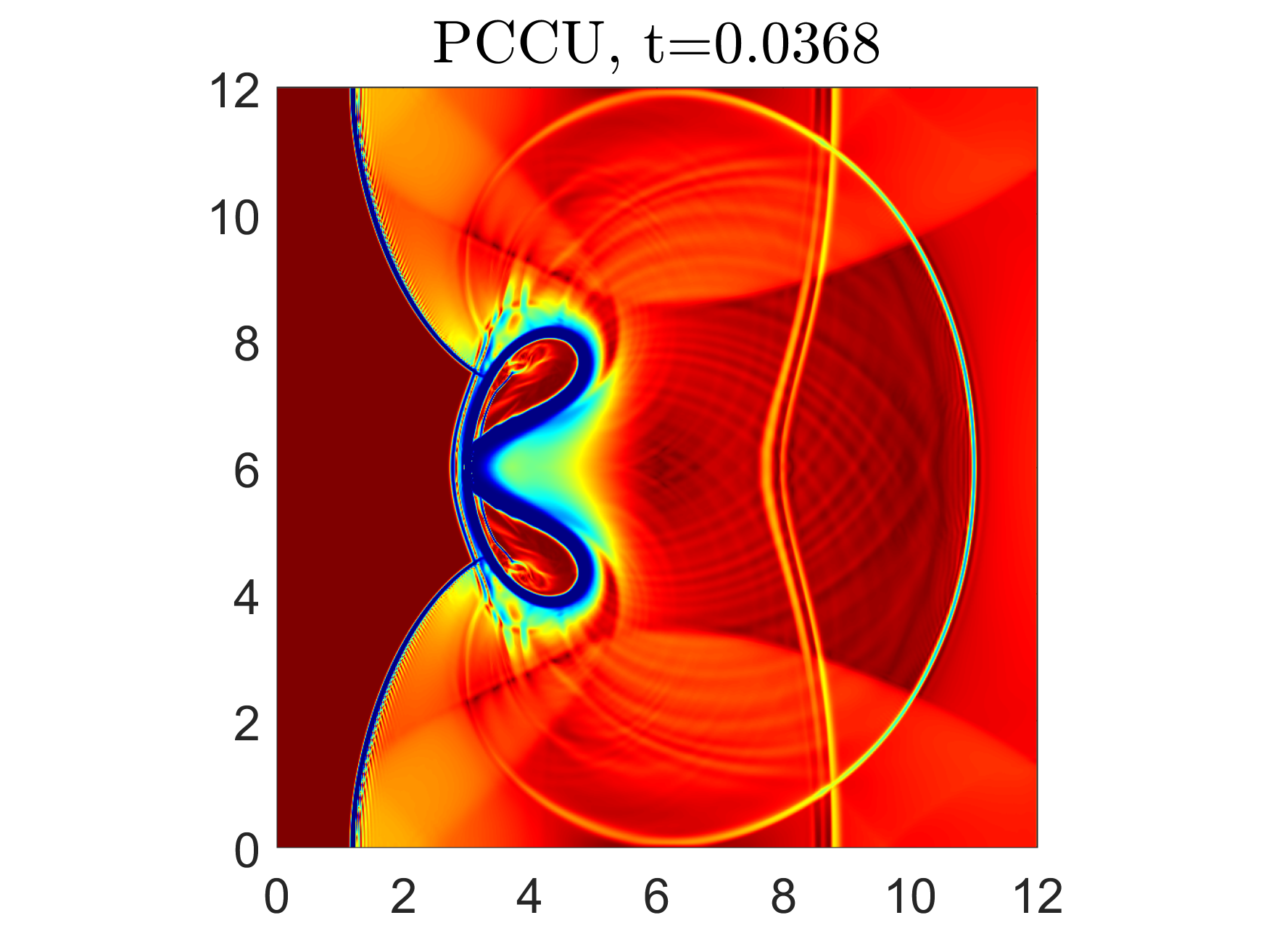}\hspace{0.5cm}
            \includegraphics[trim=2.5cm 0.4cm 2.4cm 0.2cm, clip, width=4.7cm]{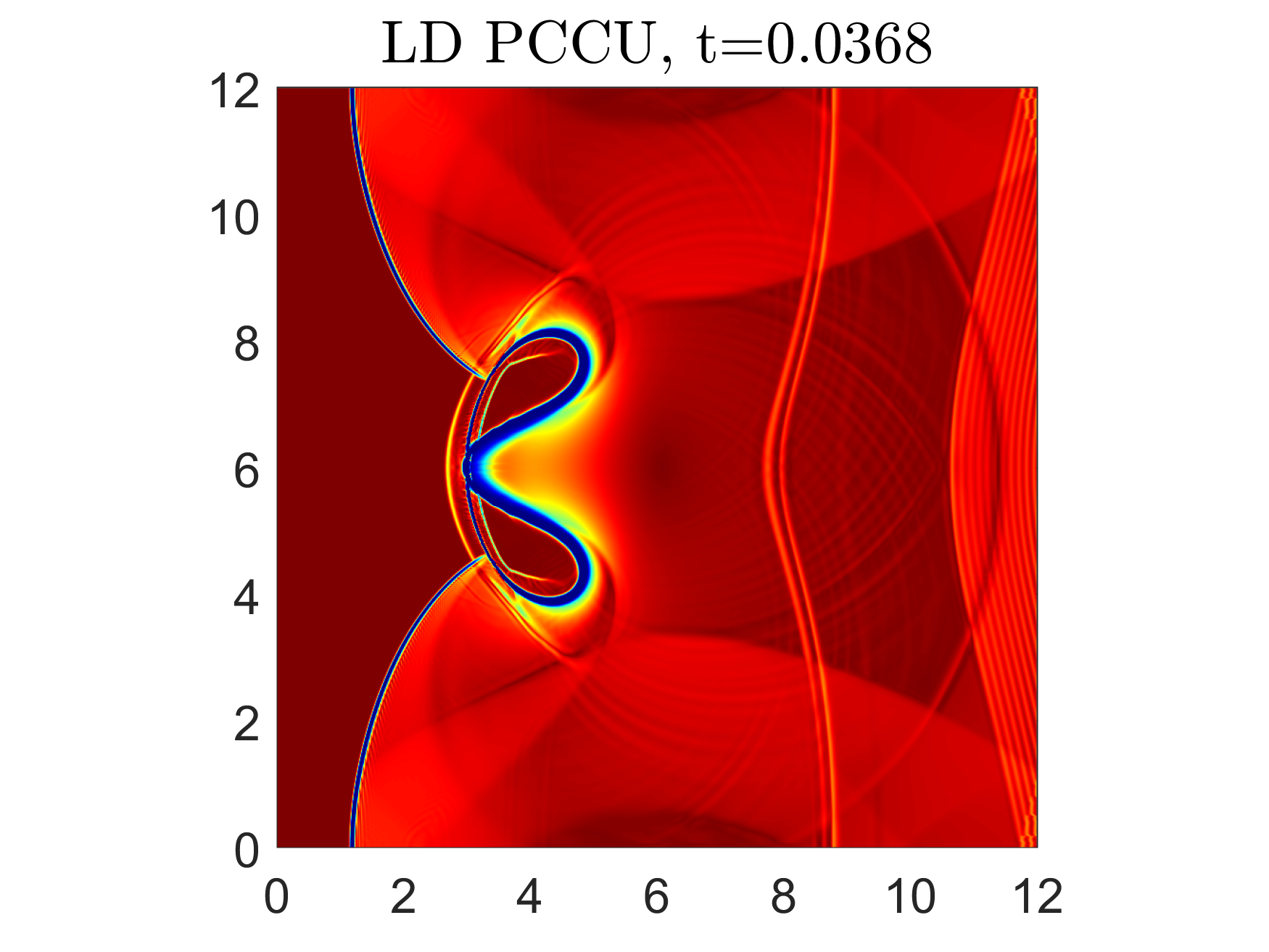}\hspace{0.5cm}
            \includegraphics[trim=2.5cm 0.4cm 2.4cm 0.2cm, clip, width=4.7cm]{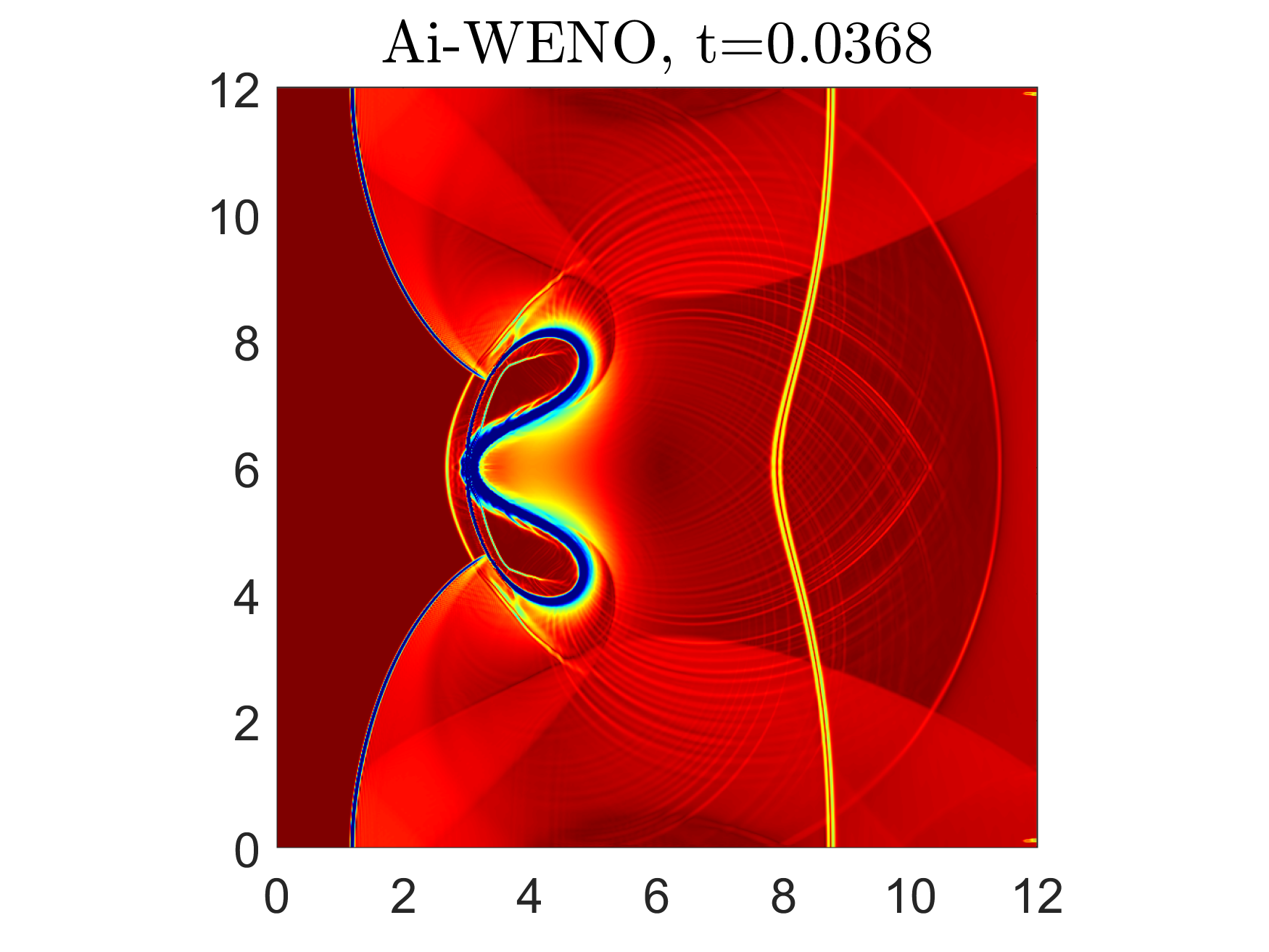}}
\caption{\sf Example 7: Shock-air bubble interaction computed by the PCCU (left column), LD PCCU (middle column), and Ai-WENO (right column)
schemes at times $t=0.0204$, 0.0305, and 0.0368.\label{fig47}}
\end{figure}

\begin{figure}[ht!]
\centerline{\includegraphics[trim=2.5cm 0.4cm 2.4cm 0.2cm, clip, width=4.7cm]{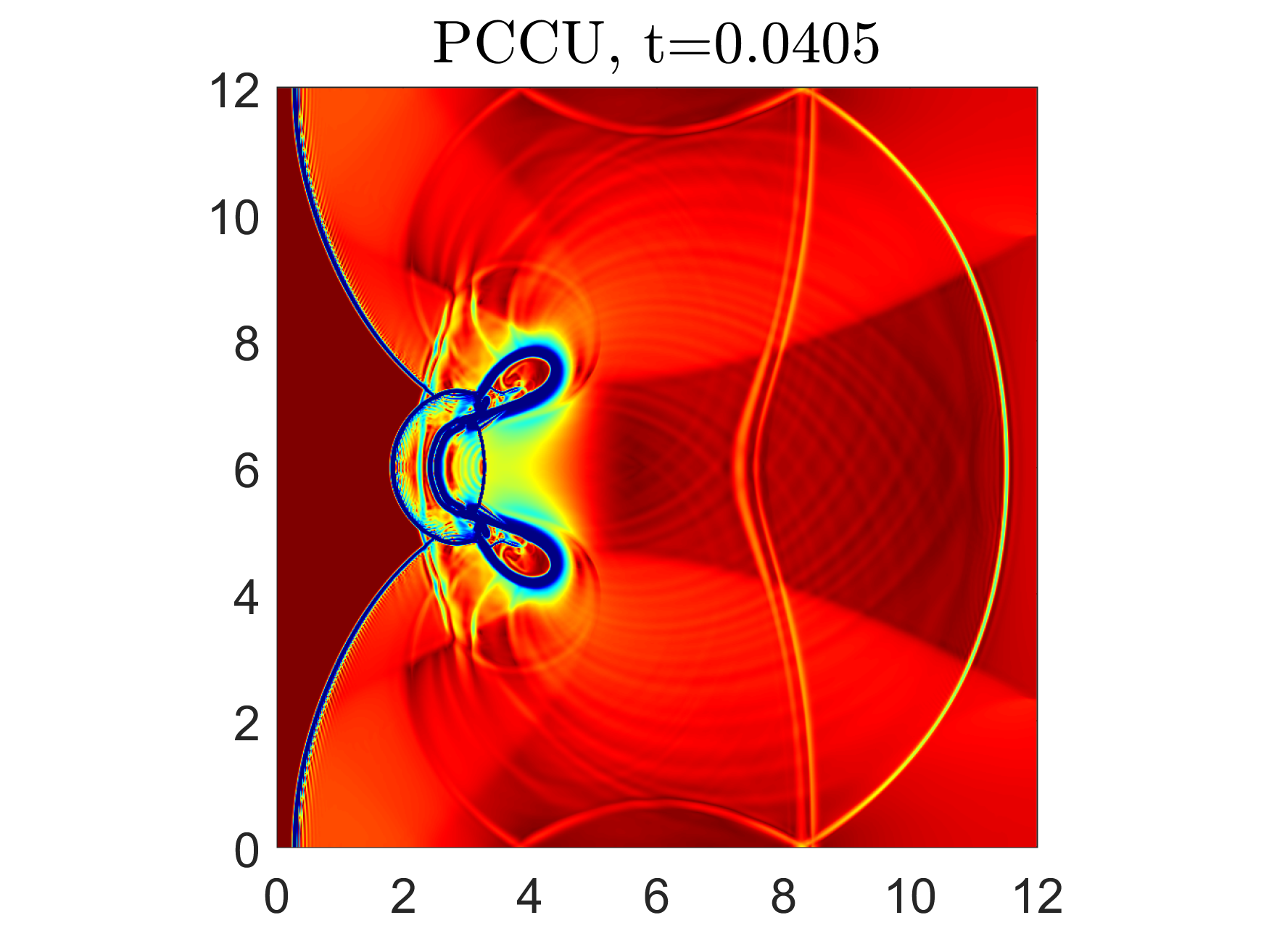}\hspace{0.5cm}
            \includegraphics[trim=2.5cm 0.4cm 2.4cm 0.2cm, clip, width=4.7cm]{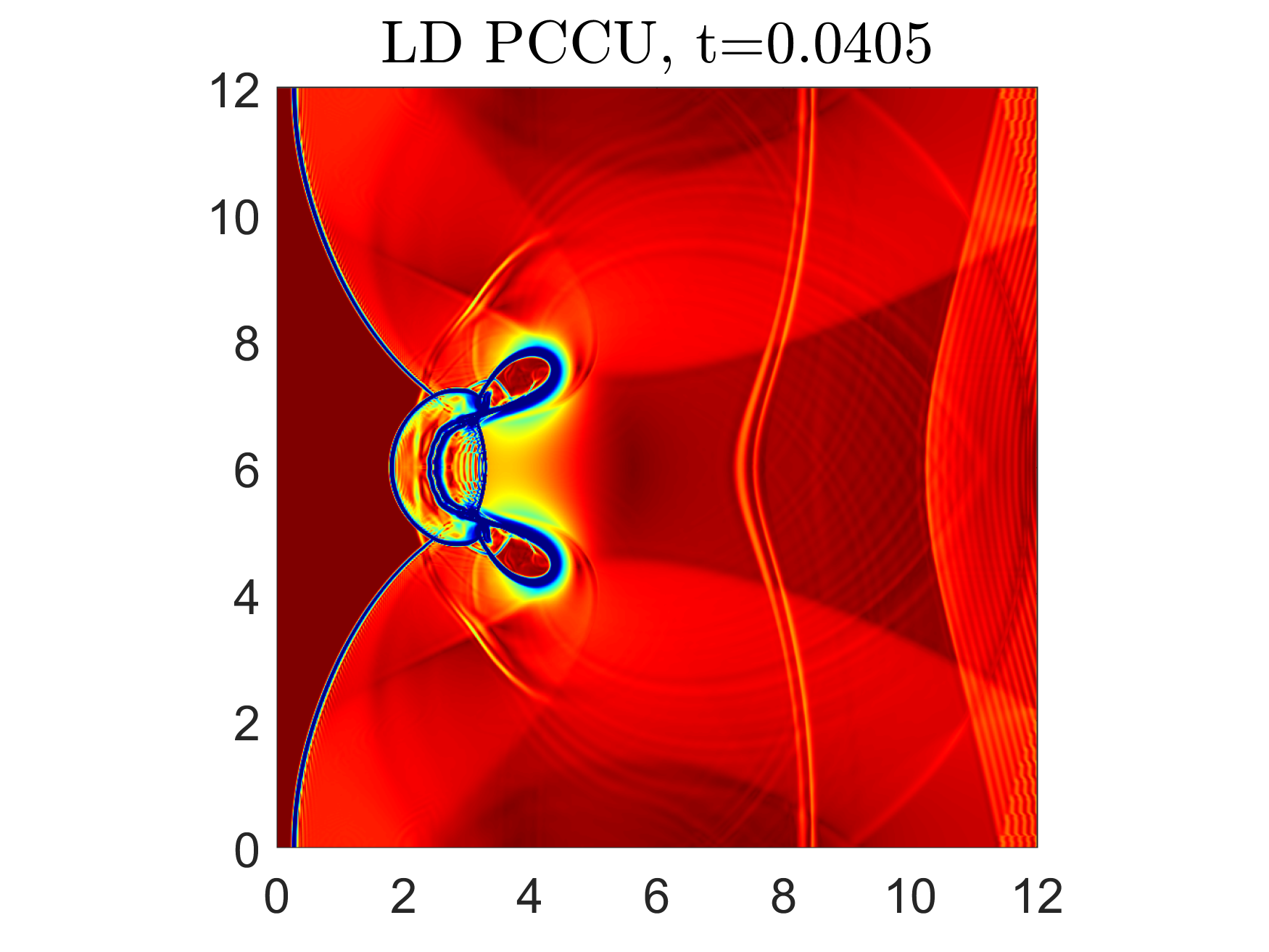}\hspace{0.5cm}
            \includegraphics[trim=2.5cm 0.4cm 2.4cm 0.2cm, clip, width=4.7cm]{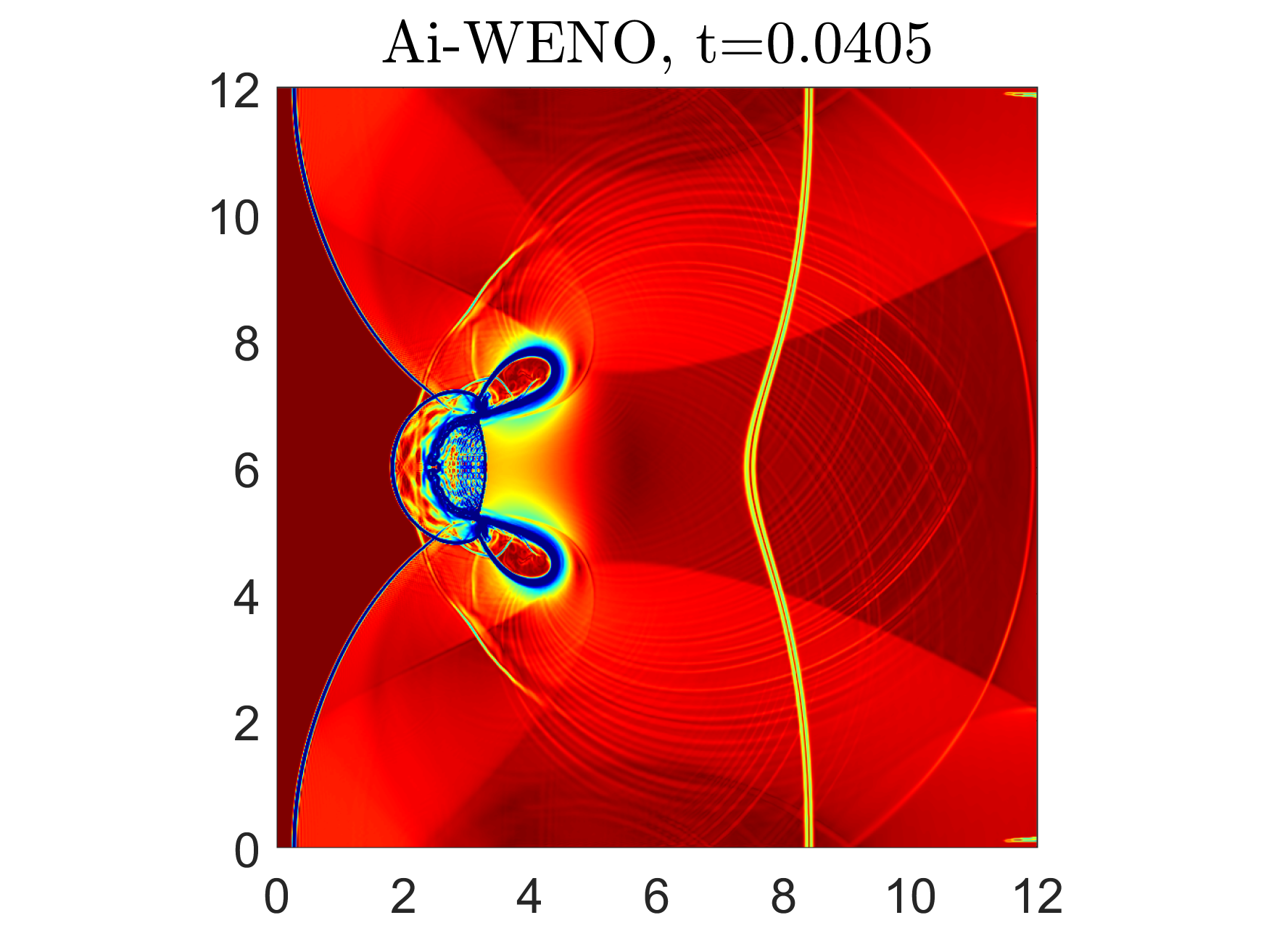}}
\vskip10pt
\centerline{\includegraphics[trim=2.5cm 0.4cm 2.4cm 0.2cm, clip, width=4.7cm]{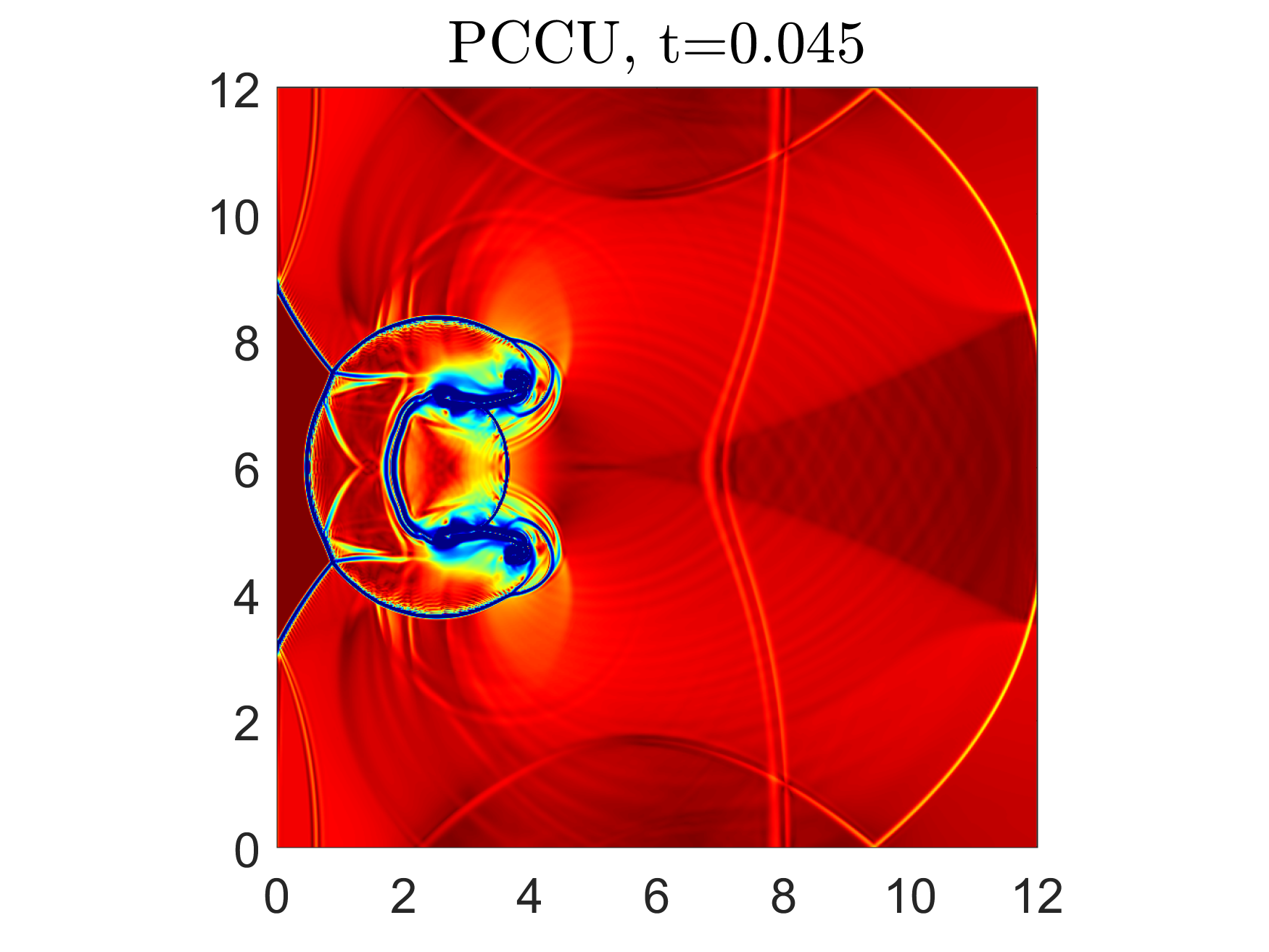}\hspace{0.5cm}
            \includegraphics[trim=2.5cm 0.4cm 2.4cm 0.2cm, clip, width=4.7cm]{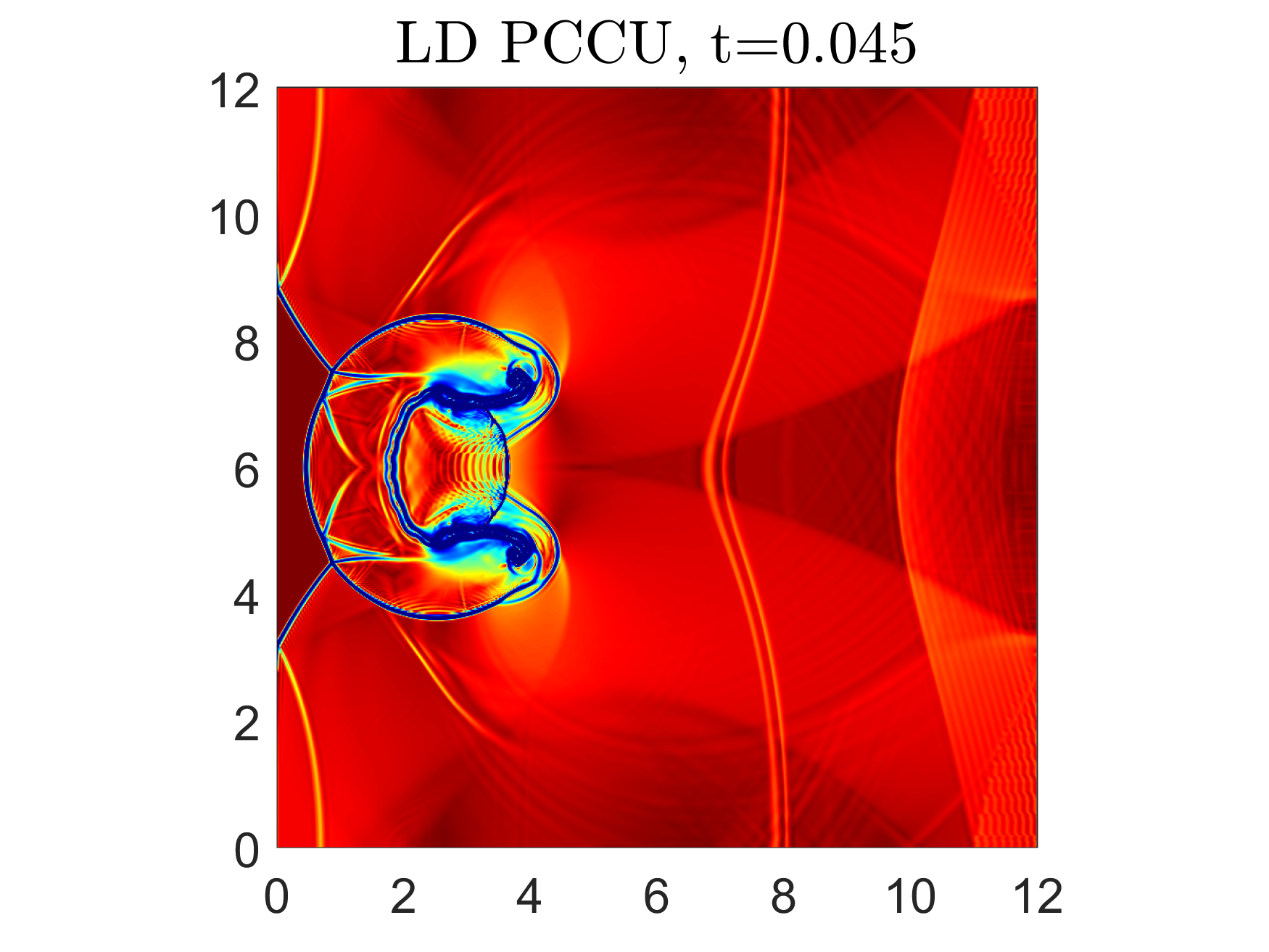}\hspace{0.5cm}
            \includegraphics[trim=2.5cm 0.4cm 2.4cm 0.2cm, clip, width=4.7cm]{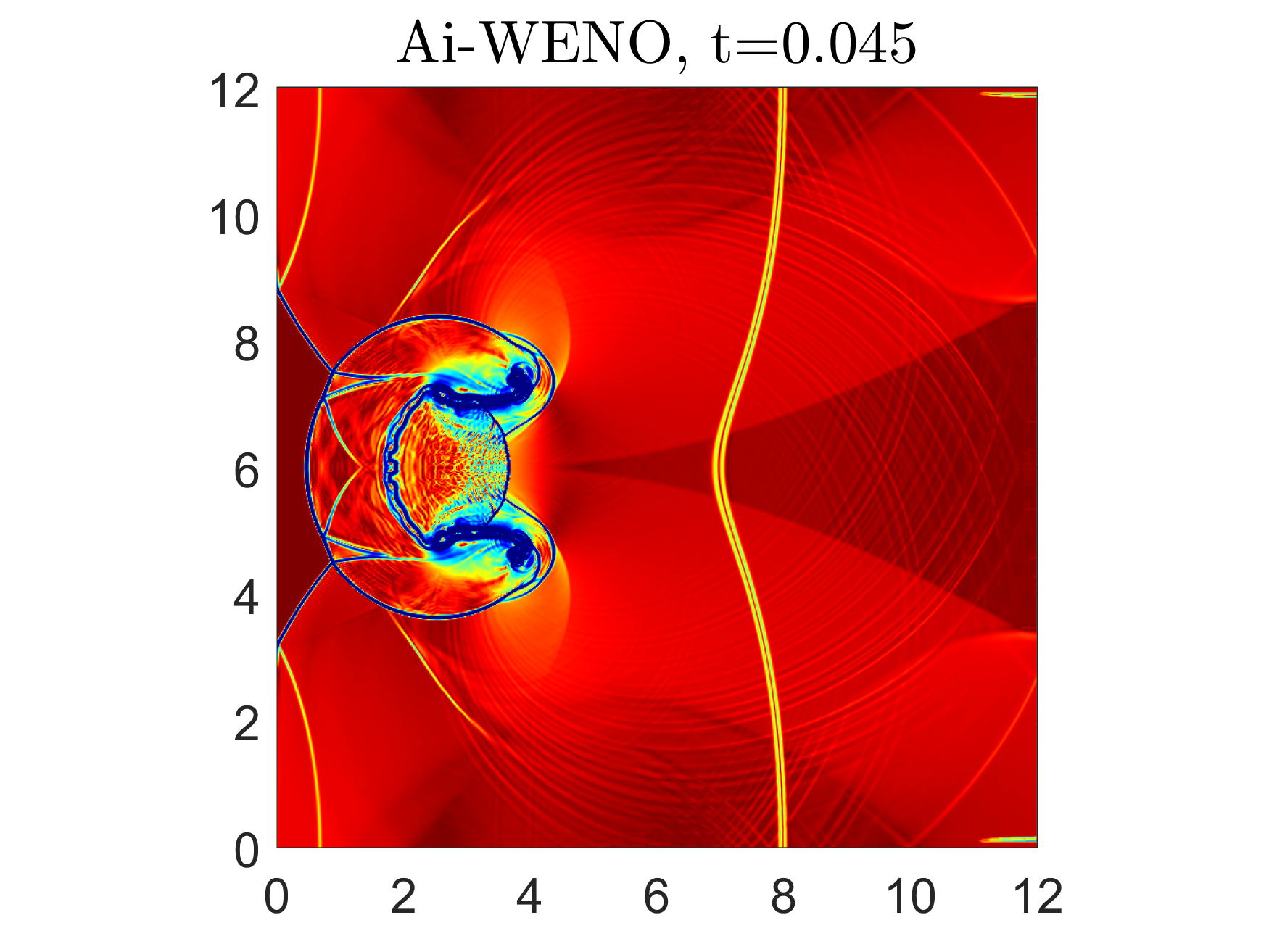}}
\caption{\sf Example 7: Same as in Figure \ref{fig47}, but at larger times $t=0.0405$ and 0.045.\label{fig48}}
\end{figure}

We compute the numerical solutions by the studied PCCU, LD PCCU, and Ai-WENO schemes until the final time $t=0.045$ on a uniform mesh with
$\dx=\dy=3/200$. In Figures \ref{fig47} and \ref{fig48}, we present different stages of the interaction process. As one can see, the bubble
containing air will be compressed by the water, propagates to the left, and changes its shape until losing its integrity and breaking up.
The obtained results are in a good qualitative agreement with the numerical results reported in \cite{CZL20}. As in the previous examples,
the resolution of the bubble interface is significantly improved by the use of the LD PCCU and Ai-WENO schemes, especially for the small
times $t=0.0204$, 0.0305, and 0.0368; see Figure \ref{fig47}. At the same time, the differences near the bubble interfaces between the LD
PCCU and Ai-WENO solutions are minor.
\begin{rmk}
In order to preserve the positivity of the pressure computed by the Ai-WENO scheme, we have used the following algorithm. We first detect
the material interfaces by \eref{3.7a} and \eref{3.7af} and then switch to the second-order LD PCCU scheme (with the one-sided point values
computed by the SBM limiter \eref{2.1.15c}, \eref{3.7}, \eref{3.7b}, \eref{3.7bf} in the overcompressive regime with $\theta=1.3$ and
$\tau=-0.5$) in the neighboring four cells at each side of the interface. This results in a hybrid ``mixed-order'' scheme, which still
achieves higher resolution compared with the second-order LD PCCU scheme.
\end{rmk}

\section{Conclusion}\label{sec5}
In this paper, we have developed flux globalization based low-dissipation (LD) path-conservative central-upwind (PCCU) schemes for one-
(1-D) and two-dimensional (2-D) compressible multifluids. The LD PCCU schemes are based on the flux globalization based PCCU schemes and
employ the recently proposed LDCU fluxes to reduce the numerical dissipation present in the original PCCU schemes. In order to further
enhance the resolution of material interfaces, we track their locations and use the overcompressive SBM limiter in their neighborhoods,
while in the rest of the computational domain, a dissipative generalized minmod limiter is utilized. The new second-order finite-volume
method is then extended to the fifth order of accuracy via the finite-difference A-WENO framework. We have applied the developed schemes to
a number of 1-D and 2-D examples and the obtained numerical results clearly demonstrate that both of the LD {\color{red}PCCU} and LD A-WENO schemes
outperform the flux globalization based PCCU scheme that employs the original central-upwind numerical flux from \cite{KLin}. At the same
time, these examples show that the fifth-order LD Ai-WENO scheme enhances the resolution achieved by the second-order LD PCCU scheme.

\subsection*{Acknowledgment}
The work of A. Kurganov was supported in part by NSFC grant 12171226 and by the fund of the Guangdong Provincial Key
Laboratory of Computational Science and Material Design (No. 2019B030301001).

\appendix
\section{1-D Fifth-Order Ai-WENO-Z Interpolation}\label{appa}
In this appendix, we briefly describe the fifth-order Ai-WENO-Z interpolation recently introduced in \cite{DLWW22,LLWDW23,WD22}.

Assume that the point values $W_{j+\ell}$ of a certain quantity $W$ at the uniform grid points $x=x_{j+\ell}$, $\ell=-2,\ldots,3$ are
available. We now show how to obtain an interpolated left-sided value of $W$ at $x=x_\jph$, denoted by $W^-_\jph$. The right-sided value
$W^+_\jph$ can then be obtained in the mirror-symmetric way.

The value $W^-_\jph$ is calculated using a weighted average of the three parabolic interpolants ${\cal P}_0(x)$, ${\cal P}_1(x)$, and
${\cal P}_2(x)$ obtained using the stencils $[x_{j-2},x_{j-1},x_j]$, $[x_{j-1},x_j,x_{j+1}]$, and $[x_j,x_{j+1},x_{j+2}]$, respectively:
\begin{equation*}
W^-_\jph=\sum_{k=0}^2\omega_k{\cal P}_k\big(x_\jph\big),
\end{equation*}
where
$$
\begin{aligned}
&{{\cal P}}_0(x_\jph)=\frac{3}{8}\,W_{j-2}-\frac{5}{4}\,W_{j-1}+\frac{15}{8}\,W_j,\quad
{{\cal P}}_1(x_\jph)=-\frac{1}{8}\,W_{j-1}+\frac{3}{4}\,W_j+\frac{3}{8}\,W_{j+1},\\
&{{\cal P}}_2(x_\jph)=\frac{3}{8}\,W_j+\frac{3}{4}\,W_{j+1}-\frac{1}{8}\,W_{j+2},
\end{aligned}
$$
and the Ai-weights $\omega_k$ are computed by
\begin{equation}
\omega_k=\frac{\alpha_k}{\alpha_0+\alpha_1+\alpha_2},\quad
\alpha_k=d_k\left[1+\bigg(\frac{\tau_5}{\beta_k+\varepsilon\mu^2_j}\bigg)^{\!r\,}\right],\quad k=0,1,2,
\label{A.1}
\end{equation}
with $d_0=\frac{1}{16}$, $d_1=\frac{5}{8}$, and $d_2=\frac{5}{16}$. The smoothness indicators $\beta_k$ for the corresponding parabolic
interpolants ${\cal P}_k(x)$ are given by
\begin{equation*}
\begin{aligned}
&\beta_0=\frac{13}{12}\big(W_{j-2}-2W_{j-1}+W_j\big)^2+\frac{1}{4}\big(W_{j-2}-4W_{j-1}+3W_j\big)^2,\\
&\beta_1=\frac{13}{12}\big(W_{j-1}-2W_j+W_{j+1}\big)^2+\frac{1}{4}\big(W_{j-1}-W_{j+1}\big)^2,\\
&\beta_2=\frac{13}{12}\big(W_j-2W_{j+1}+W_{j+2}\big)^2+\frac{1}{4}\big(3W_j-4W_{j+1}+W_{j+2}\big)^2.
\end{aligned}
\end{equation*}
Finally, in formula \eref{A.1}, $\tau_5=|\beta_2-\beta_0|$, $\mu_j=\frac{1}{5}\sum^{j+2}_{\ell=j-2}|W_\ell-\widehat W_j|+10^{-40}$ with
$\widehat W_j:=\frac{1}{5}\sum^{j+2}_{\ell=j-2}W_\ell$, and in all of the numerical examples, we have chosen $r=2$ and
$\varepsilon=10^{-12}$.

\subsection{1-D Local Characteristic Decomposition}\label{appb}
In \S\ref{sec2.3}, the Ai-WENO-Z interpolant is applied to the local characteristic variables, which are obtained using the LCD. To this
end, we first rewrite the studied $\gamma$-based multifluid system in terms of the primitive variables $\bm V=(\rho,u,p,\Gamma,\Pi)^\top$:
$$
\begin{aligned}
\bm V_t+{\cal A}\bm V_x=\bm0,\quad{\cal A}:=\begin{pmatrix}u&\rho&0&0&0\\0&u&\dfrac{1}{\rho}&0&0\\0&\gamma(p+\pi_\infty)&u&0&0\\
0&0&0&u&0\\0&0&0&0&u\end{pmatrix},
\end{aligned}
$$
and introduce the locally averaged matrices
\begin{equation}
\begin{aligned}
\widehat{\cal A}_\jph:=\begin{pmatrix}\hat u&\hat\rho&0&0&0\\0&\hat u&\dfrac{1}{\hat\rho}&0&0\\
0&\hat\gamma(\hat p+\hat\pi_\infty)&\hat u&0&0\\0&0&0&\hat u&0\\0&0&0&0&\hat u\end{pmatrix},
\end{aligned}
\label{A.1.1}
\end{equation}
where $\hat{(\cdot)}$ stands for the following averages (see \cite{Karni93}):
\begin{equation*}
\begin{aligned}
&\hat\rho=\sqrt{\rho_j\rho_{j+1}},\quad\hat u=\frac{\sqrt{\rho_j}u_j+\sqrt{\rho_{j+1}}u_{j+1}}{\sqrt{\rho_j}+\sqrt{\rho_{j+1}}},\quad
\hat p=\frac{\sqrt{\rho_j}p_j+\sqrt{\rho_{j+1}}p_{j+1}}{\sqrt{\rho_j}+\sqrt{\rho_{j+1}}},\\
&\hat\gamma=\frac{\sqrt{\rho_j}\gamma_j+\sqrt{\rho_{j+1}}\gamma_{j+1}}{\sqrt{\rho_j}+\sqrt{\rho_{j+1}}},\quad
\hat\pi_\infty=\frac{\sqrt{\rho_j}(\pi_\infty)_j+\sqrt{\rho_{j+1}}(\pi_\infty)_{j+1}}{\sqrt{\rho_j}+\sqrt{\rho_{j+1}}},
\end{aligned}
\end{equation*}
where $\gamma_j=1+1/\Gamma_j$ and $(\pi_\infty)_j={\Pi_j}/(1+\Gamma_j)$.

We then compute the matrix $R_\jph$ composed of the right eigenvectors of $\widehat{\cal A}_\jph$ and obtain
\begin{equation}
R_\jph=\begin{pmatrix}\dfrac{1}{\hat c^2}&0&0&1&\dfrac{1}{\hat c^2}\\[1.2ex]
-\dfrac{1}{\hat\rho\hat c}&0&0&0&\dfrac{1}{\hat\rho\hat c}\\1&0&0&0&1\\0&0&1&0&0\\0&1&0&0&0\end{pmatrix}\quad\mbox{and}\quad
R^{-1}_\jph=\begin{pmatrix}0&-\dfrac{\hat\rho\hat c}{2}&\dfrac{1}{2}&0&0\\0&0&0&0&1\\0&0&0&1&0\\1&0&-\dfrac{1}{\hat c^2}&0&0\\[1.2ex]
0&\dfrac{\hat\rho\hat c}{2}&\dfrac{1}{2}&0&0\end{pmatrix},
\label{A.1.3}
\end{equation}
where $\hat c=\sqrt{\gamma(\hat p+\hat\pi_\infty)/\hat\rho}$. Notice that all of the $\hat{(\cdot)}$ quantities in
\eref{A.1.1}--\eref{A.1.3} have to have a subscript index, that is, $\hat{(\cdot)}=\hat{(\cdot)}_\jph$, but we have omitted it for the sake
of brevity for all of the quantities except for $\widehat{\cal A}_\jph$.

Finally, we introduce the local characteristic variables in the neighborhood of $x=x_\jph$:
$$
\bm W_{j+\ell}=R^{-1}_\jph\mV_{j+\ell},\quad\ell=-2,\ldots,3,
$$
and apply the Ai-WENO-Z interpolation to every component of $\bm W$ to obtain $\bm W^\pm_\jph$, and then we end up with
\begin{equation*}
\mV^\pm_\jph=R_\jph\bm W^\pm_\jph.
\end{equation*}

\section{2-D Local Characteristic Decomposition}\label{appc}
In this appendix, we extend the 1-D LCD described in Appendix \ref{appb} to the system \eref{3.17}, which can be rewritten in terms of the
primitive variables $\bm V=(\rho,u,v,p,\Gamma,\Pi)^\top$ as
$$
\begin{aligned}
\bm V_t+{\cal A}\bm V_x=\bm0,\quad{\cal A}:=\begin{pmatrix}u&\rho&0&0&0&0\\0&u&0&\dfrac{1}{\rho}&0&0\\0&0&u&0&0&0\\
0&\gamma(p+\pi_\infty)&0&u&0&0\\0&0&0&0&u&0\\0&0&0&0&0&u\end{pmatrix},
\end{aligned}
$$
and introduce the locally averaged matrices
$$
\begin{aligned}
\widehat{\cal A}_{\jph,k}=\begin{pmatrix}\hat u&\hat\rho&0&0&0&0\\0&\hat u&0&\dfrac{1}{\hat\rho}&0&0\\0&0&\hat u&0&0&0\\
0&\hat\gamma(\hat p+\hat\pi_\infty)&0&\hat u&0&0\\0&0&0&0&\hat u&0\\0&0&0&0&0& \hat u\end{pmatrix},
\end{aligned}
$$
where $\hat{(\cdot)}$ stands for the following averages:
\begin{equation*}
\begin{aligned}
&\hat\rho=\sqrt{\rho_{j,k}\,\rho_{j+1,k}},\quad\hat u=\frac{\sqrt{\rho_{j,k}}u_{j,k}+\sqrt{\rho_{j+1,k}}u_{j+1,k}}
{\sqrt{\rho_{j,k}}+\sqrt{\rho_{j+1,k}}},\quad\hat p=\frac{\sqrt{\rho_{j,k}}p_{j,k}+\sqrt{\rho_{j+1,k}}p_{j+1,k}}
{\sqrt{\rho_{j,k}}+\sqrt{\rho_{j+1,k}}},\\
&\hat\gamma=\frac{\sqrt{\rho_{j,k}}\,\gamma_{j,k}+\sqrt{\rho_{j+1,k}}\,\gamma_{j+1,k}}{\sqrt{\rho_{j,k}}+\sqrt{\rho_{j+1,k}}},\quad
\hat\pi_{\infty}=\frac{\sqrt{\rho_{j,k}}(\pi_\infty)_{j,k}+\sqrt{\rho_{j+1,k}}(\pi_\infty)_{j+1,k}}{\sqrt{\rho_{j,k}}+\sqrt{\rho_{j+1,k}}},
\end{aligned}
\end{equation*}
with $\gamma_{j,k}=1+1/{\Gamma_{j,k}}$ and $(\pi_\infty)_{j,k}=\Pi_{j,k}/(1+\Gamma_{j,k})$.

We then compute the matrices $R_{\jph,k}$ and $R^{-1}_{\jph,k}$ such that the matrix $R^{-1}_{\jph,k}\widehat{\cal A}_{\jph,k}R_{\jph,k}$ is
diagonal and obtain
\begin{equation*}
R_{\jph,k}=\begin{pmatrix}\dfrac{1}{\hat c^2}&0&0&0&1&\dfrac{1}{\hat c^2}\\[1.2ex]
-\dfrac{1}{\hat\rho\hat c}&0&0&0&0&\dfrac{1}{\hat\rho\hat c}\\0&0&0&1&0&0\\1&0&0&0&0&1\\0&0&1&0&0&0\\0&1&0&0&0&0\end{pmatrix}\quad\mbox{and}
\quad R^{-1}_{\jph,k}=\begin{pmatrix}0&-\dfrac{\hat\rho\hat c}{2}&0&\dfrac{1}{2}&0&0\\0&0&0&0&0&1\\0&0&0&0&1&0\\0&0&1&0&0&0\\
1&0&0&-\dfrac{1}{\hat c^2}&0&0\\[1.2ex]0&\dfrac{\hat\rho\hat c}{2}&0&\dfrac{1}{2}&0&0\end{pmatrix}.
\end{equation*}
Appendix \ref{appb}, we have omitted the $(\jph,k)$ indices for all of the $\hat{(\cdot)}$ quantities except for
$\widehat{\cal A}_{\jph,k}$.

Finally, given the matrices $R^{-1}_{\jph,k}$ and $R_{\jph,k}$, we introduce the local characteristic variablesnin the neighborhood of
$(x,y)=(x_\jph, y_k)$:
$$
\bm W_{j+\ell,k}=R^{-1}_{\jph,k}\mV_{j+\ell,k},\quad\ell=-2,\ldots,3,
$$
apply the Ai-WENO-Z interpolation to every component of $\bm W$ to obtain $\bm W^\pm_{\jph,k}$, and end up with
\begin{equation*}
\mV^\pm_{\jph,k}=R_{\jph,k}\bm W^\pm_{\jph,k}.
\end{equation*}

Notice that the point values $\mV^\pm_{j,\kph}$ are obtained in a similar manner and we omit the details for the sake of brevity.

\bibliography{reference}
\bibliographystyle{siam}
\end{document}